\documentclass
[
    fontsize = 11 pt,        
    american,                
    captions = tableheading, 
    numbers = noenddot,      
    footheight = 35 pt,      
]
{scrartcl}
\usepackage[left=2cm,right=2cm, top=3cm, bottom=3.75cm]{geometry}

\usepackage
[
    babel = true, 
]
{microtype}           
\usepackage{csquotes} 

\usepackage
[
    dvipsnames,
    table,
]{xcolor} 

\definecolor{HPIblue}{rgb}{0,0.25,0.5} 
\definecolor{HPIorange}{rgb}{1,0.5,0} 

\usepackage{amsmath}
\usepackage{amssymb}
\usepackage{amsthm}
\usepackage{thmtools}
\usepackage{mathtools}
\usepackage{thm-restate}
\usepackage{dsfont}        

\usepackage
[
    ttscale = 0.85, 
]
{libertine} 
\usepackage
[
    libertine,    
    slantedGreek, 
    vvarbb,       
    libaltvw,     
]
{newtxmath} 
\usepackage{url} 
\usepackage{bm}  

\usepackage{graphicx} 
\usepackage
[
    subrefformat = simple, 
    labelformat = simple,  
]
{subcaption}         

\usepackage{array}     
\usepackage{booktabs}  
\usepackage{longtable} 
\usepackage{pdflscape} 

\usepackage[shortlabels]{enumitem} 

\usepackage
[
    ngerman,
    ruled,         
    vlined,        
    linesnumbered, 
    algo2e
]
{algorithm2e} 
\DontPrintSemicolon
\SetAlgoCaptionLayout{leftline}
\LinesNumbered
\let\oldnl\nl
\newcommand{\nonl}{\renewcommand{\nl}{\let\nl\oldnl}}

\usepackage{xspace}   
\usepackage
[
    shortcuts, 
]
{extdash}             
\usepackage{setspace} 

\usepackage{xparse}    
\usepackage{footnote}  
\usepackage{afterpage} 
\usepackage
[
    textsize = scriptsize, 
]
{todonotes}            

\definecolor{darkblue}{rgb}{0, 0, 0.5}

\usepackage
[
    bookmarks = true,                 
    bookmarksopen = false,            
    bookmarksnumbered = true,         
    pdfstartpage = 1,                 
    colorlinks = true,                
    allcolors = HPIblue,              
    breaklinks = true,
    linkcolor = darkblue,
    anchorcolor = darkblue,
    citecolor = darkblue,
    filecolor = darkblue,
    menucolor = darkblue,
    urlcolor = darkblue
]
{hyperref} 

\usepackage
[
    noabbrev,   
    nameinlink, 
]
{cleveref} 

\usepackage{nicefrac}

\usepackage{scrlayer-scrpage}
\usepackage{lastpage}

\usepackage{mdframed}

\usepackage[normalem]{ulem}

\usepackage{listings}
\usepackage[theorems,breakable]{tcolorbox}

\newcommand*{\ind}[1]{\mathOrText{\mathbf{1}{\{#1\}}}}

\newcommand{\referDefined}[1]{\textcolor{HPIblue}{#1}}

\newcommand{\natnum}{\mathbb{N}}

\newcommand{\integers}{\mathbb{Z}}
\newcommand{\realnum}{\mathbb{R}}

\newcommand{\CalX}{\mathcal{X}}

\newcommand{\eqnComment}[2]{\underset{{\scriptstyle \text{#1}}}{#2}}
\newcommand{\set}[2]{\{ #1 \mid #2 \}}

\newcommand{\bigO}[1]{\mathcal{O}\mathclose{\left(#1\right)}}

\newcommand{\lt}{<}
\newcommand{\gt}{>}

\renewcommand{\Pr}[1]{\operatorname{Pr}\left[#1\right]}
\newcommand{\Ew}[1]{\operatorname{E}\left[#1\right]}
\newcommand{\Var}[1]{\operatorname{Var}\left[#1\right]}

\newcommand*{\printAuthor}{}
\newcommand*{\authorName}[1]{\renewcommand*{\printAuthor}{#1}}

\newcommand*{\printTitle}{}
\newcommand*{\titleName}[1]{\renewcommand*{\printTitle}{#1}}

\newcommand*{\printSubTitle}{}
\newcommand*{\subTitleName}[1]{\renewcommand*{\printSubTitle}{#1}}

\newcommand*{\printmakerVersion}{}
\newcommand*{\makerVersionName}[1]{\renewcommand*{\printmakerVersion}{#1}}

\date{}

\makeatletter
    \def\IfEmptyTF#1%
    {%
        \if\relax\detokenize{#1}\relax%
            \expandafter\@firstoftwo%
        \else%
            \expandafter\@secondoftwo%
        \fi%
    }
\makeatother

\NewDocumentCommand{\mathOrText}{m}
{%
    \ensuremath{#1}\xspace%
}

\let\originalleft\left
\let\originalright\right
\renewcommand{\left}{\mathopen{}\mathclose\bgroup\originalleft}
\renewcommand{\right}{\aftergroup\egroup\originalright}

\makeatletter
    \DeclareRobustCommand{\bfseries}%
    {%
        \not@math@alphabet\bfseries\mathbf%
        \fontseries\bfdefault\selectfont%
        \boldmath%
    }

\xspaceaddexceptions{]\}}

\allowdisplaybreaks

\crefname{ineq}{inequality}{inequalities}
\creflabelformat{ineq}{#2{\upshape(#1)}#3}

\crefname{term}{term}{terms}
\creflabelformat{term}{#2{\upshape(#1)}#3}

\crefname{cond}{condition}{conditions}
\creflabelformat{term}{#2{\upshape(#1)}#3}

\setfootnoterule{0 cm}

\deffootnote[1.2 em]{1.2 em}{0 em}{\makebox[1.4 em][l]{\textbf{\thefootnotemark}}}

\let\oldfootnote\footnote

\newlength{\spaceBeforeFootnote} 
\newlength{\spaceAfterFootnote}  

\RenewDocumentCommand{\footnote}{o o o m}%
{%
    \IfNoValueTF{#1}%
    {%
        \oldfootnote{#4}%
    }%
    {%
        \setlength{\spaceBeforeFootnote}{\IfEmptyTF{#1}{0}{#1} em}%
        \IfNoValueTF{#2}%
        {%
            \hspace*{\spaceBeforeFootnote}\oldfootnote{#4}%
        }%
        {%
            \setlength{\spaceAfterFootnote}{\IfEmptyTF{#2}{0}{#2} em}%
            \hspace*{\spaceBeforeFootnote}\IfNoValueTF{#3}{\oldfootnote{#4}}{\oldfootnote[#3]{#4}}\hspace*{\spaceAfterFootnote}%
        }%
    }%
}

\makesavenoteenv{figure}
\makesavenoteenv{table}
\makesavenoteenv{tabular}

\lohead[\normalfont\printTitle]{}
\rohead[\normalfont\printAuthor]{}
\cohead{}
\cofoot*{}
\cofoot*{\normalfont\thepage{}~/~\pageref{LastPage}}

\newcommand{\makeHeader}{
\thispagestyle{scrheadings}

\includegraphics[height = 22 mm]{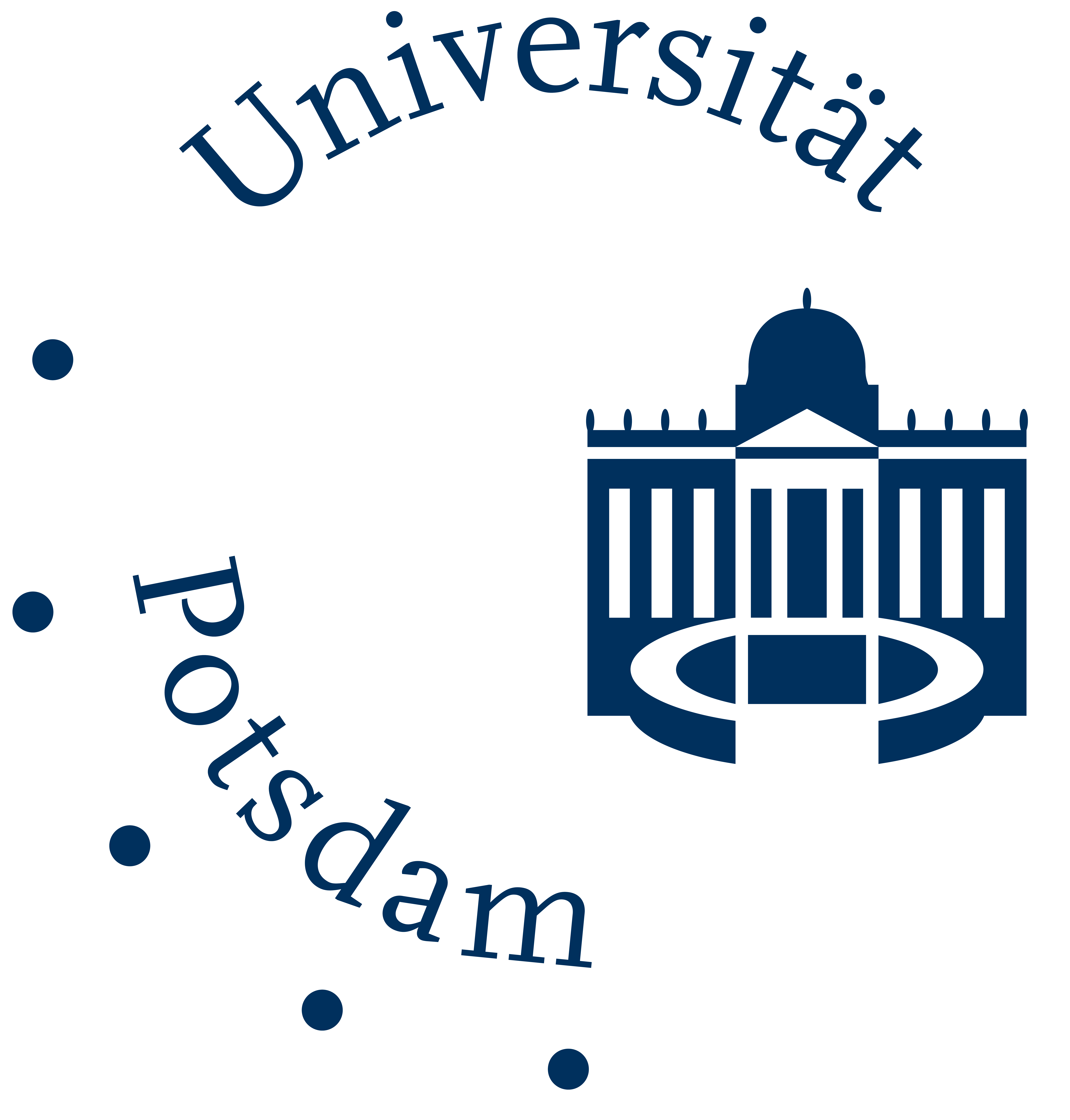}
\hfill\includegraphics[height = 1.9 cm]{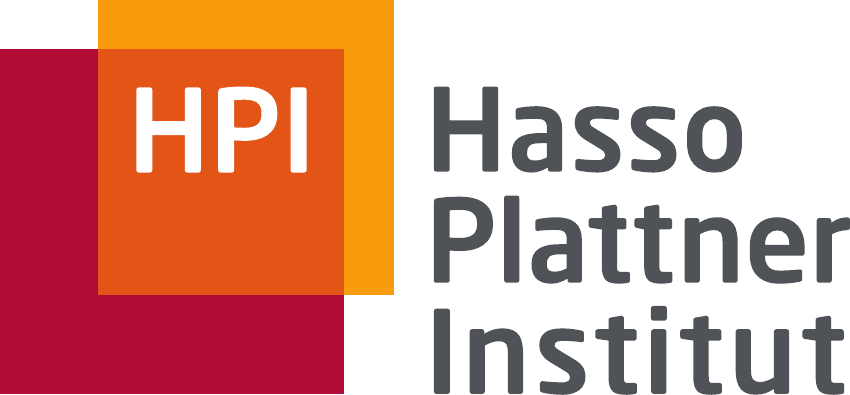}

\vfill

\begin{center}
{\Huge \printTitle}\\[10mm]
{\Large \printSubTitle\\[30mm]
\printAuthor}\\[50mm]
\end{center}

\vfill

\noindent%
\rule{\textwidth}{1 pt}\vspace*{1 ex}

\noindent%
\begin{center}
Compiled on \today\ 
\end{center}

\noindent%
\rule{\textwidth}{1 pt}

\clearpage}

\deffootnote[1.2 em]{1.2 em}{0 em}{\makebox[1.4 em][l]{\textbf{\thefootnotemark}}}

\newcommand{\ignore}[1]{}

\newcommand{\inlineComment}[1]{\textcolor{darkblue}{[Kommentar: #1]}}

\newcommand{\buildRef}[1]{\Cref{#1} \setbox0=\hbox{\nameref{#1}\unskip}\ifdim\wd0=0pt
\else
 [\nameref{#1}]%
\fi}

\newtcbtheorem[auto counter,number within=section,Crefname={Satz}{Sätze}]{theoremDE}%
  {Satz}{fonttitle=\bfseries\upshape, fontupper=\slshape\upshape,
     arc=1mm, colback=blue!5!white,colframe=HPIorange,breakable=true}{}
\newtcbtheorem[auto counter,number within=section,Crefname={Theorem}{Theorems}]{theoremEN}%
  {Theorem}{fonttitle=\bfseries\upshape, fontupper=\slshape\upshape,
     arc=1mm, colback=blue!5!white,colframe=HPIorange,breakable=true}{}

\newtcbtheorem[use counter from=theoremDE,number within=section,Crefname={Proposition}{Propositions}]{propositionDE}%
  {Propostion}{fonttitle=\bfseries\upshape, fontupper=\slshape\upshape,
     arc=1mm, colback=blue!5!white,colframe=HPIorange,breakable=true}{}
\newtcbtheorem[use counter from=theoremEN,number within=section,Crefname={Proposition}{Propositionen}]{propositionEN}%
  {Propostion}{fonttitle=\bfseries\upshape, fontupper=\slshape\upshape,
     arc=1mm, colback=blue!5!white,colframe=HPIorange,breakable=true}{}

\newtcbtheorem[use counter from=theoremDE,number within=section,Crefname={Korollar}{Korollare}]{corollaryDE}%
  {Korollar}{fonttitle=\bfseries\upshape, fontupper=\slshape\upshape,
     arc=1mm, colback=blue!5!white,colframe=HPIorange,breakable=true}{}		
\newtcbtheorem[use counter from=theoremEN,number within=section,Crefname={Corollary}{Corollaries}]{corollaryEN}%
  {Corollary}{fonttitle=\bfseries\upshape, fontupper=\slshape\upshape,
     arc=1mm, colback=blue!5!white,colframe=HPIorange,breakable=true}{}
			
\newtcbtheorem[use counter from=theoremDE,Crefname={Definition}{Definitionen}]{definitionDE}%
  {Definition}{fonttitle=\bfseries\upshape, fontupper=\slshape\upshape,
     arc=1mm, colback=blue!5!white,colframe=HPIblue,breakable=true}{}		
\newtcbtheorem[use counter from=theoremEN,Crefname={Definition}{Definitions}]{definitionEN}%
  {Definition}{fonttitle=\bfseries\upshape, fontupper=\slshape\upshape,
     arc=1mm, colback=blue!5!white,colframe=HPIblue,breakable=true}{}

\newtcbtheorem[use counter from=theoremDE,Crefname={Beispiel}{Beispiele}]{beispielDE}%
  {Beispiel}{fonttitle=\bfseries\upshape, fontupper=\slshape\upshape,coltitle=black,
    arc=0mm, colback=blue!5!white,colframe=HPIblue,detach title,boxrule=0pt,leftrule=3pt,
		attach title to upper={\ --- \ },breakable=true}{}		
\newtcbtheorem[use counter from=theoremEN,Crefname={Example}{Examples}]{beispielEN}%
  {Example}{fonttitle=\bfseries\upshape, fontupper=\slshape\upshape,coltitle=black,
    arc=0mm, colback=blue!5!white,colframe=HPIblue,detach title,boxrule=0pt,leftrule=3pt,
		attach title to upper={\ --- \ },breakable=true}{}

\newtcbtheorem[number within=section,Crefname={Bemerkung}{Bemerkungen}]{noteDE}%
  {}{theorem name,fonttitle=\bfseries\upshape, fontupper=\slshape\upshape,coltitle=black,
    arc=0mm, colback=blue!5!white,colframe=HPIblue,detach title,boxrule=0pt,leftrule=3pt,breakable=true}{}	
\newtcbtheorem[number within=section,Crefname={Note}{Notes}]{noteEN}%
  {}{theorem name,fonttitle=\bfseries\upshape, fontupper=\slshape\upshape,coltitle=black,
    arc=0mm, colback=blue!5!white,colframe=HPIblue,detach title,boxrule=0pt,leftrule=3pt,breakable=true}{}	
		
\newtcbtheorem[use counter from=theoremDE,number within=section,Crefname={Konvention}{Konventionen}]{conventionDE}%
  {Konvention}{fonttitle=\bfseries\upshape, fontupper=\slshape\upshape,coltitle=black,
    arc=0mm, colback=blue!5!white,colframe=HPIblue,detach title,boxrule=0pt,leftrule=3pt,
		attach title to upper={\ --- \ },breakable=true}{}	
\newtcbtheorem[use counter from=theoremEN,number within=section,Crefname={Convention}{Conventions}]{conventionEN}%
  {Convention}{fonttitle=\bfseries\upshape, fontupper=\slshape\upshape,coltitle=black,
    arc=0mm, colback=blue!5!white,colframe=HPIblue,detach title,boxrule=0pt,leftrule=3pt,
		attach title to upper={\ --- \ },breakable=true}{}

\newtcbtheorem[use counter from=theoremDE,Crefname={Lemma}{Lemmas}]{lemmaDE}%
  {Lemma}{fonttitle=\bfseries\upshape, fontupper=\slshape\upshape,
     arc=1mm, colback=blue!5!white,colframe=HPIblue,breakable=true}{}		
\newtcbtheorem[use counter from=theoremEN,Crefname={Lemma}{Lemmas}]{lemmaEN}%
  {Lemma}{fonttitle=\bfseries\upshape, fontupper=\slshape\upshape,
     arc=1mm, colback=blue!5!white,colframe=HPIblue,breakable=true}{}

\usetikzlibrary {arrows.meta}

\usepackage[USenglish]{babel}

\authorName{Timo Kötzing}
\titleName{Theory of Stochastic Drift}
\subTitleName{A Guided Tour}
\makerVersionName{2.4}

\newcommand{\BitStrings}{\{0,1\}^n}
\newcommand{\OneMax}{\textsc{OneMax}}
\newcommand{\LeadingOnes}{\textsc{LeadingOnes}}
\newcommand{\Needle}{\textsc{Needle}}
\newcommand{\Plateau}{\textsc{Plateau}}
\newcommand{\OneOneEA}{$(1+1)$ EA\xspace}
\newcommand{\RLS}{RLS\xspace}
\newcommand{\pplus}{p^{\scriptscriptstyle\leftarrow}}
\newcommand{\pminus}{p^{\scriptscriptstyle\rightarrow}}

\begin{document}

\makeHeader

	\textbf{Abstract.} In studying randomized search heuristics, a frequent quantity of interest is the first time a (real-valued) stochastic process obtains (or passes) a certain value. Commonly, the processes under investigation show a bias towards this goal, the \emph{stochastic drift}. Turning an iteration-wise expected bias into a first time of obtaining a value is the main result of \emph{drift theorems}. This thesis gives an introduction into the theory of stochastic drift, providing examples and reviewing the main drift theorems available. Furthermore, the thesis explains how these methods can be applied in a variety of contexts, including those where drift theorems seem a counter intuitive choice. Later sections examine related methods and approaches.

	This document is available as \href{https://hpi.de/friedrich/docs/scripts/22_AllDrift/index.html}{HTML}; furthermore, a copy was uploaded to \href{}{arXiv (TBD)}.

\tableofcontents

\clearpage

\clearpage

\section*{Preface}
\addcontentsline{toc}{section}{Preface}

	This document gives an introduction to the theory of stochastic drift, as developed by the community researching the theory of randomized search heuristics. For researchers new to the area (but with some basic familiarity with probability theory and random processes), the early sections provide a gentle introduction into the main theorems and sample applications. Later sections give more specialized theorems for particular applications. Seasoned researchers might turn directly to later sections, browsing the list of drift theorems for many settings which provides further pointers to the literature, as well as remarks on details of the techniques and their relation to similar approaches.

	Furthermore, this document serves as ``a summarized and systematic presentation of the candidate's own work'' in partial fulfillment of the requirements for \emph{Habilitation} at the Digital Engineering Faculty of the University of Potsdam, Germany.

\clearpage

\section{What is Stochastic Drift?}
\label{sec:intro}

	Suppose that you win a million dollars in a lottery and that you start spending your winnings.
	You observe that you spend on average 10.000 dollars per day. How long will your lottery winnings last?
	Intuitively, you would divide the million you won by 10.000 and estimate that your winnings would last for 100 days.
	But that feels like confusing a random process with a deterministic one.
	Well, yes, but the good news is: There is a theorem that tells us that 100 days is the mathematically precise answer, even when the process is randomized.
	Even better, if you gain money on some days (say, by playing in a casino) but still, in expectation, your balance goes down by 10.000 per day, the conclusion still holds.
	There can even be dependencies between the earnings and spendings of different days (say, on a day after a big earning, you gamble higher amounts than otherwise).
	The theorem showing that this is the case is called the \emph{additive drift theorem} (see \buildRef{thm:classicDrift:additiveDriftUpper}).
	The term \emph{drift} refers to the difference between two successive values of the process, and the term \emph{additive} refers to the requirement that the drift is, in expectation, bounded by an additive constant.

	A similar setting to that of the process described above is the well-known \emph{coupon collector} process, defined as follows.
	Suppose you want to collect coupons until you have one of each of $n$ different colors. Each day you get one coupon the color which is chosen uniformly at random and independent of the other days; in particular, you may gain coupons of a color which you already received before. How long does it take until you have a complete set of at least one coupon of each color? For the analysis, note that in the first iteration you get a new color of coupon with certainty (since you do not have any coupons yet). This changes over time: once you already collected exactly half of the colors, the probability of gaining a new one is only one in two. Once you have already gained 90 percent, it is down to one in ten, and so on. Or, flipped around: if you are only half the way from your goal, you only have half the chance of making progress, and if you are 10 percent away from the goal, you make 10 percent of the progress.
	This is a multiplicative expected progress (the progress is a multiple of the current state of the process) and the \emph{multiplicative drift theorem} (see \buildRef{thm:classicDrift:multiplicativeDrift}) can be used to analyze exactly this setting.	
	The theorem also holds when the number of coupons gained in each iteration is random (for example, if I gain every color of coupon with a random chance of some value $p$ in each iteration), and it even holds when there is a possibility of losing coupons. Furthermore, it also gives an upper concentration bound.

	As these two examples show, we analyze the expected progress of a single step of the given random process in order to find the first time the random process reaches a target state (the so-called ``first-hitting time''). This brings us to the following description of drift theory:

\begin{noteEN}{}{}
	\textbf{Drift theory} is a collection of theorems to turn iteration-wise expected gains into expected first-hitting times.
\end{noteEN}

	The first drift theorem, the \emph{additive} drift theorem, was introduced by He and Yao \cite{DBLP:journals/ai/HeY01}, based on an intricate theorem by Hajek~\cite{hajek_1982}.
	He and Yao applied their theorem in the context of analyzing randomized search heuristics (RSHs), such as evolutionary algorithms (EAs), which work by the principle of variation (mutating solutions by random changes) and selection (accepting improvements and rejecting worsenings). Drift theory gained a lot of traction in the EA theory community after the \emph{multiplicative} drift theorem was introduced by Doerr, Johannsen, and Winzen \cite{DBLP:conf/gecco/DoerrJW10}.
	Their proof used additive drift, but a proof not relying on Hajek's result was given shortly after by Doerr and Goldberg \cite{DBLP:conf/ppsn/DoerrG10a}.
	Since then, drift theory has been the dominant method for formally analyzing RSHs, easing their analysis significantly over analyses not arguing via drift.
	For example, the main result of Droste \cite{DBLP:conf/gecco/Droste04} on noisy optimization, spanning an entire paper, was reproven in a more general form by Giessen and Kötzing \cite{DBLP:journals/algorithmica/GiessenK16} on a single page.

	Drift theorems find application in the analysis of a plethora of different settings, ranging from randomized optimization over approximation algorithms to further stochastic processes (see \buildRef{sec:classicDrift} and the second part of \buildRef{sec:driftWithoutDrift}). Key to the applicability is to model the problem as a search for the time until a random process reaches a target state. However, in spite of the versatility of drift theorems, there are only very few results outside of the theory of RSHs applying drift theory \cite{DBLP:conf/esa/BertschingerLMM20,DBLP:journals/rsa/GoldbergLR20,DBLP:conf/podc/KosowskiU18}. Given the versatility of the approach within the area of randomized optimization, it is likely that a higher visibility of these theorems could benefit further research communities.

\subsection{A Guide to this Document}

	This document presents an overview over drift theory. What theorems are available? How can they be applied? What pitfalls abound when using drift theory? Concretely, the contents of this document are as follows.

	\buildRef{sec:classicDrift} gives the two already mentioned drift theorems (additive and multiplicative) formally. These are by far the two most important drift theorems and the section includes many examples of their use from a diverse range of settings.

	\buildRef{sec:potentialFunctions} discusses the most important technique of making drift theorems applicable: With potential functions, random processes can be mapped to fulfill the requirements of drift theorems, and this section discusses heuristics of how to do this.

	One main example for how potential functions can be used to make a process exhibit drift is the analysis of \emph{unbiased} random walks. These walks have, by definition, a drift of $0$, but nonetheless drift theory can be used to analyze such processes. This is detailed in \buildRef{sec:driftWithoutDrift}.

\begin{noteEN}{}{}
	For researchers new to the area, Sections~2 to~4 give a brief but well-rounded introduction to the field. Further sections provide deeper material, extending the applicability and discussing the finer points of drift theory.
\end{noteEN}

	\buildRef{sec:advancedDriftTheorems} provides a long list of available drift theorems, including the famous \emph{variable} drift theorem, as well as many other drift theorems tailored to various settings of drift. Some example applications and discussions on the relation between the different theorems give an overview of the currently available drift theorems.

	A special case of drift is exhibited by monotone processes (processes that cannot go back). This is a specific branch of analysis which was developed independently of the other drift theorems; we discuss the corresponding theorems in \buildRef{sec:fitnessLevelMethod}.

	While classic drift theorems give statements about how long it takes for a process to reach a certain state, the dual question is to ask what state to expect after a given number of iterations. Also this area has theorems not unlike drift theorems, and we discuss them in \buildRef{sec:fixedBudget}.

	In \buildRef{sec:averagingDrift} we consider the very technical side of drift theorems. We contrast and discuss different ways in stating the drift theorems and point to pitfalls in applying drift theorems without checking all conditions.

	Finally, \buildRef{sec:notation} introduces some notation used in this work, before the author gives acknowledgments in \buildRef{sec:acknowledgments}.

\clearpage

\section{A Gentle Introduction to Classic Drift Theorems}
\label{sec:classicDrift}

	In this section we present two classic drift theorems, the additive and the multiplicative drift theorem (see \buildRef{sec:classicDrift:basicTheorems}). These two theorems are the basis for most analyses made by drift theory, and many more advanced drift theorems are variations of these two core examples. We give a number of instructive examples for how and when the theorems are helpful in \buildRef{sec:classicDrift:simpleApplications}, followed by two more complex examples in \buildRef{sec:classicDrift:moreComplexExamples}; at the end of this section, in \buildRef{sec:classicDrift:EAExamples}, we provide two classic applications from the theory of randomized search heuristics. Note that, for this section, we refrain from diving into the technically most powerful statements and present simpler versions of these theorems; stronger versions can be found in \buildRef{sec:advancedDriftTheorems}.

\subsection{The Additive and the Multiplicative Drift Theorems}
\label{sec:classicDrift:basicTheorems}

	We state the two most commonly used drift theorems. The first drift theorem is the additive drift theorem, which requires a uniform bound on the expected change of a random process. It is due to \cite{DBLP:journals/ai/HeY01, DBLP:journals/nc/HeY04}. The very general version given here is due to \cite{DBLP:journals/tcs/KotzingK19}, where also an instructive proof can be found. We give another proof on \buildRef{thm:advancedDriftTheorems:additiveTimeConditionedUpper}.

\begin{theoremEN}[label=thm:classicDrift:additiveDriftUpper]{Additive Drift, Upper Bound}{}
	
		Let $(X_t)_{t \in \natnum}$ be an integrable process over $\realnum$, and let $T = \inf\set{t \in \natnum}{X_t \leq 0}$.
		Furthermore, suppose the following two conditions hold (non-negativity, drift).
		\begin{description}
			\item[(NN)] For all $t \leq T$, $X_t \geq 0$.
			\item[(D)] There is a $\delta > 0$ such that, for all $t \lt T$, it holds that $\Ew{X_t - X_{t+1} \mid X_0,\ldots, X_t} \geq \delta$.
		\end{description}
		Then
		$$
			\Ew{T} \leq \frac{\Ew{X_0}}{\delta}.
		$$
	
\end{theoremEN}

	The condition \referDefined{(D)} is the \emph{drift condition}, this is where we require the additive progress towards the target state~$0$. Note that we require to have drift for all possible histories $X_0,\ldots, X_t$ of the process. In many applications, we have a Markov chain, which implies that conditioning on the history is equivalent on conditioning on $X_t$ only. See \buildRef{sec:averagingDrift} for a detailed discussion on what to condition on.

	The condition \referDefined{(NN)} requires \emph{non-negativity} of the process. We cannot allow the process to assume smaller values than the target $0$ as demonstrated by the following example.

\begin{beispielEN}{Additive Drift and Processes Reaching Negative Numbers}{}
	Suppose our process starts with $X_0 = 5$ and, in each iteration deterministically, the process decreases by $2$. Then the expected time (in fact, the deterministic time) until $X_t \leq 0$ is exactly $3$. If we want to apply \buildRef{thm:classicDrift:additiveDriftUpper} we use that the expected gain is $2$, so the conclusion suggests an expected time of $2.5$. This incorrect conclusions comes from the disregard for \referDefined{(NN)}.

	We can amplify this effect with the following example. Let $(X_t)_{t \in \natnum}$ be a random process with $X_0 = 1$ and, for all $t$, with probability $1-1/n$, $X_{t+1} = X_t$ and otherwise $X_{t+1} = -n+1$; the expected time until $X_t \leq 0$ is $n$ (since it follows a geometric distribution with probability $1/n$), while the expected gain is $1$, for which the additive drift theorem would suggest an expected time of $1$.

\end{beispielEN}

	Note that there are also additive drift theorems that remove the condition \referDefined{(NN)} and instead incorporate an additional term in the conclusion, see \buildRef{thm:advancedDrift:additiveDriftUpperBetter}.

	The additive drift theorem also allows for a corresponding lower bound as follows \cite{DBLP:journals/ai/HeY01, DBLP:journals/nc/HeY04,DBLP:journals/tcs/KotzingK19}.
	
	In Theorem~3 of \cite{DBLP:conf/gecco/KotzingST11}, this theorem was used to show a lower bound to derive an asymptotically tight run time analysis of an evolutionary algorithm. Another application can be found in Theorem~4 of \cite{DBLP:conf/foga/0001KLNS17,DBLP:journals/tcs/FriedrichKLNS20}.

\begin{theoremEN}[label=thm:classicDrift:additiveDriftLower]{Additive Drift, Lower Bound}{}
	
		Let $(X_t)_{t \in \natnum}$ be an integrable process over $\realnum$, and let $T = \inf\set{t \in \natnum}{X_t \leq 0}$.
		Furthermore, suppose the following conditions (bounded steps, drift).
		\begin{description}
			\item[(B)] There is a $c \gt 0$ such that, for all $t \lt T$, it holds that $\Ew{|X_t - X_{t+1}| \mid  X_0,\ldots, X_t} \leq c$.
			\item[(D)] There is a $\delta \gt 0$ such that, for all $t \lt T$, it holds that $\Ew{X_t - X_{t+1} \mid  X_0,\ldots, X_t} \leq \delta$.
		\end{description}
		Then
		$$
			\Ew{T} \geq \frac{\Ew{X_0}}{\delta}.
		$$
	
\end{theoremEN}

	For this lower bound we need to require \referDefined{(B)}, a bounded expected step size. This is to avoid counterexamples like the following process.

\begin{beispielEN}{Additive Drift and Unbounded Step Size}{}
	Let $(X_t)_{t \in \natnum}$ with $X_0 = 1$ and, for all $t$, with probability $1/2$, $X_{t+1} = 0$ and otherwise $X_{t+1} = 2X_t - 2\delta$. This process exhibits a drift of $\delta$, suggesting an expected time of $1/\delta$, but the true time until $X_t \leq 0$ is again geometrically distributed, this time with probability $1/2$, giving an expected time of $2$.
\end{beispielEN}

	In order to apply an additive drift theorem, one has to find a single constant $\delta$ bounding drift uniformly. However, for processes where large parts of the state space exhibit a drift very differnt from this uniform bound, stronger results can be obtained by using a drift theorem which allows for a different drift in different states of the process.

	The multiplicative drift theorem covers the case where the drift is proportional to the current value of the process. It is due to \cite{DBLP:conf/gecco/DoerrJW10}, with tail bounds given in \cite{DBLP:conf/ppsn/DoerrG10a}.

\begin{theoremEN}[label=thm:classicDrift:multiplicativeDrift]{Multiplicative Drift}{}
	
		Let $(X_t)_{t \in \natnum}$ be an integrable process over $\{0, 1\} \cup S$, where $S \subset \realnum_{\gt 1}$, and let $T = \inf\set{t \in \natnum}{X_t \leq 0}$.
		
		Assume that there is a $\delta \in \realnum_{+}$ such that, for all $s \in S \cup \{1\}$ and all $t \lt T$, it holds that
			$$\Ew{X_t - X_{t+1} \mid X_0,\ldots, X_t} \geq \delta X_t.$$
		Then
			$$\Ew{T} \leq \frac{1 + \ln \Ew{X_0}}{\delta}.$$
		Further, for all $k \gt 0$ and $s \in S \cup \{1\}$ with $\Pr{X_0 \leq s} \gt 0$, it holds that
		$$
			\Pr{T \gt \frac{k + \ln s}{\delta} \mid X_0 \leq s} \leq \mathrm{e}^{-k}\,.
		$$
	
\end{theoremEN}

	The condition \referDefined{(D)} gives a bound dependent on the history, specifically dependent on the ``current'' value of the process. Intuitively it requires that, if the process has a current value of $X_t = s$, then the drift is at least $\delta s$. In fact, the multiplicative drift theorem is frequently stated with a condition of $X_t=s$ instead of $X_0,\ldots, X_t$, and the upper bound is written as $\delta s$ instead of $\delta X_t$. See \buildRef{sec:averagingDrift} for a detailed discussion on the different ways to write a drift theorem.

	These drift theorems cover a lot of applications; the remainder of this section gives a range of usecases. Most scientists consider the drift theorems stated above first before turning to other drift theorems (see \buildRef{sec:advancedDriftTheorems} for a list and discussion of such alternatives). An incomplete list of some applications of these basic theorems, regarding the analysis of evolutionary algorithms, is as follows.
	\begin{itemize}
			\item Theorem~15 of \cite{DBLP:conf/gecco/KotzingSNO12} uses it for a simple $O(n \log n)$ bound.
			\item Similarly easy argument are given in Theorems~14 and~17 of \cite{DBLP:conf/gecco/DoerrDK15}.
			\item A number of applications is given in \cite{DBLP:conf/foga/0001KLNS17,DBLP:journals/tcs/FriedrichKLNS20}.
			\item Lemma~2 of \cite{DBLP:conf/ppsn/KotzingM12} uses the concentration bound of \buildRef{thm:classicDrift:multiplicativeDrift}.
			\item So does Theorem 9 of \cite{DBLP:conf/isaac/FriedrichKKS15} (see also Theorem 9 of \cite{DBLP:journals/tec/0001KKS17}) for the analysis of the cGA.
			\item The application in Theorem~9 of \cite{DBLP:conf/foga/FeldmannK13} is a bite more intricate argument for an upper bound via multiplicative drift.
			\item Similarly in Theorems~8 and~15 of \cite{DBLP:conf/gecco/FriedrichKKS15} for the analysis of an ant colony optimization (ACO) algorithm optimizing noisy OneMax.
	\end{itemize}

	We note that, in our applications of the drift theorems in the following, we do not show that the random processes under consideration are integrable, since this is easily observed from the context that they are defined in.

\subsection{Some Simple Applications}
\label{sec:classicDrift:simpleApplications}

	We will start by lookig at the process from \buildRef{sec:intro} about collecting coupons, a classic process analyzed in many text books on random processes. We start with a (suboptimal) analysis via additive drift.

\begin{theoremEN}[label=thm:classicDrift:couponCollectorAdditive]{Coupon Collector with Additive Drift}{}
	
		Suppose we want to collect at least one of each color of $n \in \natnum_{\geq 1}$ coupons. Each round, we are given one coupon chosen uniformly at random from the~$n$ colors. Then, in expectation, we have to collect for at most $n^2$ rounds.

\end{theoremEN}
\begin{proof}
		Let $X_t$ be the number of coupons missing after $t$ rounds and let $t \lt T$. The probability of making progress (of $1$) with coupon $t+1$ is at least $X_t/n$. In the worst case, when only one color is missing, this is still $1/n$. Thus, $\Ew{X_t - X_{t+1} \mid X_0,\ldots, X_t} = X_t/n \geq 1/n$. Since we start with $X_0=n$ missing colors, an application of \buildRef{thm:classicDrift:additiveDriftUpper} gives the desired upper bound of $n^2$ rounds, using $\delta = 1/n$.
	\end{proof}

	The analysis with additive drift completely disregards the very high probability of finding new colors while still a lot of colors are missing. Thus, the analysis with multiplicative drift gives a much better bound, as the following theorem shows. In a sense, the multiplicative drift theorem is a generalization of the classic analysis of the coupon collector process; or, vice versa, the analysis of the coupon collector process follows directly from the multiplicative drift theorem.

\begin{theoremEN}[label=thm:classicDrift:couponCollector]{Coupon Collector with Multiplicative Drift}{}
	
		Suppose we want to collect at least one of each color of $n \in \natnum_{\geq 1}$ coupons. Each round, we are given one coupon chosen uniformly at random from the~$n$ colors. Then, in expectation, we have to collect for at most $n(1+ \ln n)$ rounds. Furthermore, for all $k \in \realnum_{\gt 0}$, overshooting this time by $kn$ has a probability of at most $\mathrm{e}^{-(k+1)}$.

\end{theoremEN}
\begin{proof}
		Let $X_t$ be the number of coupons missing after $t$ rounds and let $t \lt T$. The probability of making progress (of $1$) with coupon $t+1$ is $X_t/n$. Thus, $\Ew{X_t - X_{t+1} \mid X_0,\ldots, X_t} = X_t/n$. An application of \buildRef{thm:classicDrift:multiplicativeDrift} gives the desired result.
	\end{proof}

	A lower bound can be derived with an appropriate lower bounding multiplicative drift theorem (see \buildRef{thm:advancedDrift:couponCollectorLowerBound}). Since the process is monotone, both an upper and a lower bound can be derived with the \emph{fitness level method}, see \buildRef{thm:fitnessLevelMethod:couponCollectorLowerBoundFLM}.

	Using an analogous proof as in \buildRef{thm:classicDrift:couponCollector}, one can directly analyze a \emph{generalized version} of the coupon collector process as follows.

\begin{theoremEN}[label=thm:classicDrift:generalizedCouponCollector]{Generalized Coupon Collector}{}
	
		Suppose we want to collect at least one of each color of $n \in \natnum_{\geq 1}$ coupons. For each color of coupon and each round, we get this color of coupon with probability at least $p \in (0, 1]$. Then, in expectation, we have to wait for at most $(1+ \ln n)/p$ rounds. Furthermore, for all $k \in \realnum_{\gt 0}$, overshooting this time by $k/p$ has a probability of at most $\mathrm{e}^{-(k+1)}$.

\end{theoremEN}
\begin{proof}
		Let $X_t$ be the number of coupons missing after $t$ rounds and let $t \lt T$. The expected progress is $\Ew{X_t - X_{t+1} \mid X_0,\ldots,X_t} \geq p X_t$, since the expected number of missing coupons that we get in the next iteration is $p X_t$. An application of \buildRef{thm:classicDrift:multiplicativeDrift} gives the desired result.
	\end{proof}

	Note that this generalized version does not make any assumptions on how many coupons we get per iteration, or whether these indicator random variables are in any way correlated.

	We now turn to the well-known geometric distribution. The typical computation for its expectation involves modifying infinite sums. Using drift, the computation is rather simple. Furthermore, our analysis allows for processes where the probability of success changes over time and depends on the history, but a uniform bound on this probability is known.

\begin{theoremEN}[label=thm:classicDrift:geometricDistribution]{Geometric Distribution}{}
	
		Let $(X_t)_{t \in \natnum}$ be some random process where, in each iteration, a \emph{success} event happens with some probability, possibly dependent on the history of the process; we let $S(X_0,\ldots,X_t)$ denote the success event. Then the following estimates hold for all $p \in (0,1]$.
		\begin{enumerate}
				\item If, for each $t$, $\Pr{S(X_0,\ldots,X_t) \mid X_0,\ldots,X_t} \leq p$, then the expected time until any success event happened is at least $1/p$.
				\item If, for each $t$, $\Pr{S(X_0,\ldots,X_t) \mid X_0,\ldots,X_t} \geq p$, then the expected time until any success event happened is at most $1/p$.
				\item If, for each $t$, $\Pr{S(X_0,\ldots,X_t) \mid X_0,\ldots,X_t} = p$, then the expected time until any success event happened is exactly $1/p$.
		\end{enumerate}

\end{theoremEN}
\begin{proof}
		For all $t \in \natnum$, let $X_t$ be $0$ if a success event has happened within the first $t$ iterations, and $1$ otherwise. Let $t \lt T$. Then $\Ew{X_t - X_{t + 1} \mid X_0,\ldots,X_t} = p$. Thus, \buildRef{thm:classicDrift:additiveDriftUpper} and \buildRef{thm:classicDrift:additiveDriftLower} give us the corresponding bounds of $1/p$ iterations until the first success event.
	\end{proof}

	Next we consider a sequence of fair coin tosses. Known as the \href{https://en.wikipedia.org/wiki/Gambler%27s_fallacy}{gambler's fallacy} is the believe that a sequence of ``heads'' makes the occurrence of ``tails'' more likely. Quite in contrast to this, for any given $k \in \natnum$ there will be an occurrence of $k$ ``heads'' in a row if the coin is tossed sufficiently often. In the following theorem we derive exactly how long we have to wait in expectation for such an event to happen. In the proof we apply the additive drift theorem not going down towards $0$, but going up to a value of $k$. Since the additive drift is symmetrical, we can use it in either direction equally.

\begin{theoremEN}[label=thm:classicDrift:winningStreaks]{Winning Streaks}{}
	
		Let $k \in \natnum$ be given. Consider flipping a fair coin indefinitely. Then the expected number of iterations until the first time that \emph{heads} comes up $k$ times in a row is (exactly) $f(k) = 2^{k+1} - 2$.

\end{theoremEN}
\begin{proof}
		For all $t \in \natnum$, let $R_t$ be the length of the current streak of heads after $t$ iterations ($R_t = 0$ if in iteration $t$ we got tails, as well as before any coin flip at $t=0$). In the following computation, we will condition on a value for the current search point, which is equivalent to conditioning on the history since our process is a discrete Markov chain (see \buildRef{sec:averagingDrift} for details). 
		Let $X_t = f(R_t)$ be our process for which we aim to show drift. Let $i \in \natnum$ be given. If our current streak of heads is $i$, then in the next iteration one of two things happens: either we lose all progress, falling to a potential of $f(0) = 0$, or we gain $f(i+1) - f(i)$. Each happens with probability $1/2$, so we have
		\begin{align*}
			\Ew{X_{t+1} - X_t \mid X_t = f(i)}
			&= \frac{1}{2}f(i+1)  - f(i)\\
			&= 2^{i+2}/2 - 2/2 - (2^{i+1} - 2)\\
			&= 1.
		\end{align*}
		Thus, using \buildRef{thm:classicDrift:additiveDriftUpper} and \buildRef{thm:classicDrift:additiveDriftLower} together, going up instead of down, we get an expected number of iterations of $f(k) = 2^{k+1}-2$ to reach a streak of $k$ heads.
	\end{proof}

	Note that the potential function in the last proof, as in many places where potential functions are used, is not intuitive, so let us discuss where this potential function comes from. We decide we want to set up for additive drift, since the additive drift theorem gives both lower and upper bounds. Since any potential function that gives an additive drift can be normalized to give an additive drift of $1$, we search for a potential function that gives a drift of exactly $1$. From the two possible outcomes of the coin flipping process in each iteration, we now get the condition of $f(i+1)/2 - f(i) = 1$ for the potential $f$. In this case, this is a straightforward and easy to solve recurrence relation, so that with the (arbitrary) setting of $f(0) = 0$ we get the desired formula for $f$.

	For a more in-depth discussion of potential functions and their use for the application of drift theorem, see \buildRef{sec:potentialFunctions}.

\subsection{More Complex Problems}
\label{sec:classicDrift:moreComplexExamples}

	In contrast to the previous applications of drift theorems, the following examples consider processes that are not Markov. This is no problem for the drift theorems, but the user now has to make sure that all bounds hold regardless of history, not just with respect to the current value of the process.

	Our next example is a randomized algorithm for finding, in expectation, a $2$-approximation of the classical vertex cover problem. For an undirected graph $(V, E)$, a subset $C \subseteq V$ such that, for all $\{u, v\} \in E$, $u$ or $v$ is in~$C$ is called a \emph{vertex cover}. By \buildRef{thm:classicDrift:additiveDriftUpper}, we easily bound the expected size of the vertex cover that the algorithm constructs.

\begin{theoremEN}[label=thm:classicDrift:vertexCover]{Vertex Cover Approximation}{}
	
		Given an undirected graph, iteratively choose an uncovered edge and add uniformly at random an endpoint to the cover. Then, in expectation, the resulting cover is a $2$-approximation of an optimal vertex cover of the given graph.

\end{theoremEN}
\begin{proof}
		Let a graph $G$ be given. Furthermore, fix a minimum vertex cover $C$. For all $t$, let $D_t$ be the set of vertices chosen by the algorithm after $t$ iterations. Let $X_t$ be $0$ if $D_t$ is a vertex cover, and otherwise let $X_t$ be the number of vertices of $C$ that are not in $D_t$. Clearly, the algorithm terminates exactly when $X_t = 0$. Furthermore, in each step and regardless of history, the algorithm selects a vertex from $C$ with probability at least $1/2$, since, for every edge of $G$, at least one of the endpoints is in $C$. We get $\Ew{X_t - X_{t + 1} \mid X_0,\ldots,X_t} \geq \frac{1}{2}$. Hence, using \buildRef{thm:classicDrift:additiveDriftUpper}, we get that the algorithm terminates in expectation after choosing $2|C|$ vertices.
	\end{proof}

	The next example considers a simple randomized sorting algorithm. This and similar sorting algorithms were considered by Scharnow, Tinnefeld, and Wegener~\cite{DBLP:journals/jmma/ScharnowTW04} (before the advent of drift theory). The analysis via the multiplicative drift theorem is short, easy and intuitive.

\begin{theoremEN}[label=thm:classicDrift:randomSorting]{Random Sorting}{}
	
		Consider the sorting algorithm which, given an input array $A$ of length $n \in \natnum_{\geq 1}$, iteratively chooses two different positions of the array uniformly at random and swaps them if and only if they are out order. Then the algorithm obtains a sorted array after $\Theta(n^2 \log n)$ iterations in expectation.

\end{theoremEN}
\begin{proof}
		For all $i, j \in [n]$ with $i \lt j$, an ordered pair $(i,j)$ is called an \emph{inversion} if and only if $A[i] \gt A[j]$. Note that the maximum number of inversions is $\binom{n}{2}$. Let $X_t$ be the number of inversions after $t \in \natnum$ iterations, and let $A_t$ denote the array after that iteration. If the algorithm chooses a pair which is not an inversion, nothing changes. If the algorithm chooses an inversion $(i,j)$, then this inversion is removed; for any other inversion, only indices $k \in [i..j]$ are relevant.
		If $A_t[k] \lt A_t[j]$ ($\lt A_t[i]$), then $(i,k)$ is an inversion before and after the swap, while $(k,j)$ is neither an inversion before nor after the swap; similarly for $A_t[k] \gt A_t[i]$ ($\gt A_t[j]$). Finally, if $A_t[j] \lt A_t[k] \lt A_t[i]$, then $(i,k)$ and $(k,j)$ are inversions before the swap but are not afterwards. Overall, this shows that the number of inversions goes down by at least $1$ whenever the algorithm chooses an inversion for swapping, regardless of history.
		
		Let $t$ be such that $A_t$ is not sorted. Since the probability of the algorithm choosing an inversion is $X_t/\binom{n}{2}$, we get $\Ew{X_t - X_{t+1} \mid X_0,\ldots,X_t} \geq X_t/\binom{n}{2}$. An application of \buildRef{thm:classicDrift:multiplicativeDrift} gives the desired upper bound.

		Regarding the lower bound, consider the array $A$ which is almost sorted but the first and second element are swapped, the third and fourth, and so on. Then the algorithm effectively performs a coupon collector process on $n/2$ coupons, where each has a probability of $1/\binom{n}{2}$ to be collected. This takes an expected time of $\Omega(n^2 \log n)$ with a proof analogous to that of \buildRef{thm:advancedDrift:couponCollectorLowerBound}.
	\end{proof}

\subsection{Classic Results for Evolutionary Algorithms}
\label{sec:classicDrift:EAExamples}

	The basic algorithm we want to analyze is the \OneOneEA; it proceeds as follows (see also \buildRef{sec:algorithms}).

\begin{center}
\includegraphics[width=80mm]{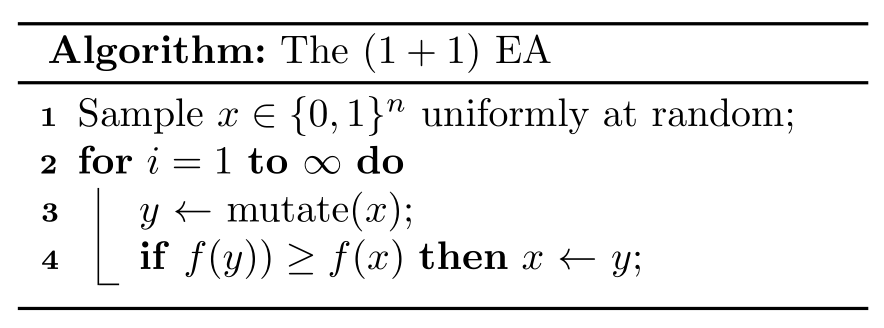}

\end{center}

	The \OneOneEA minimizing a function $f\colon \BitStrings \rightarrow \realnum$. Mutation flips each bit independently with probability $1/n$.

	The algorithm is set up to maximize the given function $f$; by turning the inequality around, we get the analogous algorithm for minimization.

\begin{theoremEN}[label=thm:classicDrift:OneOneEAOneMax]{\OneOneEA on \textsc{OneMax}}{}

		Consider the \OneOneEA maximizing the fitness function \textsc{OneMax}. Then the expected time for the algorithm to find the optimum is $\bigO{n \log n}$ iterations.

\end{theoremEN}
\begin{proof}
		For all $t$, let $X_t$ be the Hamming distance to the target of the current individual after $t$ iterations. We want to use the multiplicative drift theorem and estimate drift as follows. If the currently best search point has a Hamming distance of $s$, then, for each bit $i$ of the $s$ missing positions, the event of flipping position $i$ and no other when producing offspring will result in an accepted offspring with a distance of $1$ less to the target. These events are disjoint (since only one bit flips) and each has a probability of $1/n \cdot (1-1/n)^{n-1} \geq 1/(en)$; each decreases the distance to the target by $1$. Since the \OneOneEA does not accept worsenings, no other event can contribute negatively to the drift, so we can pessimistically assume a contribution of $0$ to the drift in all other cases. Thus, we get
		\begin{align*}
		\Ew{X_{t} - X_{t+1} \mid X_0,\ldots,X_t} \geq X_t / (en).
		\end{align*}
		An application of \buildRef{thm:classicDrift:multiplicativeDrift} gives the desired upper bound since $X_0 \leq n$.
	\end{proof}

\begin{theoremEN}[label=thm:classicDrift:OneOneEALeadingOnes]{\OneOneEA on \textsc{LeadingOnes}}{}

		Consider the \OneOneEA maximizing the fitness function \textsc{LeadingOnes}. Then the expected time for the algorithm to find the optimum is $\bigO{n^2}$ iterations.

\end{theoremEN}
\begin{proof}
		For all $t$, let $X_t$ be the number of leading ones of the current individual after $t$ iterations (i.e.~the fitness). We want to use the additive drift theorem and estimate drift as follows. Improving the fitness of the current individual requires flipping its first $0$ and none of the previous positions. There are at most $n-1$ previous positions, so the probability is at least $1/n \cdot (1-1/n)^{n-1} \geq 1/(en)$. An improvement is an improvement by at least $1$. Since the \OneOneEA does not accept worsenings, no event can contribute negatively to the drift, so we can pessimistically assume a contribution of $0$. Thus, we get
		\begin{align*}
		\Ew{X_{t} - X_{t+1} \mid X_0,\ldots,X_t} \geq 1 / (en).
		\end{align*}
		An application of \buildRef{thm:classicDrift:additiveDriftUpper} gives the desired upper bound since $X_0 \geq 0$ and the target is at $n$.
	\end{proof}

	Note that much more precise bounds are known for the optimization time of the \OneOneEA on \textsc{LeadingOnes}, see \buildRef{thm:fitnessLevelMethod:leadingOnes}.

\clearpage

\section{The Art of Potential Functions}
\label{sec:potentialFunctions}

	Drift theorems can be applied to random processes on $\realnum$. For the analysis of randomized algorithms, this typically means that one has to map the state of the algorithm to a real number, so that the resulting process will be a process on $\realnum$. Such a mapping is called \emph{potential function} and we already saw multiple in \buildRef{sec:classicDrift}. Especially in the proof of \buildRef{thm:classicDrift:winningStreaks} we saw that sometimes unintuitive potential functions can lead to very strong results.	In fact, one could say that the art of applying drift theorems is in choosing the right potential function. Later in this work, for example in the proof of \buildRef{thm:driftWithoutDrift:variance-two-barrier-hitting-time}, we see yet more intricate potential functions. In this section, we want to discuss a few cardinal examples.

\subsection{A Simple Heuristic for Choosing Potential Functions}

	There is a very important rule of thumb to designing potential functions: \textbf{better search points} should have \textbf{better potential}. The following example showcases this.

\begin{beispielEN}[label=ex:potentialFunctions:betterMeansBetter]{Better Search Points with Better Potential}{}
	Let $(X_t)_{t \in \natnum}$ be a Markov chain with $X_0 = 2$ and, for all $t$, if $X_t = 2$ then, with probability $1-1/100$, $X_{t+1} = 0$ and otherwise $X_{t+1} = 1$; if $X_t=1$ then $X_{t+1}=0$ with probability $1/100$ and $X_{t+1}=1$ otherwise. Additive drift provides an upper bound of an expected $200$ iterations to reach $0$ (since the lowest drift of $1/100$ is encountered in state $1$ and we start in state $2$).

	This process is very much misleading in that state $1$ sounds like it is closer than $2$ to the target of $0$, when actually it is not. Consider the potential function $f(0)=0$, $f(1) = 100$ and $f(2) = 2$. Now both states $1$ and $2$ have an expected drift of exactly $1$ towards the target $0$, and we start in a state with potential $2$, so we get an expected time of $2$ to reach $0$ from the additive drift theorem.
\end{beispielEN}

	From this example we see that what seems ``natural'' (because some process on the reals presents itself) might not be the best for drift. In fact, as we will see later in this section, a potential can turn drift away from the optimum into drift towards the optimum.

	The example can be generalized to arbitrary processes: on time-homogeneous Markov chains, the best potential for getting a tight bound with the additive drift theorem is the potential which assigns each state the time until finding the target from starting in that state. The next theorem from \cite{DBLP:journals/nc/HeY04} makes this formal.

\begin{theoremEN}[label=thm:potentialFunctions:useExpectedTime]{Expected Time as Potential}{}
	
		Let $\CalX$ be some state space and let $(X_t)_{t \in \natnum}$ be a time-homogeneous Markov chain on $\CalX$ and let $O \subseteq \CalX$ be a set of targets. For any $x \in \CalX$, let $T(x)$ be the random variable describing the number of steps until reaching an element in $O$ when starting in $x$, and suppose that all such $T_x$ have finite expectation. We define
		$$
		g\colon \CalX \rightarrow \realnum, x \mapsto \Ew{T(x)}.
		$$
		Then, for all $t$ with $X_t \not\in O$,
		$$
		\Ew{g(X_{t}) - g(X_{t+1}) \mid t \lt T(X_0)} = 1.
		$$

\end{theoremEN}
\begin{proof}
		Since $(X_t)_{t \in \natnum}$ is a time-homogeneous Markov-chain, let an operator $\theta$ be given such that, for all $t \in \natnum$, $X_{t+1} = \theta(X_t)$. For all $i \in \natnum$, we use $\theta^i$ to denote the $i$-times self-composition of $\theta$. In particular, for all $t \in \natnum$, we have 
		$$
		T(X_t) = \min_{i \in \natnum} \theta^i(X_t) \in O = \min_{i \in \natnum} X_{t+i} \in O.
		$$
		For all $t \in \realnum$, conditional on $X_t \not\in O$ (equivalently: $t \lt T(X_0)$) we thus have
		$$
		T(X_{t+1}) = \min_{i \in \natnum} X_{t+i+1} \in O = \left(\min_{i \in \natnum} X_{t+i} \in O\right) - 1 = T(X_{t}) - 1.
		$$
		In particular,
		\begin{align*}
		\Ew{g(X_{t}) - g(X_{t+1}) \mid t \lt T(X_0))} & = \Ew{\Ew{T(X_t)} - \Ew{T(X_{t+1})} \mid t \lt T(X_0) }\\
		& = \Ew{\Ew{T(X_t)} - \Ew{T(X_{t})-1} \mid t \lt T(X_0)}\\
		& = 1.
		\end{align*}
		This shows the claim.
	\end{proof}

	The theorem is interesting for understanding what a good potential should be; in order to apply a drift theorem it is, however, completely useless: We could now use upper and lower additive drift theorems and the proven drift of $1$ to derive an expected time of $\Ew{g(X_0)}$ to find an element of $O$ when starting in $X_0$. But $g(X_0)$ is defined to be the expected time to find an element from $O$ when starting in $X_0$, so we arrived where we started.

	Note that, in general, after mapping a Markov chains with a potential function, the resulting process is not necessarily a Markov chain any more. This is not a problem at all, since a well-formulated drift theorem does not need the requirement of the process being Markov chain, see \buildRef{sec:averagingDrift} for a discussion.

	When considering optimization algorithms, it is sometimes easy to show that the \emph{distance to the optimum} or directly the \emph{fitness} decreases in expectation in each step. This means that this would make for a good potential function; but sometimes this expected change is either too weak, too hard to analyze or even negative. In this case, more inventive potential functions are sought, which is the concern of the remainder of this section.

	Before we dive into developing concrete potential functions, we discuss normalizing potential functions. This will later make one decision very easy for us.

\begin{theoremEN}[label=thm:potentialFunctions:normalizeDrift]{Normalizing Additive Drift}{}
	
		Let $\CalX$ be some state space and let $(X_t)_{t \in \natnum}$ be a random process on $\CalX$. Let $c \in \realnum_{\gt 0}$ and let $g\colon \CalX \rightarrow \realnum$ be any potential function such that 
		  $$\Ew{g(X_{t}) - g(X_{t+1}) \mid g(X_0), \ldots, g(X_t)} \geq c.$$
		Then there is a potential function $\overline{g}\colon \CalX \rightarrow \realnum$ such that 
		  $$
		  \Ew{\overline{g}(X_{t}) - \overline{g}(X_{t+1}) \mid \overline{g}(X_0), \ldots, \overline{g}(X_t)} \geq 1,
		  $$
		and $\Ew{\overline{g}(X_{0})} = \Ew{g(X_{0})}/c$.

\end{theoremEN}
\begin{proof}
		We choose $\overline{g} = g/c$.
	\end{proof}

	The theorem shows that, whenever there is any potential function at all amenable to analysis by additive drift, there is one with an expected drift of $1$. Note that normalization has no impact on the resulting time bound, since the starting value is scaled correspondingly (and now equals the time bound derived by the additive drift theorem).

\subsection{Potential Functions for Two-Part Drift}
\label{sec:potentialFunctions:glueTogether}

	In our first example, we want to ``glue together'' two drift regimes.

\begin{beispielEN}[label=thm:potentialFunctions:glueTogether]{Gluing Together Fitness Functions}{}
	Let $(X_t)_{t \in \natnum}$ be a discrete integrable process on $[0,n]$ with $X_0 = n$ and let $k \in [0..n]$. Let $T$ be the first time $t$ such that $X_t = 0$. Suppose that we have 
	$$
	\Ew{X_{t} - X_{t+1} \mid X_t \geq k} \geq 2
	$$
	and
	$$
	\Ew{X_{t} - X_{t+1} \mid 0 \lt X_t \lt k} \geq 1.
	$$
	In other words: while potential is high, we get a drift of $2$, for small values only a drift of $1$.
	Further assume that, if $X_t \lt k$, then also $X_{t+1}\lt k$. We want to show that 
	$$
	\Ew{T} \leq \frac{n+k}{2}.
	$$
	
	In order to make use of the stronger drift for large values of the process, we choose the potential function 
	$$
	g\colon \realnum \rightarrow \realnum, x \mapsto \begin{cases}
	x, 			&\mbox{if }x\lt k;\\
	(x+k)/2,	&\mbox{otherwise.}
	\end{cases}
	$$
	We have $g(X_0) = (n+k)/2$, so it is sufficient to show a drift of at least $1$ to get our desired bound. This holds trivially for $X_t \lt k$, heavily relying on the fact that this implies $X_{t+1} \lt k$. For the following reasoning, note that, for all $x \in \realnum$, $g(x) \leq (x+k)/2$. For $X_t \geq k$, we see
	\begin{align*}
	\Ew{g(X_{t}) - g(X_{t+1}) \mid X_t \geq k} & \geq \Ew{(X_{t}+k)/2 - (X_{t+1}+k)/2 \mid X_t \geq k}\\
	&	=  \Ew{X_{t} - X_{t+1} \mid X_t \geq k} / 2\\
	& \geq 1.
	\end{align*}
	This shows a drift of $1$ in potential and thus gives the desired bound using \buildRef{thm:classicDrift:additiveDriftUpper}.
\end{beispielEN}

	Note that the potential function $g$ used in the proof above is concave, so in the main derivation in the proof we could have reasoned with Jensen's Inequality. Thus, this approach generalizes to other concave potential functions.

\subsection{$(1+1)$ EA on Linear Functions}

	One of the most famous examples of an analysis with potential functions and drift theory is the analysis of the \OneOneEA (see \buildRef{sec:algorithms}) on linear functions. In fact, the paper introducing the multiplicative drift analysis \cite{DBLP:journals/algorithmica/DoerrJW12} used this drift theorem with a suitable potential function to show an upper bound of $(1+o(1))\; 1.39e \; n \ln(n)$ on arbitrary linear functions. This bound was later improved to $(1\pm o(1))\; e \; n \ln(n)$, including a matching lower bound, in \cite{DBLP:journals/cpc/Witt13}, with a more intricate potential function that crucially depended on the concrete linear function.

	Here we give a simple proof from \cite{DBLP:journals/algorithmica/DoerrJW12}, showcasing the use of potential functions which achieves a bound of $(1+o(1))\; 4e \; n \ln(n)$. We will further restrict the linear functions to have no duplicate weights, avoiding to treat this edge case.

\begin{theoremEN}[label=thm:potentialFunctions:oneoneLinFunNoDuplicates]{$(1+1)$ EA on Linear Functions, no duplicated weights}{}
	
		Let $f$ be any linear function without duplicate weights. The expected time until the \OneOneEA on $f$ samples the optimum for the first time is $(1+o(1))\; 4e \; n \ln(n)$.

\end{theoremEN}
\begin{proof}
		Let $w_1, \ldots, w_n$ be the weights of $f$.
		We can assume, without loss of generality, that the weights are ordered decreasingly and are positive, that is,
		$$
		w_1 \gt w_2 \gt \ldots \gt w_n \gt 0.
		$$

		Now we define our potential function as follows. Let $g\colon \{0,1\}^n \rightarrow \realnum$ be such that, for all $x \in \{0,1\}^n$, 
		$$
		g(x) = \sum_{i=1}^n (2-i/n)(1-x_i).
		$$
		This potential essentially awards a potential weight between $1$ and $2$ to any incorrectly set bit, where higher potential weight for a bit position corresponds with higher weight of this bit position in the objective function.

		For each $t \geq 0$, let $X_t$ be the current best bit string found by the \OneOneEA after $t$ iterations. We want to show that there is multiplicative drift in $(g(X_t))_{t \in \natnum}$. To that end, fix $t \in \natnum$. Let $I = \set{i \leq n}{x_i=0}$ be the set of positions where the current best bit string has a $0$. We define a number of events that we want to distinguish; these events will be a partition of the entire event space at iterations $t$. 
		\begin{itemize}
				\item For each $i \in I$, let $A_i$ be the event that bit $i$ is the only $0$-bit which is flipped by mutation and the final offspring is accepted.
				\item Let $C$ be the event that at least 2 of the $0$ bits are flipped.
				\item Let $D$ be the event that no $0$ bit is flipped.
		\end{itemize}

		Let $\Delta = g(X_t) - g(X_{t+1})$ be the drift. We can now get the following breakup of the expected drift by the law of total expectation.
		\begin{align*}
		\Ew{\Delta \mid g(X_t)} = \;\; &\Ew{\Delta \mid g(X_t), C}\Pr{C \mid g(X_t)}\\
			&+ \Ew{\Delta \mid g(X_t), D}\Pr{D \mid g(X_t)}\\
			&+ \sum_{i \in I} \Ew{\Delta \mid g(X_t), A_i}\Pr{A_i \mid g(X_t)}.
		\end{align*}

		If no $0$ bit flips (event $D$), then any flip of a $1$ bit will result in worse offspring, which will be discarded; thus, $\Ew{\Delta \mid g(X_t), D} = 0$.

		Now consider event $C$. At least two $0$ bits flip, leading to an increase in potential of at least $2$. There are at most $n$ many $1$ bits, each flipping with a probability of $1/n$ and the potential associated with that bit is strinctly less than $2$. Thus we lose (in expectation) strictly less than a potential of $2$ from flipping $1$ bits, but we gain a potential of at least $2$ from flipping $0$ bits. Furthermore, the result might be accepted or not, and if it is not accepted, then $\Delta = 0$. Note that, the more $1$ bits are flipped, the less likely the offspring is accepted, so there is a negative correlation between number of $1$ bits flipped and the probability of acceptance. This shows $\Ew{\Delta \mid g(X_t), C} \geq 0$.

		Now we consider the events $A_i$. Note that, for all $i \in I$, $P(A_i) \geq (1-1/n)^{n-1}/n \geq 1/en$, since the event that bit $i$ flips and no other is a subevent of $A_i$. Not flipping $n-1$ bits has a probability of $(1-1/n)^{n-1}$, and flipping a specific bit has a probability of $1/n$. Crucially, it is impossible that any bit with position $\lt i$ flips and the offspring is accepted, since it has a (strictly!) higher weight than bit $i$ (and bit $i$ is the only $0$ bit that flips).

		Overall, we have the following.
		\begin{eqnarray*}
		\Ew{\Delta \mid g(X_t)} 
		  & \geq & \sum_{i\in I} \Ew{\Delta \mid g(X_t), A_i}\Pr{A_i \mid g(X_t)}\\
		  & \geq & \sum_{i\in I} \frac{1}{en} \Ew{\Delta \mid g(X_t), A_i}\\
		  & \geq & \frac{1}{en} \sum_{i\in I} \left[ \left(2- \frac{i}{n} \right) - \sum_{j = i+1}^n\frac{1}{n} \left(2 - \frac{j}{n} \right)\right]\\
		  & = & \frac{1}{en} \sum_{i\in I} \left[ \left(2- \frac{i}{n} \right) - \frac{1}{n} \left(2(n-i) -  \frac{\sum_{j = i+1}^n j}{n} \right) \right]\\
		  & = & \frac{1}{en} \sum_{i\in I} \left[ 2- \frac{i}{n}  - \frac{2(n-i)}{n} +  \frac{\sum_{j = i+1}^n j}{n^2} \right]\\
		  & = & \frac{1}{en} \sum_{i\in I} \left[ \frac{- i + 2i}{n} +  \frac{n(n+1) - (i+1)i}{2n^2} \right]\\
		  & \geq & \frac{1}{en} \sum_{i\in I} \left[ \frac{i}{n} +  \frac{n(n+1) - ni}{2n^2} \right]\\
		  & = & \frac{1}{en} \sum_{i\in I} \left[\frac{i + (n+1)/2-i/2}{n} \right]\\
		  & = & \frac{1}{en} \sum_{i\in I} \left[\frac{1}{2} + \frac{i+1}{2n} \right]\\
		  & \geq & \frac{1}{en} \sum_{i\in I} \frac{1}{2}\\
		  & = & \frac{|I|}{2en}.
		\end{eqnarray*}
		Since we have $|I| \geq g(X_t)/2$, we get a multiplicative drift with drift constant $\delta = 4e\; n$. Using $g(X_0) \leq 2n$, we can apply \buildRef{thm:classicDrift:multiplicativeDrift} to get 
		$$
		\Ew{T} \leq (1+o(1)) \; 4e \; n \ln(n).
		$$
	\end{proof}

\subsection{Designing a Potential Function via Step-Wise Differences}
\label{sec:potentialFunctions:stepWisePotential}

	In this section we give a method for finding a suitable potential function by defining the potential differences of ``neighboring''  states. Note that this method was used to find the proof given in \buildRef{thm:classicDrift:winningStreaks} and is also the basis of the work in \buildRef{subsec:advancedDriftTheorems:finitSearchSpaces}, with details in \cite{DBLP:conf/ppsn/KotzingK18}. Furthermore, an early version of this method for overcoming negative drift was given in \cite{DBLP:conf/gecco/GiessenK14,DBLP:journals/algorithmica/GiessenK16}.

	We consider the optimization of a test function which looks like $\OneMax$ for most of the search space, but around the optimum is a plateau of constant fitness. This is a fitness function defined as follows, for a given parameter $k \in \natnum$ (for $x \in \BitStrings$ we use $|x|_1$ to denote the number of $1$s in $x$).
	$$
	\Plateau_k\colon \BitStrings \rightarrow \BitStrings, x \mapsto
	\begin{cases}
	|x|_1,	&\mbox{if }|x|_1 \leq n-k \mbox{ or }|x|_1 = n;\\
	n-k,	&\mbox{otherwise.}
	\end{cases}
	$$
	This function is maximized by the bit string $1^n$. All bit strings with a distance between $1$ and $k$ to the optimum have identical fitness, so there is no guiding signal towards the optimum on that so-called plateau.

	We want to study the random search heuristics \OneOneEA and \RLS on the plateau function (see \buildRef{sec:algorithms}).

	For the \OneOneEA, we can analyze the performance as follows. Within $\bigO{n \log n}$ iterations the algorithm will have found a search point on the plateau, that is, at distance at most $k$ to the optimum (this follows analogously to the analysis of \OneOneEA on $\OneMax$). From now on at most $k$ bits will be incorrect, and correcting exactly those $k$ bits and no others has a probability of 
	$$
	\left(\frac{1}{n}\right)^k\left(1-\frac{1}{n}\right)^{n-k} \geq \frac{1}{en^k}.
	$$
	Thus, for $k \geq 2$, the total expected optimization time will be $\bigO{n^k}$. This reasoning disregards the analysis of the random walk performed by the algorithm on the plateau.

	The search heuristic Random Local Search (\RLS) exchanges the mutation operator of the \OneOneEA for an operator which flips exactly one bit. The analysis of the \OneOneEA above explicitly makes use of large steps which are not performed by \RLS, so a different analysis is required. We now have to understand the random walk on the plateau as an essential part to finding the optimization time, and analyzing it with drift theory provides a nice example of the power of potential functions. In fact, it is somewhat surprising that, also for this fitness function, an analysis with drift theory can find a good bound on the expected optimization time: Using the fitness function as the potential function, the expected drift for ``inner'' points on the plateau (where all neighbors are also points on the plateau) is $0$.

	We want to develop a potential function that is $0$ at the optimum. We aim for an expected drift for the \RLS optimizing $\Plateau$ of at least $1$, given that \buildRef{thm:potentialFunctions:normalizeDrift} shows that we can always find a normalized drift function. For reasons of symmetry, all bit strings at the same distance to the optimum should have the same potential, so we now wonder what should be the potential of a bit string with exactly $d$ many $0$s (so $d$ is the Hamming distance to the optimum). Let us simplify and consider first the case of $d=1$.

	On the plateau, any change is accepted by \RLS. Thus, if there is only one incorrect bit, \RLS will correct it with probability $1/n$ and otherwise lose a different bit with probability $1-1/n$. If we think about the potential difference between $d=0$ and $d=1$ as $a(0)$, and the potential difference between $d=1$ and $d=2$ as $a(1)$, then the expected gain in potential is given by
	$$
	\frac{1}{n} \cdot a(0) - \left(1 - \frac{1}{n}\right) \cdot a(1) = \frac{a(0) - (n-1)a(1)}{n}.
	$$
	We want this quantity to be a least $1$, so, for fixed $a(0)$, we get $a(1) \leq (a(0)-n) / (n-1)$. This is a rather complex term, but note that for $a(0) \geq 2n$, we can choose $a(1) = a(0)/(2n)$, a much simpler term.

	Turning to the general case of arbitrary $d$, we get an expected drift of
	$$
	\frac{d}{n} \cdot a(d-1) - \left(1 - \frac{d}{n}\right) \cdot a(d) = \frac{d\cdot a(d-1) - (n-d)a(d)}{n}.
	$$
	Thus, if again $a(d-1) \geq 2n$, we could work with $a(d) = a(d-1) / (2n)$. Inductively, we now have $a(d) = a(0) / (2n)^d$.

	Note that we need this to hold for $d$ starting at $d=0$ up until $d=k-1$, since on the plateau we cannot go outward from being exactly $k$ away. We can now choose $a(0)$ to suit all requirements. Concretely, for all $d \in [0..k-1]$ we need
	$$
	\frac{a(0)}{(2n)^d} = a(d) \geq 2n,
	$$
	so $a(0) \geq (2n)^{d+1}$. This restriction is strongest for $d = k-1$, leading to $a(0) \geq (2n)^k$.

	After we established the size of the different \emph{gaps} between different states of the algorithm, we now define the potential function (denoting the number of $0$s in a bit string $x$ as $|x|_0$) as
	$$
	g\colon \BitStrings \rightarrow \realnum_{\geq 0}, x \mapsto \sum_{i=0}^{|x|_0-1} a(i).
	$$
	This way of defining a potential function as a sum of ``gaps'' has a the advantage that the difference in potential of similar search points is easy to compute.

	We now get to the final proof, where we will also need to worry about the ``easy'' part of the search space (again ``gluing together'' the drift regimes as in \buildRef{sec:potentialFunctions:glueTogether}).

\begin{theoremEN}[label=thm:potentialFunctions:plateau:upperBound]{\RLS on Plateau, upper bound}{}
	
		Let $k \geq 2$. The expected time for \RLS to optimize $\Plateau$ is $\bigO{n^{k}}$.

\end{theoremEN}
\begin{proof}
		Let $(X_t)_{t \in \natnum}$ be the current search point of \RLS after $t$ iterations. For all $d \in \natnum$ we define $a(d) = (2n)^{k-d}$ and $g_0 = \sum_{i=0}^{k-1} a(i)$. Now we define a potential function $g$ as
		$$
		g: \BitStrings \rightarrow \realnum_{\geq 0}, x \mapsto 
		\begin{cases}
		\sum_{i=0}^{|x|_1-1} a(i), 	&\mbox{if }|x|_0 \leq k;\\
		g_0 + (|x|_0 - k)n, &\mbox{if }|x|_0 \gt k.
		\end{cases}		
		$$
		Intuitively, we artificially distort the gaps where the plateau does not provide a fitness signal and use unit gap sizes in the easy part. Note that this is suboptimal for the easy part, but the impact on the overall bound of the theorem will only be in lower order terms, since the time to cross the plateau dominates.
		
		Let now $t$ be given and let $x = X_t$ and $x' = X_{t+1}$. We are interested in bounding $\Ew{g(x)-g(x')}$. We will use the law of total expectation and make a case distinction on $|x|_0$.

		First, let $d \lt k$ be given and consider an iteration of \RLS on a bit string with exactly $d$ many $0$s (where it flips exactly $1$ bit). \RLS either gains a potential of $a(d-1)$, with probability $d/n$, or loses a potential of $a(d)$, otherwise. We have
		\begin{align*}
		\Ew{g(x) - g(x') \mid |x|_0 = d} & = \frac{d}{n} a(d-1) - \frac{n-d}{n}a(d)\\
		& = \frac{d}{n} \cdot n^{k-d+1} - \frac{n-d}{n} \cdot n^{k-d}\\
		& = (dn-n+d) \cdot n^{k-d-1}.
		\end{align*}
		For $d=1$ this equals $n^{k-2} \geq 1$; for $d \gt 1$ this is at least $(d-1)n \cdot n^{k-d-1} \geq n^{k-d}$. Since $d \leq k$, this value is at least $1$.

		We now consider the case of $d = k$. Note that, in this case, we cannot lose potential, as the selection of \RLS discards any strictly worse search point. Thus, we have
		\begin{align*}
		\Ew{g(x) - g(x') \mid |x|_0 = k} & = \frac{k}{n} a(k-1)\\
		& = \frac{k}{n} \cdot n^{k-k+1}\\
		& = k \geq 1.
		\end{align*}
		
		Finally, we consider the case of $d \gt k$. Also in this case we cannot lose potential; we have
		\begin{align*}
		\Ew{g(x) - g(x') \mid |x|_0 = d} & = \frac{d}{n} \cdot n\\
		& = d \geq 1.
		\end{align*}
		Thus, \buildRef{thm:classicDrift:additiveDriftUpper} gives an upper bound on the expected optimization time of the maximal potential value of $g_0 + n(n - k) \leq n \cdot 2n^{k} + n^2 = \bigO{n^{k}}$.
	\end{proof}

	For $k$ constant, we can use drift theory with a similar potential function to find a matching lower bound. Note that this a particular strength of the additive drift theorem: it comes with a matching lower bound without additional requirements to the random process (such as concentration in each step).

\begin{theoremEN}[label=thm:potentialFunctions:plateau:lowerBound]{\RLS on $\Plateau$, lower bound}{}
	
		Let $k \geq 2$ be constant. The expected time for \RLS to optimize $\Plateau$ is $\Omega(n^{k})$.

\end{theoremEN}
\begin{proof}
		Let $(X_t)_{t \in \natnum}$ be the current search point of \RLS after $t$ iterations. For all $d \in \natnum$ we define $a(d) = ((n-k)/k)^{k-d}$ and $g_0 = \sum_{i=0}^{k-1} a(i)$. Now we define a potential function $g$ as
		$$
		g: \BitStrings \rightarrow \realnum_{\geq 0}, x \mapsto 
		\begin{cases}
		\sum_{i=0}^{|x|_0-1} a(i),	&\mbox{if }|x|_0 \leq k;\\ 
		g_0, 						&\mbox{if }|x|_0 \gt k.
		\end{cases}
		$$
		Intuitively, we ignore the run time outside of the plateau, since it does not contribute to the asymptotic bound.
		
		Let now $t$ be given and let $x = X_t$ and $x' = X_{t+1}$. We are interested in bounding $E[g(x)-g(x')]$, this time we want an upper bound. We will again use the law of total expectation and make a case distinction on $|x|_0$.
		
		First, let $d \lt k$ be given and let the current bit string $x$ have $|x|_0 = d$. Since \RLS will flip exactly $1$ bit, we either gain a potential of $a(d-1)$ (with probability $d/n$) or lose a potential of $a(d)$, otherwise. From $d \leq k$ we see $d(n-k) \leq k(n-d)$ and thus $$
		\frac{d}{n} \; \frac{n-k}{k} \leq \frac{n-d}{n}.
		$$
		Now we can derive
		\begin{align*}
		\Ew{g(x) - g(x') \mid |x|_0 = d} & = \frac{d}{n} a(d-1) - \frac{n-d}{n}a(d)\\
		& = \frac{d}{n} \cdot ((n-k)/k)^{k-d+1} - \frac{n-d}{n} \cdot ((n-k)/k)^{k-d}\\
		& \leq \frac{n-d}{n} \cdot ((n-k)/k)^{k-d} - \frac{n-d}{n} \cdot ((n-k)/k)^{k-d}\\
		& = 0.
		\end{align*}
		An upper bound of $0$ might seem surprising, but in this area of the search space the distortion by the potential function is large enough to arrive at negative drift.
		
		We now consider the case of $d = k$. In this case, we cannot lose potential. We have
		\begin{align*}
		E[g(x) - g(x') \mid |x|_0 = k] & = \frac{k}{n} a(k-1)\\
		& = \frac{k}{n} \cdot ((n-k)/k)^{k-k+1}\\
		& = (n-k)/n \leq 1.
		\end{align*}
		
		Finally, we consider the case of $d \gt k$. Since in this case all neighboring search points have the same potential, we again have a drift of $0$.
		\begin{align*}
		\Ew{g(x) - g(x') \mid |x|_0 = d} & = 0.
		\end{align*}
		The initial potential is $g_0$ with a probability of at least some constant $c$.
		Thus, \buildRef{thm:classicDrift:additiveDriftLower} gives a lower bound on the expected optimization time of the initial potential value of $c g_0 \geq c(n/k)^k = \Omega(n^{k})$.
	\end{proof}

	As can be seen from the proof, the lower bound extends to super-constant $k$ as $\Omega((n/k)^k)$. Note that both this lower bound and the upper bound of $\bigO{n^k}$ are no longer optimal for super-constant $k$. In particular, the extreme case of $k=n$ is known as the $\Needle$ function (see \cite{DBLP:journals/ec/GarnierKS99} for the first analysis on $\Needle$).

\subsection{Further Potential Functions}

	The literature knows many more example applications of potential functions in order to allow for the applications of drift theory. For example, in \cite{DBLP:conf/gecco/0001KN0R23} a clever potential function is used to incorporate a state of the algorithm into the general progress of the algorithm towards the goal. In Section~4.2 of \cite{DBLP:conf/gecco/DoerrKLL17}, the potential function essentially has two parts to allow for a unified drift argument, rather than arguing over two phases.

	Further interesting potential functions can be found in \cite{DBLP:journals/algorithmica/DoerrDK18}. One function incorporates speed (a self-adjusting parameter) of the algorithm and the distance to the optimum into a single potential. Another combines the distances in different dimensions in a suitably scaled way to arrive at a useful potential function.

\subsection{Conclusion}

	While building a potential function is more of an art than a science, there are heuristics which can help.
	\begin{itemize}
			\item As we saw in \buildRef{ex:potentialFunctions:betterMeansBetter}, states that are ``closer'' to the target should have potential ``closer'' to that of the target; in fact, as shown by \buildRef{thm:potentialFunctions:useExpectedTime}, the most accurate potential assigns each search point the ``distance'' to the target.
			\item Drift might be different in different parts of the search space; in this case, we can use potential functions to ``glue together'' these parts, as showcased by \buildRef{thm:potentialFunctions:glueTogether}.
			\item In \buildRef{sec:potentialFunctions:stepWisePotential} we saw one way of iteratively building a potential function by comparing ``neighboring'' states.
	\end{itemize}

\clearpage

\section{Going Nowhere: Drift Without Drift}
\label{sec:driftWithoutDrift}

	We encounter a surprisingly easy application of drift theory in the absence of drift. An example of a random process which does not exhibit any expected change (drift) is the \emph{Gambler's Ruin} process. By considering a transformation of the process (essentially: squaring it) the process now exhibits drift in the order of its variance, which is then amenable to analysis with drift theory.

\subsection{Unbiased Random Walks}

	We start by deriving two general corollaries, before we draw conclusions for specific random walks. The first, \buildRef{thm:driftWithoutDrift:variance-two-barrier-hitting-time}, concerns a completely unbiased random walks. The second, \buildRef{thm:driftWithoutDrift:variance-one-barrier-hitting-time}, gives the situation for random walks with one barriers.

\begin{theoremEN}[label=thm:driftWithoutDrift:variance-two-barrier-hitting-time]{Unbiased Random Walk on the Line}{}
	
		Let $n \in \natnum$, let $(X_t)_{t \in \natnum}$ be an integrable random process over $[0,n]$, and let $T = \inf\set{t \in \natnum }{ X_t \in \{0, n\}}$. Suppose that there is a $\delta \in \realnum_{+}$ such that, for all $t \lt T$, we have the following conditions (variance, drift).
		\begin{description}
			\item[(Var)] $\Var{X_{t+1} - X_t \mid X_0,\ldots,X_t} = \delta$;
			\item[(D)] $\Ew{X_{t+1} - X_t \mid X_0,\ldots,X_t} = 0$.
		\end{description}
		Then $\Ew{T} = \frac{\Ew{X_0(n - X_0)}}{\delta}$.

\end{theoremEN}
\begin{proof}
		We consider the process $Y_t = X_t(n-X_t)$.
		Note that~$T$ is the first time $t \in \natnum$ such that $Y_t = 0$. In the following, we condition on $X_0,\ldots,X_t$, a filtration that $Y_0,\ldots,Y_t$ is adapted to; this allows us to apply our drift theorems by \buildRef{thm:averagingDrift:bigImplication} while giving information not only about the value of $Y_t$, but also about $X_t$. Furthermore, for all $s \in [1..n-1]$, we have
		\begin{align*}
			\Ew{Y_t - Y_{t+1} \mid X_0,\ldots,X_t} &= \Ew{X_{t+1}^2-X_{t}^2 \mid X_0,\ldots,X_t} - n\Ew{X_{t+1}-X_{t} \mid X_0,\ldots,X_t}\\
				&= \Ew{X_{t+1}^2 \mid X_0,\ldots,X_t} - X_t^2 = \Var{X_{t+1} \mid X_0,\ldots,X_t}\\
				&= \Var{X_{t+1} - X_t \mid X_0,\ldots,X_t} = \delta.
		\end{align*}
		Thus, we have a drift of $\delta$ towards $0$. Since $Y_0 = X_0(n - X_0)$, the theorem follows from an application of \buildRef{thm:classicDrift:additiveDriftUpper}.
	\end{proof}

	Since the proof is based on the additive drift theorem, a lower bound of $\delta$ on the variance is enough for an upper bound on the expected first-hitting time and vice versa.

\begin{theoremEN}[label=thm:driftWithoutDrift:variance-one-barrier-hitting-time]{Unbiased Random Walk on the Line, One Barrier}{}
	
		Let $n \in \natnum$, let $(X_t)_{t \in \natnum}$ be an integrable random process over $[0,n]$, and let $T = \inf\set{t \in \natnum }{ X_t = n}$.
		Suppose that there is a $\delta \in \realnum_{+}$ such that, for all $t \lt T$, we have the following conditions (variance, drift).
		\begin{description}
			\item[(Var)] $\Var{X_{t+1} - X_t \mid X_0,\ldots,X_t} \geq \delta$;
			\item[(D)] $\Ew{X_{t+1} - X_t \mid X_0,\ldots,X_t} \geq 0$.
		\end{description}
		Then $\Ew{T} \leq \frac{n^2 - \Ew{X_0^2}}{\delta}$.

\end{theoremEN}
\begin{proof}
		We consider the process $Y_t = n^2 - X_t^2$.    Note that~$T$ is the first time such that $Y_t = 0$. Further note that, from \referDefined{(D)} we get, for all $t \lt T$, 
		\begin{equation}
		\Ew{X_{t + 1} \mid X_0,\ldots, X_t}^2 \geq X_t^2.\tag{$\ast$}
		\end{equation}
		As in the previous proof, we condition on $X_0,\ldots,X_t$ and implicitly use \buildRef{thm:averagingDrift:bigImplication}. We now have
		\begin{align*}
			\Ew{Y_t - Y_{t+1} \mid X_0,\ldots, X_t} & = \Ew{X_{t+1}^2 - X_{t}^2 \mid X_0,\ldots, X_t}\\
			& = \Ew{X_{t+1}^2 \mid X_0,\ldots, X_t} - X_t^2\\
			&\eqnComment{$(\ast)$}{\geq} \Ew{X_{t+1}^2 \mid X_0,\ldots, X_t} - \Ew{X_{t + 1} \mid X_0,\ldots, X_t}^2\\
			&= \Var{X_{t+1} \mid X_0,\ldots, X_t}\\
            &= \Var{X_{t+1} - X_t \mid X_0,\ldots, X_t}\\
			&\geq \delta.
		\end{align*}
		Thus, we have a drift of at least $\delta$ towards $0$. Since $Y_0 = n^2 - X_0^2$, the theorem follows from an application of \buildRef{thm:classicDrift:additiveDriftUpper}.
	\end{proof}

	Note that in neither of the two preceding theorems is the process allowed to overshoot the target. Using an additive drift theorem that allows for overshooting, like \buildRef{thm:advancedDrift:additiveDriftUpperBetter}, one can derive corresponding extensions of the above two theorems with essentially the same proof. We note that \buildRef{thm:driftWithoutDrift:variance-one-barrier-hitting-time} is tight with the following example.

\begin{beispielEN}{Fair Random Walk}{}
	Let $(X_t)_{t \in \natnum}$ be the time-homogeneous Markov-chain on $[0..n]$, where, for all $t \in \natnum$,
	\begin{enumerate}
			\item for all $i \in [1..n-1]$, $\Pr{X_{t + 1} = i+1 \mid X_t = i} = 1/2 = \Pr{X_{t + 1} = i-1 \mid X_t = i}$;
			\item the state~$0$ is reflective, that is, $\Pr{X_{t + 1} = 1 \mid X_t = 0} = 1$; and
			\item the state~$n$ is absorbing.
	\end{enumerate}
	We transform~$X$ into the fair random walk $(Y_t)_{t \in \natnum}$ on $[0..2n]$, where the states~$0$ and~$2n$ are both absorbing, such that, for all $t \in \natnum$, it holds that $X_t = |Y_t - n|$.
	
	Informally, we mirror~$X$ at~$0$ and then shift it by~$n$. Whenever this new process is at~$n$, it goes to either~$n - 1$ or~$n + 1$, each with probability~$1/2$, which results exactly in~$Y$. Note that $T = \inf\set{t \in \natnum }{ Y_t \in \{0, 2n\}} = \inf\set{t \in \natnum }{ X_t = n}$. Applying \buildRef{thm:driftWithoutDrift:variance-two-barrier-hitting-time} to~$Y$ yields $\Ew{T} = \Ew{Y_0(2n - Y_0)}$. Since $X_0 \leq n$, it holds that $Y_0 = n - X_0$. Substituting this back into the equation for~$\Ew{T}$ yields $\Ew{T} = \Ew{(n - X_0)(n + X_0)} = n^2 - \Ew{X_0^2}$, which is exactly the bound of \buildRef{thm:driftWithoutDrift:variance-one-barrier-hitting-time}.
\end{beispielEN}

\subsection{Analysis of Concrete Unbiased Random Walks}
\label{sec:driftWithoutDrift:examples}

	In this section we see several domains in which we apply our theorems about unbiased random walks. The Gambler's Ruin in \buildRef{thm:driftWithoutDrift:drunkardsWalk} is the most straightforward application of \buildRef{thm:driftWithoutDrift:variance-two-barrier-hitting-time}. A more intricate application is given in \buildRef{thm:driftWithoutDrift:recolour}, where it is used to bound the expected run time for an algorithm to find a certain coloring of a graph.

	Regarding \buildRef{thm:driftWithoutDrift:variance-one-barrier-hitting-time}, in \buildRef{thm:driftWithoutDrift:twoSat} we use it to derive an upper bound on the time for an algorithm to find a satisfying assignment for a 2-SAT formula.

	We start with the \emph{gambler's ruin}, a random walk on the line. It starts at $n$, going either one step left of one step right, each with probability $1/2$, modeling winning or losing a fair coin toss to either win or lose a coin. The question of how long it takes to either be broke ($0$ coins left) or double the starting number of coins is the simplest setting of an unbiased random walk. This process also goes by many other names, such as \emph{drunkard's walk}, \emph{random walk on a line}, or \emph{one-dimensional random walk}.

\begin{theoremEN}[label=thm:driftWithoutDrift:drunkardsWalk]{Gambler's Ruin}{}
	
		Suppose we start with $n \in \natnum$ coins and, in each iteration, uniformly at random either gain a coin or lose a coin. Then, after an expected number of exactly $n^2$ iterations, we are either broke or have reached a total of $2n$ coins.

\end{theoremEN}
\begin{proof}
		Let $(X_t)_{t \in \natnum}$ be the random sequence of the number of coins after. We have 
		$$
			\Ew{X_{t+1} - X_t \mid X_0,\ldots,X_t} = \frac{1}{2} \cdot 1 + \frac{1}{2} \cdot (-1) = 0.
		$$
		Furthermore, we have
		$$
			\Var{X_{t+1} - X_t \mid X_0,\ldots,X_t} = \frac{1}{2} \cdot 1^2 + \frac{1}{2} \cdot (-1)^2 = 1.
		$$
		Thus, we can now apply \buildRef{thm:driftWithoutDrift:variance-two-barrier-hitting-time} to get the desired result.
	\end{proof}

	Our next example considers the analysis of the run time of an algorithm. McDiarmid~\cite{mcdiarmid_1993} studies the following simple randomized algorithm called \textsc{Recolour}, for coloring a given undirected graph~$G$ with two colors such that it contains no monochromatic triangle (a subgraph on three pairwise connected vertices which are all colored with the same color). \textsc{Recolour} starts with an arbitrary $2$-coloring of~$G$. At every step, it checks whether the current coloring has a monochromatic triangle. If so, \textsc{Recolour} changes the color of one of the vertices of this triangle uniformly at random. Otherwise, the $2$-coloring has no monochromatic triangles and it is the output of \textsc{Recolour}.

	McDiarmid shows that, when \textsc{Recolour} is applied to a 3-colorable graph~$G$ (a graph that can be colored with three colors so that no two neighbors share a color), it returns a 2-coloring of~$G$ with no monochromatic triangle in expected time $\bigO{n^4}$. His analysis shows that the expected run time of the algorithm is bounded above by the expected hitting time of a random walk on the line with two absorbing states -- which is exactly the setting of \buildRef{thm:driftWithoutDrift:variance-two-barrier-hitting-time}. This analysis in turn relies on previous results on one-dimensional random walks, which usually require lengthy calculations.

	We present a simple and self-contained proof of the $\bigO{n^4}$ expected run time of the \textsc{Recolour} algorithm for finding a 2-coloring with no monochromatic triangles on 3-colorable graphs. Our proof follows the proof of McDiarmid~\cite{mcdiarmid_1993} to reduce the problem to an unbiased random walk on the line and then uses \buildRef{thm:driftWithoutDrift:variance-two-barrier-hitting-time}. A similar analysis can be used to derive an upper bound on the run time of \textsc{Recolour} on hypergraph colorings.

\begin{theoremEN}[label=thm:driftWithoutDrift:recolour]{The Recolour Algorithm}{}
	
		The expected run time of \textsc{Recolour} on a 3-colorable graph with $n \in \natnum_+$ vertices is $\bigO{n^4}$.

\end{theoremEN}
\begin{proof}
		Let $G=(V,E)$ be a 3-colorable graph, and let $\chi\colon V\rightarrow \{1, 2, 3\}$ be a 3-coloring of $G$. Let $U = \set{v\in V }{ \chi(v) \in \{1, 2\} }$  be the set of all vertices which are colored with colors $1$ and~$2$. Note that any 2-coloring of~$G$ that agrees with~$\chi$ on the vertices from~$U$ is a 2-coloring of~$G$ with no monochromatic triangles. Thus, the run time of \textsc{Recolour} is bounded from above by the expected time that \textsc{Recolour} takes to find such a coloring.

		Let $\chi_t$ be the 2-coloring found by \textsc{Recolour} at time $t \in \natnum$. Let $Y_t$ be the number of vertices~$u\in U$ such that $\chi_t(u)=\chi(u)$. The algorithm terminates when $Y_t \in \{0,|U|\}$, since agreeing on all vertices of $U$ is a coloring without monochromatic triangles, but disagreeing on all vertices from $U$ is also such a valid coloring, since the use of the colors is symmetric.

		Let $s \in [1.. |U| - 1]$ denote an outcome of $Y_t$ before the algorithm terminates. We then have that $\Pr{Y_{t+1}=Y_t+1 \mid Y_t = s} = 1/3$, as, for every monochromatic triangle, there is exactly one vertex in $u \in U$ with $\chi_t(u) \neq \chi(u)$. Similarly, $\Pr{Y_{t+1}=Y_t-1 \mid Y_t = s} = 1/3$. Thus, $Y_t$ is an unbiased random walk on the line with first-hitting time $T = \inf \set{t \in \natnum }{ Y_t \in \{0,|U|\} }$. We have
		$$
			\Var{Y_{t+1} - Y_t \mid Y_0,\ldots,Y_t} = \frac{1}{3} \cdot 1^2 + \frac{1}{3} \cdot 0^2 + \frac{1}{3} \cdot (-1)^2 = \frac{2}{3}.
		$$		
		Applying \buildRef{thm:driftWithoutDrift:variance-two-barrier-hitting-time} we get
		$$
		\Ew{T}=\frac{3\Ew{Y_0(|U| - Y_0)}}{2}\leq \frac{3n^2}{8}.
		$$
		At each step, the algorithm requires $\bigO{n^2}$ time to find a monochromatic triangle and modify this to obtain a new coloring, which concludes the proof.
	\end{proof}

	The analysis of the \textsc{Recolour} algorithm for finding 2-colorings with no monochromatic triangles appears as an exercise in~\cite{MitzenmacherUpfal:2005:ProbComp}.

	In the final example of this section e consider finding satisfying assignments of 2-SAT formulas. Papadimitriou~\cite{DBLP:conf/focs/Papadimitriou91} studies the following simple randomized algorithm that returns a satisfying assignment of a satisfiable 2-SAT formula~$\phi$ with $n$ variables and $m$ clauses within $\bigO{n^2m}$ time in expectation. The algorithm starts with a random assignment of the variables of $\phi$. At every step, the algorithm checks whether there is an unsatisfied clause for this assignment. If so, the algorithm changes the assignment of one of the variables of this assignment uniformly at random. Otherwise, the assignment is satisfying and it is the output of the algorithm.

	The analysis given is similar to the \textsc{Recolour} algorithm, and it also relies on the previous results on one-dimensional random walks. An extensive analysis of this algorithm appears in~\cite{MitzenmacherUpfal:2005:ProbComp}. Here, we present a simpler proof that uses \buildRef{thm:driftWithoutDrift:variance-one-barrier-hitting-time}.

\begin{theoremEN}[label=thm:driftWithoutDrift:twoSat]{Random 2-SAT}{}
	
		 The randomized 2-SAT algorithm, when run on a satisfiable 2-SAT formula over $n \in \natnum_+$ variables and $m \in \natnum_+$ clauses, terminates in $\bigO{n^2m}$ time in expectation.

\end{theoremEN}
\begin{proof}
		Let $\phi$ be a satisfiable 2-SAT formula and $a$  a satisfying assignment. At each time step $t \in \natnum_+$, the randomized 2-SAT algorithm finds a (not necessarily satisfying) assignment~$a_t$. Let $X_t$ be the random variable denoting the number of variables that have the same truth assignment in both $a$ and $a_t$. Let $T$ be the first time the algorithm reaches a satisfying assignment for $\phi$. Assume that a clause $x \vee y$ is not satisfied by $a_t$. Since~$a$ is a satisfying assignment, $a$ and $a_t$ differ in the assignment of at least one of the variables in this clause. Thus, 
		\begin{align*}
		\Pr{X_{t+1} = X_t+1 \mid X_0, \ldots, X_t} &\geq 1/2;\mbox{ and}\\
		\Pr{X_{t+1} = X_t-1 \mid X_0, \ldots, X_t} &\leq 1/2.
		\end{align*}
		When $a_t=a$, the algorithm terminates.
		By \buildRef{thm:driftWithoutDrift:variance-one-barrier-hitting-time} with variance bounded by $1$, we have $\Ew{T} \leq n^2$. In order to transition from~$a_t$ to $a_{t+1}$, the algorithm requires $\bigO{m}$ time (since the 2-SAT formula has $m$ distinct clauses), concluding the proof.
	\end{proof}

\clearpage

\section{The Zoo: A Tour of Drift Theorems}
\label{sec:advancedDriftTheorems}

	We have seen the basic two drift theorems, the additive drift theorem and the multiplicative drift theorem, in  \buildRef{sec:classicDrift}. In this section we provide a list of more advanced drift theorems with applications.
	\begin{enumerate}
			\item In \buildRef{subsec:advancedDriftTheorems:additive} we start by extending the additive drift theorem; we see how to avoid the requirement of non-negativity (allowing overshooting of the target) and explore different conditions for the drift.
			\item In \buildRef{subsec:advancedDriftTheorems:additiveConcentration} we give concentration results for additive drift.
			\item \buildRef{subsec:advancedDriftTheorems:additiveConcentration} provides a different view on additive drift by considering concentration.
			\item In \buildRef{subsec:advancedDriftTheorems:multiplicative} we give lower bounds for the case of multiplicative drift.
			\item While additive drift required drift to be constant and multiplicative drift required proportional drift, in \buildRef{subsec:advancedDriftTheorems:variable} we give theorems allowing for an arbitrary \emph{monotone} dependence of the drift on the current state.
			\item Somewhat different in flavor is \buildRef{sec:advancedDrift:negativeDrift}. Here we discuss drift theorems providing \emph{exponential lower bounds} given drift away from the target.
			\item In \buildRef{subsec:advancedDriftTheorems:finitSearchSpaces} we consider the special case of random processes on finite search spaces.
			\item Some settings allow drift only when far away from the target, but in the proximity of the target the drift is negative. In this case, the theorem of \buildRef{subsec:advancedDriftTheorems:headwindDrift} can offer an upper bound on the run time nonetheless.
			\item In order to derive good upper bounds even when the drift gets stronger when getting closer to the optimum, \buildRef{subsec:advancedDriftTheorems:multiUpDrift} provides a drift theorem for the case of proportionally increasing drift.
			\item A completely different approach to understanding drift is given by Wormald and briefly discussed in \buildRef{subsec:advancedDriftTheorems:wormald}.
	\end{enumerate}
	Note that there are a few novel approaches to analyzing multi-dimensional potential functions \cite{DBLP:journals/tcs/Rowe18,DBLP:conf/ppsn/JanettL22}; while the initial works are promising, they have not gained traction yet and we will not discuss them here.

\subsection{Additive Drift}
\label{subsec:advancedDriftTheorems:additive}

	We want to start with an illustrative proof for a strong version of the additive drift theorem; the proof is adapted from the proof of Theorem~2.3.1 in \cite{DBLP:series/ncs/Lengler20}.

\begin{theoremEN}[label=thm:advancedDriftTheorems:additiveTimeConditionedUpper]{Additive Drift, Upper Bound, Time Condition}{}
	
		Let $(X_t)_{t \in \natnum}$ be a stochastic process on $\realnum$ with deterministic $X_0$, and let $T = \inf\set{t \in \natnum}{X_t \leq 0}$. Suppose that there is a $\delta \gt 0$ so that we have the following conditions (drift, non-negativity).
		\begin{description}
			\item[(D)] For all $t$ with $\Pr{t \lt T} \gt 0$, $\Ew{X_t - X_{t+1} \mid t \lt T} \geq \delta$.
			\item[(NN)] For all $t \leq T$, $X_t \geq 0$.
		\end{description}
		We have
		$$
			\Ew{T} \leq X_0/\delta.
		$$

\end{theoremEN}
\begin{proof}
		We first replace the process $(X_t)_{t \in \natnum}$ with a process $(X'_t)_{t \in \natnum}$ such that, for all $t \leq T$, we have $X'_t = X_t$, and for all $t \gt T$ we have $X'_t = X'_{t-1}$. Both processes are $\leq 0$ at the same time, but $(X'_t)_{t \in \natnum}$ does not change after that.
		We have that \referDefined{(NN)} gives $X_T = 0$, and, thus, $X_t =0$ for all $t \geq T$. Together with \referDefined{(NN)} we thus have
		\begin{equation}
		\forall t \in \natnum: X_t \geq 0.\tag{NN'}
		\end{equation}
		Furthermore, for all $t \lt T$ with $\Pr{t \lt T} \gt 0$, we have $\Ew{X_t - X_{t+1} \mid t \geq T} = 0$.
		
		We now have, for all $t \in \natnum$ with $\Pr{t \lt T} \gt 0$,
		\begin{align*}
		\Ew{X_t - X_{t+1}}	& = \Ew{X_t - X_{t+1} \mid t \lt T}\Pr{t \lt T} + \Ew{X_t - X_{t+1} \mid t \geq T}\Pr{t \geq T} \\
					& = \Pr{t \lt T} \Ew{X_t - X_{t+1} \mid t \lt T}\\
					& \eqnComment{(D)}{\geq} \Pr{t \lt T} \delta \\
					& = \delta\Pr{T \gt t}.
		\end{align*}
		The first equality is the law of total expectation; the second follows from $X_{t} = X_{t+1}$ for $t \geq T$. Note that the overall inequality holds trivially for all $t \in \natnum$ such that $\Pr{t \lt T} = 0$, so it holds for all $t$.  Explicitly, for all $t \in \natnum$ we have
		\begin{equation}
		\Pr{T \gt t} \leq \frac{1}{\delta} \; \Ew{X_t - X_{t+1}}.\tag{$\ast$}
		\end{equation}

		Since $T$ takes only values in $\natnum \cup\{\infty\}$, we have
		$$
		\Ew{T} = \sum_{i=0}^\infty\Pr{T \gt i}.
		$$
		We want to use this to compute $\Ew{T}$. For all $n \in \natnum$ we have
		$$
		\sum_{t=0}^n \Pr{T \gt t} \eqnComment{($\ast$)}{\leq} \frac{1}{\delta} \sum_{t=0}^n \left( \Ew{X_t} - \Ew{X_{t+1}}\right) = \frac{1}{\delta} \left(X_0 - \Ew{X_{n+1}}\right) \eqnComment{(NN')}{\leq} \frac{X_0}{\delta}.
		$$
		Since all partial sums are upper bounded by $X_0/\delta$, so is the infinite sum.
	\end{proof}

	Note that we can turn the proof around to get the analogous version for a lower bound. Again the proof is essentially taken from the proof of Theorem~1 in \cite{DBLP:series/ncs/Lengler20}. Note that it uses the somewhat strong assumption of a bounded search space, whereas \buildRef{thm:classicDrift:additiveDriftLower} only requires a bound on the size of each step.

\begin{theoremEN}[label=thm:advancedDriftTheorems:additiveTimeConditionedLower]{Additive Drift, Lower Bound, Time Condition}{}
	
		Let $(X_t)_{t \in \natnum}$ be a stochastic process on $\realnum$ with deterministic $X_0$, and let $T = \inf\set{t \in \natnum}{X_t \leq 0}$. Suppose that there is a $\delta \in \realnum_{+}$ so that we have the following conditions (drift, upper bounded search space).
		\begin{description}
			\item[(D)] For all $t$ with $\Pr{t \lt T} \gt 0$, $\Ew{X_t - X_{t+1} \mid t \lt T} \leq \delta$.
			\item[(UB)] There is a $c \gt 0$ such that, for all $t \lt T$, $X_t \leq c$.
		\end{description}
		We have
		$$
			\Ew{T} \geq X_0/\delta.
		$$

\end{theoremEN}
\begin{proof}
		We first replace the process $(X_t)_{t \in \natnum}$ with a process $(X'_t)_{t \in \natnum}$ such that, for all $t \leq T$, we have $X'_t = X_t$, and for all $t \gt T$ we have $X'_t = X'_{t-1}$. Both processes are $\leq 0$ at the same time, but $(X'_t)_{t \in \natnum}$ does not change after that. Thus, for all $t \lt T$ with $\Pr{t \lt T} \gt 0$ we have $\Ew{X_t - X_{t+1} \mid t \geq T} = 0$.
		
		We now have, for all $t \in \natnum$ with $\Pr{t \lt T} \gt 0$,
		\begin{align*}
		\Ew{X_t - X_{t+1}}	& = \Ew{X_t - X_{t+1} \mid t \lt T}\Pr{t \lt T} + \Ew{X_t - X_{t+1} \mid t \geq T}\Pr{t \geq T} \\
					& = \Pr{t \lt T} \Ew{X_t - X_{t+1} \mid t \lt T}\\
					& \eqnComment{(D)}{\leq} \Pr{t \lt T} \delta \\
					& = \delta\Pr{T \gt t}.
		\end{align*}
		The first equality is the law of total expectation; the second follows from $X_{t} = X_{t+1}$ for $t \geq T$. Note that the overall inequality holds trivially for all $t \in \natnum$ such that $\Pr{t \lt T} = 0$, so it holds for all $t$. Explicitly, for all $t \in \natnum$ we have
		\begin{equation}
		\Pr{T \gt t} \geq \frac{1}{\delta} \; \Ew{X_t - X_{t+1}}.\tag{$\ast$}
		\end{equation}
			
		Since $T$ takes only values in $\natnum \cup\{\infty\}$, we have
		$$
		\Ew{T} = \sum_{i=0}^\infty\Pr{T \gt i}.
		$$
		We want to use this to compute $\Ew{T}$. For all $n \in \natnum$, we have
		$$
		\sum_{t=0}^n \Pr{T \gt t} \eqnComment{($\ast$)}{\geq} \frac{1}{\delta} \sum_{t=0}^n \left( \Ew{X_t} - \Ew{X_{t+1}}\right) = \frac{1}{\delta} \left(X_0 - \Ew{X_{n+1}}\right).
		$$
		It remains to be shown that $\Ew{X_{n+1}}$ converges to a value $\leq 0$ for $n$ going to infinity. Using the $c$ from \referDefined{(UB)}, we have for all $n \in \natnum$ that
		\begin{align*}
		\Ew{X_{n}} & = \Ew{X_{n}\mid n \lt T}\Pr{n \lt T} + \Ew{X_{n}\mid n \geq T}\Pr{n \geq T} \leq c \cdot \Pr{n \lt T} + 0 \cdot \Pr{n \geq T}\\
			& = c \Pr{n \lt T}.
		\end{align*}
		We distinguish two cases. If $\Pr{n \lt T}$ converges to $0$ for $n$ going to $\infty$, then $\Ew{X_{n}}$ converges to $0$ as desired. Otherwise, there is a non-zero probability of $T=\infty$, in which case the theorem follows directly from that.
	\end{proof}

	Both in \buildRef{thm:classicDrift:additiveDriftUpper}, the classic version of the additive drift theorem, as in the version just above, it required that the target of $0$ must be hit exactly and not overshot \referDefined{(NN)}. From \cite{Krejca:thesis:19} we have a stronger version that allows for overshooting. This is frequently helpful, for example for finding approximations, when potentially much better values than required can be achieved.

\begin{theoremEN}[label=thm:advancedDrift:additiveDriftUpperBetter]{Additive Drift, Upper Bound with Overshooting}{}
	
		Let $(X_t)_{t \in \natnum}$ be an integrable process over $\realnum$, and let $T = \inf\set{t \in \natnum }{ X_t \leq 0}$.
		Furthermore, suppose the following (drift).
		\begin{description}
			\item[(D)] There is a $\delta > 0$ such that, for all $t \lt T$, it holds that $\Ew{X_t - X_{t+1} \mid X_0,\ldots,X_t} \geq \delta$.
		\end{description}
		Then
		$$
			\Ew{T} \leq \frac{\Ew{X_0} - \Ew{X_T}}{\delta}.
		$$
	
\end{theoremEN}

	In a sense, this drift theorem is simpler than \buildRef{thm:classicDrift:additiveDriftUpper}: the requirement (NN) is dropped and the expected time is increased corresponding to the expected additional distance the process will have traveled (note that $\Ew{X_T}$ is not a positive value, since $T$ is the first point $t$ where $X_t \leq 0$). From the condition (NN) we could derive $X_T = 0$ and thus immediately recover \buildRef{thm:classicDrift:additiveDriftUpper}.

\subsection{Additive Drift: Concentration}
\label{subsec:advancedDriftTheorems:additiveConcentration}

	One of the reasons why the additive drift theorem is so general (we only really have a requirement on the expectation of change, the first moment, but not on the higher moments) is that we only get a conclusion about the expectation of the first hitting time of the target. With requirements on the higher moments we can derive concentration bounds on the expected first hitting time. This is provided by \cite{DBLP:journals/algorithmica/Kotzing16}, from which we give two different variants, one using an absolute bound on the step size \referDefined{(B)}, Theorem~2 in the cited work, and one requiring concentrated step size \referDefined{(C)}, combining Theorems~10 and~15 from the cited work.  Each time we get that there is only a very small probability of arriving significantly later than the expected time of $n/\delta$. An example application of such a concentration result for additive drift is in \cite{DBLP:conf/foga/KotzingLW15} regarding an analysis of the \OneOneEA on a dynamic version of $\OneMax$ (see Theorem~10 in the cited work). Another application is given in Theorem~5 of \cite{DBLP:conf/foga/0001KLNS17,DBLP:journals/tcs/FriedrichKLNS20}.

\begin{theoremEN}[label=thm:advancedDrift:additiveDriftUpperConcentrationBoundedSteps]{Additive Drift, Upper Concentration, Bounded Step Size}{}
	
		Let $(X_t)_{t \in \natnum}$ be an integrable process over $\realnum$, and let $T = \inf\set{t \in \natnum }{ X_t \leq 0}$.
		Furthermore, suppose that there is $c \gt 0$ such that we have the following conditions (drift, bounded steps).
		\begin{description}
			\item[(D)] There is a $\delta > 0$ such that, for all $t \lt T$, it holds that $\Ew{X_t - X_{t+1} \mid X_0,\ldots, X_t} \geq \delta$.
			\item[(B)] For all $t \geq 0$, $|X_{t+1}-X_t| \leq c$.
		\end{description}
		Let $n \in \natnum$ such that $X_0 \leq n$. Then, for all $s \geq 2n / \delta$,
		$$
			\Pr{T \geq s} \leq \exp\left( - \frac{s\delta^2}{8c^2} \right).
		$$
	
\end{theoremEN}

\begin{theoremEN}[label=thm:advancedDrift:additiveDriftUpperConcentrationConcentratedSteps]{Additive Drift, Upper Concentration, Concentrated Step Size}{}
	
		Let $(X_t)_{t \in \natnum}$ be an integrable random process over $\realnum$, and let $T = \inf\set{t \in \natnum }{ X_t \leq 0 }$.
		Furthermore, suppose that there are $\varepsilon \gt 0$ and $c \gt 0$ such that we have the following conditions (drift, concentration).
		\begin{description}
			\item[(D)] There is a $\delta > 0$ such that, for all $t \lt T$, it holds that $\Ew{X_t - X_{t+1} \mid X_0,\ldots, X_t} \geq \delta$.
			\item[(C)] For all $t \geq 0$ and all $x \geq 0$, $\Pr{|X_{t+1}-X_t| \geq x \mid X_t} \leq \frac{c}{(1+\varepsilon)^x}$.
		\end{description}
		Let $n \in \natnum$ such that $X_0 \leq n$. Then, for all $s \geq 2n / \delta$,
		$$
			\Pr{T \geq s} \leq \exp\left( - \frac{s\delta}{4} \min\left(\frac{\varepsilon}{4}, \frac{\delta\varepsilon^3}{256c} \right) \right).
		$$
	
\end{theoremEN}

	Also in \cite{DBLP:journals/algorithmica/Kotzing16} are analogous \emph{lower} bounds. Again we give two different variants, one using an absolute bound on the step size \referDefined{(B)}, Theorem~1 in the cited work, and one requiring concentrated step size \referDefined{(C)}, combining Theorems~10 and~14 from the cited work. Each time we get that there is only a very small probability of arriving significantly before the expected time of $n/\delta$.

\begin{theoremEN}[label=thm:advancedDrift:additiveDriftLowerConcentrationBoundedSteps]{Additive Drift, Lower Concentration, Bounded Step Size}{}
	
		Let $(X_t)_{t \in \natnum}$ be an integrable random process over $\realnum$, and let $T = \inf\{t \in \natnum \mid X_t \leq 0\}$.
		Furthermore, suppose that there is $c \gt 0$ such that we have the following conditions (drift, bounded steps).
		\begin{description}
			\item[(D)] There is a $\delta > 0$ such that, for all $t \lt T$, it holds that $\Ew{X_t - X_{t+1} \mid X_0,\ldots, X_t} \leq \delta$.
			\item[(B)] For all $t \geq 0$, $|X_{t+1}-X_t| \leq c$.
		\end{description}
		Let $n \in \natnum$ such that $X_0 \geq n$. Then, for all $s \leq n / (2\delta)$,
		$$
			\Pr{T \lt s} \leq \exp\left( - \frac{n^2}{8c^2s} \right).
		$$
	
\end{theoremEN}

\begin{theoremEN}[label=thm:advancedDrift:additiveDriftLowerConcentrationConcentratedSteps]{Additive Drift, Lower Concentration, Concentrated Step Size}{}
	
		Let $(X_t)_{t \in \natnum}$ be an integrable random process over $\realnum$, and let $T = \inf\{t \in \natnum \mid X_t \leq 0\}$.
		Furthermore, suppose that there are $\varepsilon \gt 0$ and $c \gt 0$ such that (drift, concentration)
		\begin{description}
			\item[(D)] there is a $\delta > 0$ such that, for all $t \lt T$, it holds that $\Ew{X_t - X_{t+1} \mid X_0,\ldots, X_t} \leq \delta$;
			\item[(C)] for all $t \geq 0$ and all $x \geq 0$, $\Pr{|X_{t+1}-X_t| \geq x \mid X_t} \leq \frac{c}{(1+\varepsilon)^x}$.
		\end{description}
		Let $n \in \natnum$ such that $X_0 \geq n$. Then, for all $s \leq n / (2\delta)$,
		$$
			\Pr{T \lt s} \leq \exp\left( - \frac{n}{4} \min\left(\frac{\varepsilon}{4}, \frac{n\varepsilon^3}{256cs} \right) \right).
		$$
	
\end{theoremEN}

	The overall situation depending on the strength of the drift is depicted in detail in \cite{DBLP:journals/algorithmica/Kotzing16}. In particular, there are three main regimes:
	\begin{enumerate}
			\item If the drift is at least $\delta \geq 1/n$, then we get high concentration of the first hitting time.
			\item If the drift is  $\delta \in [-1/n, 1/n]$ but the variance is significant, then we get to hit the optimum with constant chance within $\bigO{n^2}$ steps, see \buildRef{thm:driftWithoutDrift:variance-one-barrier-hitting-time}.
			\item If the drift is much smaller than $-1/n$, then we have negative drift and only a superpolynomially small chance to reach the optimum in polynomial time, see \buildRef{thm:advancedDrift:negativeDriftConstantStep}.
	\end{enumerate}

	The literature knows also the following theorem for bounding additive drift only relying on the variance, given by Semenov and Terkel \cite{DBLP:journals/ec/Semenov03}.

\begin{theoremEN}[label=thm:advancedDrift:additiveDriftVarianceConcentration]{Additive Drift, Upper Concentration, Bounded Variance}{}
	
		Let $(X_t)_{t \in \natnum}$ be an integrable random process over $\realnum$ with $X_0=0$.
		Furthermore, suppose the following (drift, variance).
		\begin{description}
			\item[(D)] There is $\delta > 0$ such that, for all $t \lt T$, it holds that $\Ew{X_t - X_{t+1} \mid X_0,\ldots, X_t} \geq \delta$.
			\item[(Var)] There is $c \gt 0$ such that, for all $t \geq 0$, $\Var{X_{t+1} \mid  X_0,\ldots, X_t} \leq c$.
		\end{description}
		Then, for all $\varepsilon \gt 0$, the following holds with probability $1$.
		$$
			X_t \geq t \delta - o\left(t^{0.5+\varepsilon}\right).
		$$
	
\end{theoremEN}

\subsection{Multiplicative Drift}
\label{subsec:advancedDriftTheorems:multiplicative}

	The plain multiplicative drift theorem (see \buildRef{thm:classicDrift:multiplicativeDrift}) is already very strong, in that it requires few conditions on the search space and even gives a concentration (in one direction). What it does not provide is a lower bound. One possible such bound can be found in \cite{DBLP:journals/cpc/Witt13} which we state here.

\begin{theoremEN}[label=thm:advancedDrift:multiplicativeDriftLowerMonotone]{Multiplicative Drift, Lower Bound, Monotone}{}
	
		Let $(X_t)_{t \in \natnum}$ be a discrete, integrable process over $\{0, 1\} \cup S$, where $S \subset \realnum_{> 1}$ is finite, and let $T = \inf\set{t \in \natnum }{ X_t \leq 0}$.
		
		We assume that there are $\beta,\delta \in (0,1)$ such that the following conditions (drift, monotonicity, concentration) hold for all $s \gt 1$ and $t \geq 0$ with $\Pr{X_t = s} \gt 0$.
		\begin{description}
			\item[(D)] $\Ew{X_t - X_{t+1} \mid X_t = s} \leq \delta s$.
			\item[(M)] $X_{t+1} \leq X_t$.
			\item[(C)] $\Pr{X_t - X_{t+1} \geq \beta s \mid X_t = s} \leq \beta\delta /\ln(s)$.
		\end{description}
		Then
		$$\Ew{T \mid X_0} \geq \frac{\ln (X_0)}{\delta} \cdot \frac{1-\beta}{1+\beta} \geq \frac{\ln (X_0)}{\delta} \cdot (1-2\beta).$$
	
\end{theoremEN}

	From \cite{DBLP:journals/algorithmica/DoerrDK18} we have a variant which allows for non-monotone drift. It substitutes the monotonicity with the requirement that we cannot expect more progress from first returning to bigger values of the process. Turned around, progress in any state $s$ cannot be bigger than in a state $s' \lt s$. We use the notation $(x)_+ := \max(0,x)$.

\begin{theoremEN}[label=thm:advancedDrift:multiplicativeDriftLower]{Multiplicative Drift, Lower Bound}{}
	
		Let $(X_t)_{t \in \natnum}$ be a discrete random process over $\{0, 1\} \cup S$, where $S \subset \realnum_{> 1}$ is finite, and let $T = \inf\set{t \in \natnum}{X_t \leq 1}$.
		
		We assume that there are $\beta,\delta \in (0,1)$ such that the following conditions (drift, concentration) hold for all $s \gt 1$ and $t \geq 0$ with $\Pr{X_t = s} \gt 0$.
		\begin{description}
			\item[(D)] For all $s'$ with $1 \lt s' \leq s$: $\Ew{(s' - X_{t+1})_+ \mid X_0,\ldots, X_t, X_t = s} \leq \delta s'$.
			\item[(C)] For all $s'$ with $1 \lt s' \leq s$: $\Pr{s' - X_{t+1} \geq \beta s' \mid X_0,\ldots, X_t, X_t = s} \leq \beta\delta /\ln(s')$.
		\end{description}
		Then
		$$\Ew{T \mid X_0} \geq \frac{\ln (X_0)}{\delta} \cdot \frac{1-\beta}{1+\beta} \geq \frac{\ln (X_0)}{\delta} \cdot (1-2\beta).$$
	
\end{theoremEN}

	As an alternative, we can find a lower bound when the step size is bounded. The following theorem is given in \cite{DBLP:journals/tcs/DoerrKLL20}. A further version can be found in \cite{DBLP:journals/tcs/KotzingK19}.

\begin{theoremEN}[label=thm:advancedDrift:multiplicativeDriftLowerBoundedStepSize]{Multiplicative Drift, Lower Bound, Bounded Step Size}{}
	
		Let $(X_t)_{t \in \natnum}$ be a random process over $\realnum_+$, let $x_{\mathrm{min}} \gt 0$, and let $T = \inf\set{t \in \natnum}{X_t \leq x_{\mathrm{min}}}$.
		
		We assume that there are $c,\delta \in \realnum_+$ with $x_{\mathrm{min}} \geq \sqrt{2}c$ such that the following conditions (drift, bounded step size) hold for all $t \lt T$.
		\begin{description}
			\item[(D)] $\Ew{X_t - X_{t+1} \mid X_0,\ldots,X_t} \leq \delta X_t$.
			\item[(B)] $|X_t-X_{t+1}| \leq c$.
		\end{description}
		Then
		$$\Ew{T \mid X_0} \geq \frac{1+\ln (X_0)-\ln(x_{\mathrm{min}})}{2\delta + \frac{c^2}{x_{\mathrm{min}}^2-c^2}}.$$
	
\end{theoremEN}

	Note that, for typical applications, $\delta$ is small; yet the term $2\delta$ should dominate the term $\frac{c^2}{x_{\mathrm{min}}^2-c^2}$ to give a tight bound. But this is typically not a problem: consider the setting of $\delta = \Theta(1/n)$ and $X_0 = \Theta(n)$. We can let $x_{\mathrm{min}} = \Theta(\sqrt{n})$ and suppose we can bound $c = o(\sqrt{n})$ with sufficiently high probability (which would be typical). Then the theorem lets us derive the asymptotically optimal bound of $\Omega(n \log n)$.

	As an example application we provide a lower bound for the coupon collector process (see the upper bound proven in \buildRef{thm:classicDrift:couponCollector}).

\begin{theoremEN}[label=thm:advancedDrift:couponCollectorLowerBound]{Coupon Collector, Lower Bound}{}
	
		Suppose we want to collect at least one of each color of $n \in \natnum_{\geq 1}$ coupons. Each round, we are given one coupon with a color chosen uniformly at random from the~$n$ kinds. Then, in expectation, we have to collect for at least $\Omega(n \ln n)$ iterations.

\end{theoremEN}
\begin{proof}
		Let $X_t$ be the number of coupons missing after $t$ iterations. We want to apply \buildRef{thm:advancedDrift:multiplicativeDriftLowerBoundedStepSize} and note that, since each iteration at most one coupon is gained and none is lost, we can use $c=1$ to satisfy \referDefined{(B)}. Furthermore, regarding \referDefined{(D)}, the probability of making progress (of $1$) with coupon $t+1$ is $X_t/n$. Thus, $\Ew{X_t - X_{t+1} \mid X_0,\ldots, X_t} = X_t/n$ and we can use $\delta=1/n$. We set $x_{\mathrm{min}}=\sqrt{n}$ and an application of \buildRef{thm:advancedDrift:multiplicativeDriftLowerBoundedStepSize} gives an upper bound of
		$$
		\frac{1+\ln (n)-\ln(\sqrt{n})}{\frac{2}{n} + \frac{1}{n-1}} \geq \frac{1+\ln(n)/2}{\frac{3}{n-1}} = \frac{1}{6} \cdot (n-1)\ln(n) = \Omega(n \ln n).
		$$
	\end{proof}

\subsection{Variable Drift}
\label{subsec:advancedDriftTheorems:variable}

	A more general version of \buildRef{thm:classicDrift:additiveDriftUpper} and \buildRef{thm:classicDrift:multiplicativeDrift} is the variable drift theorem, allowing for any \emph{monotone} dependency of the drift on the current state (meaning that a larger distance to the target has to imply a larger drift). It is due to \cite{DBLP:journals/ijicc/MitavskiyRC09, Johannsen:thesis:10} and was improved in \cite{DBLP:journals/tcs/RoweS14}. We give here the version from \cite{DBLP:journals/tcs/KotzingK19}, where the random process is not assumed to be discrete or Markov.

\begin{theoremEN}[label=thm:advancedDrift:variableDrift]{Variable Drift}{}
	
		Let $(X_t)_{t \in \natnum}$ be a random process over $\realnum$, $x_{\mathrm{min}} \in \realnum_{+}$, and let $T = \inf\set{t \in \natnum }{ X_t \lt x_{\mathrm{min}}}$. Additionally, let $I$ denote the smallest real interval that contains at least all values $x \geq  x_{\mathrm{min}}$ that, for all $t \leq T$, any $X_t$ can take. Furthermore, suppose that there is a function $h\colon I \rightarrow \realnum_+$ such that the following conditions (drift, monotonicity, start, non-negativity) hold for all $t \leq T$.
		\begin{description}
			\item[(D)] $\Ew{X_t - X_{t+1} \mid X_0,\ldots,X_t} \geq h(X_t)$.
			\item[(M)] The function $h$ is monotonically non-decreasing.
			\item[(S)] $X_0 \geq x_{\mathrm{min}}$.
			\item[(NN)] $X_t \geq 0$.
		\end{description}
		Then
		$$
		\Ew{T} \leq \frac{1}{h(x_{\mathrm{min}})} + \int_{x_{\mathrm{min}}}^{X_0} \frac{1}{h(z)}\,\mathrm{d}z.
		$$
	
\end{theoremEN}

	Note that the additive drift theorem is the special case of constant $h$ and the multiplicative drift theorem is the special case of linear $h$. It is surprising that the cases of additive and multiplicative drift are sufficient in many applications, but the variable drift theorem can also in these cases sometimes give tighter bounds.

	Concentration bounds for variable drift are also available \cite{DBLP:conf/isaac/LehreW14}. The same paper also gives a lower bounding variable drift theorem, which requires $h$ to be monotonically non-\emph{decreasing}, the opposite as for the upper bound. Further variants can be found in \cite{DBLP:journals/tcs/KotzingK19}, including lower bounds for step-size bounded settings.

	An example application of \buildRef{thm:advancedDrift:variableDrift} is the optimization of \textsc{LeadingOnes} by the \OneOneEA. It is known \cite{DBLP:journals/tcs/DrosteJW02} that the expected gain in fitness value per iteration, given that the current fitness value is $n-s$ (and thus $s$ away from the optimum), is (essentially) at least
	$$
	h(s) = 2 \cdot \left(1-1/n\right)^{n-s} \cdot \frac{1}{n}.
	$$
	The middle term is the probability to not lose a bit already gained; the $1/n$ is the probability to flip the left-most $0$ and the $2$ is the expected fitness gained when the two just mentioned events happen (one bit flipped, plus an expected one more bit that happens to be correctly set). The middle term can be lower-bounded by $1/e$, which allows for applying \buildRef{thm:classicDrift:additiveDriftUpper}, giving a total run time of at most 
	$$
	\frac{e}{2} \cdot n^2.
	$$
	Using the variable drift theorem directly on the bound given by $h$, a simple integration gives an upper bound on the optimization time of
	$$
	\frac{e-1}{2} \cdot n^2.
	$$
	This optimization time was first established in \cite{DBLP:conf/ppsn/BottcherDN10}. In this example, the use of the variable drift theorem improved the leading constant. Note that the bound can be derived also by other means, for example by the fitness level method, which can also be used to show tightness of this bound, see \buildRef{thm:fitnessLevelMethod:leadingOnes}.

	An essential application of a variable drift is given in the proof of Theorem 17 in \cite{DBLP:journals/algorithmica/DoerrFFFKS19}, considering the optimization of \textsc{OneMax} by an islands-based evolutionary algorithm, employing $\lambda$ islands and an exchange of individuals every $\tau$ rounds. In particular, the considered drift function is $h\colon \realnum \rightarrow \realnum$ such that, for all $s > 0$,
	$$
	h(s) = \left. \ln \lambda \middle/ \ln \!\left(\!\frac{n \; \ln \lambda}{\tau s} \!\right) \right.
	$$
	The final bound on the run time is shown to be asymptotically tight, thanks to using both upper and lower bounding variable drift theorems.

	Further uses of the variable drift theorem are given in Theorem~7 of \cite{DBLP:conf/foga/0001KLNS17,DBLP:journals/tcs/FriedrichKLNS20} and in Theorem~6 of \cite{DBLP:conf/gecco/DoerrDK16,DBLP:journals/algorithmica/DoerrDK18}.

\subsection{Negative Drift}
\label{sec:advancedDrift:negativeDrift}

	When the drift goes away from the target, we speak of \emph{negative drift}. The negative drift theorem \cite{DBLP:journals/algorithmica/OlivetoW11,DBLP:journals/corr/OlivetoW12} gives an exponential lower bound in this setting.

\begin{theoremEN}[label=thm:advancedDrift:negativeDrift]{Negative Drift}{}
	
		Let $(X_t)_{t \in \natnum}$ be a stochastic process over $\realnum$. Suppose there is an interval $[a,b] \subseteq \realnum$, two constants $\delta, \varepsilon \gt 0$ and, possibly depending on $\ell = b-a$, a function $r(\ell)$ satisfying $1 \leq r(\ell) = o(\ell / \log \ell)$ such that, for all $t \geq 0$, the following conditions (drift, concentration) hold.
		\begin{description}
			\item[(D)] $\Ew{X_{t+1} - X_t \mid X_0,\ldots, X_t; a \lt X_t \lt b} \geq \delta$.
			\item[(C)] For all $j \geq 0$, $\Pr{|X_{t+1} - X_t| \geq j \mid X_0,\ldots, X_t; a \lt X_t} \leq \frac{r(\ell)}{(1+\varepsilon)^{j}}$.
		\end{description}
		Then there is a constant $c$ such that, for $T = \min\set{t \geq 0}{X_t \leq a}$, we have
		$$
		\Pr{T \leq 2^{c \ell/ r(\ell)} \; \Big| \; X_0 \geq b} = 2^{-\Omega(\ell/ r(\ell))}.
		$$
	
\end{theoremEN}

	Note that drift goes with a strength \emph{independent of the width $\ell = b-a$ of the interval} away from the target $a$ (while the process is in the interval). A version with scaling which allows for more flexibility in this dependence is given in \cite{DBLP:journals/tcs/OlivetoW14}.

	A variant that allows for arbitrary $\varepsilon$ (with decaying guarantees) is given in \cite{DBLP:journals/algorithmica/Kotzing16} as follows. It requires a bounded step size, but in return gives a very simple and easy-to-apply bound.

\begin{theoremEN}[label=thm:advancedDrift:negativeDriftConstantStep]{Negative Drift, Bounded Step Size}{}
	
		Let $(X_t)_{t \in \natnum}$ be a stochastic process over $\realnum$, each with finite expectation, and let $n \gt 0$. Let $T = \min\set{t \geq 0}{X_t \geq n}$ and suppose there are $0 \lt c \lt n$ and $\varepsilon \lt 0$ such that, for all $t \geq 0$, the following conditions hold (drift, boundedness).
		\begin{description}
			\item[(D)] $\Ew{X_{t+1} - X_t \mid X_0,\ldots, X_t} \leq \varepsilon$.
			\item[(B)] $|X_{t+1} - X_t| \lt c$.
		\end{description}
		Then, for all $s \geq 0$, we have
		$$
		\Pr{T \leq s} = s\exp\left( - \frac{n |\varepsilon|}{2c^2} \right).
		$$
	
\end{theoremEN}

	Given as Corollary~22 in \cite{DBLP:journals/algorithmica/Kotzing16} is a second variant of the negative drift theorem. It allows for very large $r$ while still giving a super-polynomial bound for finding the target in polynomial time.

\begin{theoremEN}[label=thm:advancedDrift:negativeDriftTwo]{Negative Drift II}{}
	
		Let $(X_t)_{t \in \natnum}$ be a stochastic process over $\realnum$. Suppose there is an interval $[a,b] \subseteq \realnum$, two constants $\delta, \varepsilon > 0$ and, possibly depending on $\ell = b-a$, a function $r(\ell)$ satisfying $1 \leq r(\ell) = \exp(o(\sqrt[4]{\ell}))$ such that, for all $t \geq 0$, the following conditions hold (drift, concentration).
		\begin{description}
			\item[(D)] $\Ew{X_{t+1} - X_t \mid X_0,\ldots, X_t; a \lt X_t \lt b} \geq \varepsilon$.
			\item[(C)] For all $j \geq 0$, $\Pr{|X_{t+1} - X_t - \varepsilon| \geq j \mid X_0,\ldots, X_t; a \lt X_t} \leq \frac{r(\ell)}{(1+\delta)^{j}}$.
		\end{description}
		Then there is a constant $c$ such that, for $T = \min\set{t \geq 0}{X_t \leq a}$, we have
		$$
		\Pr{T \leq 2^{c \sqrt{\ell}} \; \Big| \; X_0 \geq b} = 2^{-\Omega(\sqrt[4]{\ell})}.
		$$
	
\end{theoremEN}

	Example applications of negative drift theorems for the analysis of evolutionary algorithms are given in the proofs of the following statements. Lemma~3 of \cite{DBLP:conf/ppsn/KotzingM12}; in Lemma~8 of \cite{DBLP:conf/isaac/FriedrichKKS15} (see also Lemma 5 of \cite{DBLP:journals/tec/0001KKS17}); Theorem 3 of \cite{DBLP:journals/tec/0001KKS17}; Lemma~6 \cite{DBLP:conf/ppsn/FriedrichKS16}; and Lemma~13 \cite{DBLP:conf/gecco/FriedrichKK16}.

\subsection{Finite State Spaces}
\label{subsec:advancedDriftTheorems:finitSearchSpaces}

	Most drift theorems consider a random walk on the real numbers, sometimes restricted to non-negative numbers. For the analysis of discrete algorithms, frequently the state space is even more restricted, in particular finite. By numbering the successive states, we can assume the state space to be $[0..n]$. For this setup we have the following drift theorem from \cite{DBLP:conf/ppsn/KotzingK18}. Note that the proof given in the paper is derived by the method of step-wise differences, see \buildRef{sec:potentialFunctions:stepWisePotential}.

\begin{theoremEN}[label=thm:advancedDrift:finiteStatesUpper]{Finite State Spaces, Upper Bound}{}
	
		Let $(X_t)_{t \in \natnum}$ be a time-homogeneous Markov chain on $[0..n]$ and let $T$ be the first time $t$ such that $X_t = 0$. Suppose there are two functions $\pplus\colon [1..n] \to [0, 1]$ and $\pminus\colon [0..n-1] \to [0, 1]$ such that, for all $t \lt T$ and all $s \in [1..n]$,
        \begin{enumerate}
            	\item $\pplus(s) > 0$,
            	\item $\Pr{X_t - X_{t+1} \geq 1 \mid X_t = s} \geq \pplus(s)$,
            	\item $\Pr{X_t - X_{t+1} = -1 \mid X_t = s} \leq \pminus(s)$ (for $s \neq n$), and
            	\item $\Pr{X_t - X_{t+1} \lt -1 \mid X_t = s} = 0$ (for $s \neq n$).
		\end{enumerate}
        Then
        $$
            \Ew{T \mid X_0} \leq \sum_{s = 1}^{X_0}\sum_{i = s}^{n}\frac{1}{\pplus(i)}\prod_{j = s}^{i - 1}\frac{\pminus(j)}{\pplus(j)}.
        $$
	
\end{theoremEN}

	A special case of this theorem is the fitness level method (see \buildRef{thm:fitnessLevelMethod:upper}), where the process is monotone; we can recover this setting by setting $\pminus$ to be constantly $0$, significantly simplifying the above formula.

	We also have the corresponding lower bound.

\begin{theoremEN}[label=thm:advancedDrift:finiteStatesLower]{Finite State Spaces, Lower Bound}{}
	
		Let $(X_t)_{t \in \natnum}$ be a time-homogeneous  Markov chain on $[0..n]$ and let $T$ be the first time $t$ such that $X_t = 0$. Suppose there are two functions $\pplus\colon \{1, \ldots, n\} \to [0, 1]$ and $\pminus\colon [0..n-1] \to [0, 1]$ such that, for all $t \lt T$ and all $s \in [1..n]$,
        \begin{enumerate}
            	\item $\pplus(s) > 0$,
            	\item $\Pr{X_t - X_{t+1} = 1 \mid X_t = s} \leq \pplus(s)$,
            	\item $\Pr{X_t - X_{t+1} \gt 1 \mid X_t = s} = 0$, and
            	\item $\Pr{X_t - X_{t+1} \leq -1 \mid X_t = s} \geq \pminus(s)$ (for $s \neq n$).
		\end{enumerate}
        Then
        $$
            \Ew{T \mid X_0} \geq \sum_{s = 1}^{X_0}\sum_{i = s}^{n}\frac{1}{\pplus(i)}\prod_{j = s}^{i - 1}\frac{\pminus(j)}{\pplus(j)}.
        $$
	
\end{theoremEN}

	Note that for processes which make steps of at most $1$ and given exact $\pminus$ and $\pplus$, the two bounds coincide.

\subsection{Headwind Drift}
\label{subsec:advancedDriftTheorems:headwindDrift}

	Sometimes drift only carries until shortly before the target, but then, close to the target, turns negative. In case only a small remaining distance needs to be bridged, and the probability of going the right way is still sufficiently high, the following \emph{Headwind} drift theorem can be used to directly get a decent bound without relying on hand crafted potential functions. The theorem was developed and applied in \cite{DBLP:conf/foga/KotzingLW15}.

\begin{theoremEN}[label=thm:advancedDrift:headwindDrift]{Headwind Drift}{}
	
		Let $(X_t)_{t \in \natnum}$ be a time-homogeneous Markov chain on $[0..n]$. Let bounds 
		$$
		p^{-}(i)\leq \Pr{X^{t+1} \leq i-1 \mid X_t = i}
		$$
		and   
		$$
		p^{+}(i) \geq \Pr{X^{t+1} \geq i+1 \mid X_t = i},
		$$
		where $0 \leq i \leq n$, be given, and define 
		$$
		\delta(i):= p^{-}(i)-\Ew{(X_{t+1}-i) \cdot \mathbb{1}[X_{t+1} > i] \mid X_t=i}.
		$$
		Assume that $\delta(i)$ is monotone increasing with respect to $i$ and let $\kappa \geq \max\set{i \ge 0 }{ \delta(i) \leq 0 }$ (noting that $\delta(0) \leq 0$). 
		The function $g\colon [0..n+1]\to \realnum_+$ is defined by
		$$
		g(i) := \sum_{k=i+1}^n \frac{1}{\delta(k)}
		$$
		for $i \geq \kappa$ (in particular, $g(n)=g(n+1)=0$), and inductively by 
		$$
		g(i) := \frac{1+(p^{+}(i+1)+p^{-}(i+1))g(i+1)}{p^{-}(i+1)}
		$$
		for $i \lt \kappa$. 

		Then it holds for the first hitting time~$T := \min\set{t \geq 0 }{ X_t=0 }$ of state~$0$ that
		$$
		\Ew{T \mid X_0} \leq g(0)-g(X_0).
		$$
	
\end{theoremEN}

	We can also get a closed expression for the expected first hitting time $\Ew{T \mid X_0}$. This expression involves the factor $\sum_{k=\kappa+1}^N \frac{1}{\delta(k)}$ that is reminiscent of the formula for the expected first hitting time of state~$\kappa$ under variable drift towards the target (see \buildRef{thm:advancedDrift:variableDrift}). For the states less than~$\kappa$, where drift away from the target holds, the product $\prod_{k=1}^{\kappa} \frac{p^{+}(k)+p^{-}(k)}{p^{-}(k)}$ comes into play. Intuitively, it represents the waiting time for the event of taking $\kappa$ consecutive steps against the drift. Since the product involves probabilities conditioned on leaving the states, which effectively removes self-loops, another sum of products must be added. This sum, represented by the second line of the expression for $\Ew{T \mid X_0}$, intuitively accounts for the self-loops.

\begin{corollaryEN}[label=cor:advancedDrift:headwindDriftClosedForm]{Headwind Drift, Closed Form}{}
	
		Let the assumptions of \buildRef{thm:advancedDrift:headwindDrift} hold. Then
		\begin{align*}
			\Ew{T \mid X^0} & \leq \left(\left(\sum_{k=\kappa+1}^N \frac{1}{\delta(k)}\right) \left(\prod_{k=1}^{\kappa} \frac{p^{+}(k)+p^{-}(k)}{p^{-}(k)} \right)\right)
			\\
			& \qquad + \left(\sum_{k=1}^{\kappa} \frac{1}{p^{-}(k)}\prod_{j=1}^{k-1} \frac{p^{+}(j)+p^{-}(j)}{p^{-}(j)}\right).
		\end{align*}
	
\end{corollaryEN}

	Note that there is a similarity between the theorems for headwind drift and those from \buildRef{subsec:advancedDriftTheorems:finitSearchSpaces}. This is because the analysis for the last steps of headwind drift is essentially an analysis brute-forcing the small interval of negative drift, which also happens in \buildRef{subsec:advancedDriftTheorems:finitSearchSpaces}.

\subsection{Multiplicative Up-Drift}
\label{subsec:advancedDriftTheorems:multiUpDrift}

	The idea of \buildRef{thm:classicDrift:multiplicativeDrift} was to have multiplicative drift going \emph{down} towards $0$. While this has many applications (owning to the fact that progress in optimization typically gets harder as better and better solutions are found), there are also a number of processes that gain in speed over time, typically making progress proportional to the current state of the process, such as rumor spreading, epidemics and population take-over. This is known as multiplicative up-drift and was studied in depth in \cite{DBLP:journals/algorithmica/DoerrK21}. The main theorem is the following.

\begin{theoremEN}[label=thm:advancedDrift:MultiUpDrift]{Multiplicative Up-Drift}{}
	
		Let $(X_{t})_{t \in \natnum}$ be a random process over  $\integers_{\geq 0}$. Let $n,k \in \integers_{\geq 1}$, $E_0 \gt 0$, $\gamma_0 \lt 1$, and $\delta \gt 0$ such that $n - 1  \leq \min\{\gamma_0 k, (1+\delta)^{-1} k\}$. Let $D_0 = \min (\lceil 100/\delta \rceil,n)$ when $\delta \le 1$ and $D_0 = \min (32,n)$ otherwise. Assume that, for all $t \geq 0$ and all $x \in [0..n-1]$ with $\Pr{X_{t} = x} \gt 0$, the following two conditions hold (binomial distribution, gain at $0$); note that we use the concept of stochastic dominance.
		\begin{description}
			\item[(Bin)]  If $x \geq 1$, then $(X_{t+1} \mid X_0,\ldots, X_t; X_{t} = x) \succeq \mathrm{Bin}(k,(1+\delta) x/k)$.
			\item[(0)]  $\Ew{ \min(X_{t+1}, D_0) \mid X_{t} = 0} \geq E_0$.
		\end{description}
		Let $T := \min\set{t \geq 0}{X_{t} \geq n}$.

		Then, \textbf{if $\delta \leq 1$},
		\begin{align*}
		\Ew{T} 
		& \leq \frac{4D_0 }{0.4088 E_0} + \frac{15}{1-\gamma_0} D_0 \ln(2 D_0)+ 2.5 \log_2(n) \lceil 3 / \delta \rceil.
		\end{align*}
		
		In particular, when $\gamma_0$ is bounded away from $1$ by a constant, then $E[T] = \bigO{\frac{1}{E_0\delta} + \frac{\log (n)}{\delta}}$, where the asymptotic notation refers to $n$ tending to infinity and where $\delta=\delta(n)$ may be a function of $n$.
		Furthermore, if $n \gt 100/\delta$, then we also have that once the process has reached state of at least $100/\delta$, the probability to ever return to a state of at most $50/\delta$ is at most $0.5912$.

		\textbf{If $\delta \gt 1$}, then we have
		\begin{align*}
		\Ew{T}
		& \leq \frac{128}{0.78 E_0} + 2.6 \log_{1+\delta}(n) + 81\\
		& = \bigO{\frac{1}{E_0} + \frac{\log (n)}{\log(\delta)}}.
		\end{align*}
		In addition, once the process has reached state $32$ or higher, the probability to ever return to a state lower than $32$ is at most $\tfrac{1}{e(e-1)} \lt 0.22$.
	
\end{theoremEN}

	Note that this drift theorem is essentially restricted to processes based on the binomial distribution. For many applications this restriction is satisfied, particularly for the level-based theorem introduced in \cite{DBLP:conf/gecco/Lehre11} and refined in \cite{DBLP:journals/algorithmica/DangL16,DBLP:journals/tec/CorusDEL18}. We now discuss the currently strongest version in terms of the asymptotics in $\delta$, given in \cite{DBLP:journals/algorithmica/DoerrK21} as a consequence to the multiplicative up-drift theorem.

	The general setup of level-based theorems is as follows. There is a ground set $\mathcal{X}$, which in typical applications is the search space of an optimization problem. On this ground set, a Markov chain $(P_t)$ induced by a population-based EA is defined. We consider populations of fixed size $\lambda$, which may contain elements several times (multi-sets). We write $\mathcal{X}^\lambda$ to denote the set of all such populations. We only consider Markov chains where each element of the next population is sampled independently with repetition. That is, for each population $P \in \mathcal{X}^\lambda$, there is a distribution $D(P)$ on $\mathcal{X}$ such that given $P_t$, the next population $P_{t+1}$ consists of $\lambda$ elements of $\mathcal{X}$, each chosen independently according to the distribution $D(P_t)$. As all our results hold for any initial population $P_0$, we do not make any assumptions on $P_0$.

	In the level-based setting, we assume that there is a partition of $\mathcal{X}$ into levels $A_1, \dots, A_m$ (leading to the name of a \emph{level-based} theorem). Based on information in particular on how individuals in higher levels are generated, we aim for an upper bound on the first time such that the population contains an element of the highest level $A_m$.

\begin{theoremEN}[label=thm:advancedDrift:LevelBasedTheorem]{Level-Based Theorem}{}
	
		Consider a population-based process as described above.

		Let $(A_1,\ldots,A_m)$ be a partition of $\mathcal{X}$. Let $A_{\ge j} := \bigcup_{i=j}^m A_i$ for all $j \in [1..m]$. Let $z_1,\ldots,z_{m-1},\delta \in (0,1]$, and let $\gamma_0 \in (0,\frac{1}{1+\delta}]$  with $\gamma_0 \lambda \in \integers$. Let $D_0 = \min\{\lceil 100/\delta \rceil,\gamma_0 \lambda\}$ and $c_1 = 56 \, 000$. Let
		$$
		t_0 = \frac{7000}{\delta} \left(m + \frac{1}{1-\gamma_0} \sum_{j=1}^{m-1} \log^0_2\left(\frac{2\gamma_0\lambda}{1+\frac{z_j \lambda}{D_0}}\right) + \frac{1}{\lambda} \sum_{j=1}^{m-1}\frac{1}{z_j} \right),
		$$
		where $\log^0_2(x) := \max(0,\log_2(x))$ for all $x \in \realnum_+$. Assume that, for any population $P \in \mathcal{X}^\lambda$, the following three conditions are satisfied (drift, zero condition, population size).
		\begin{description}
			\item[(D)] For each level $j \in [1..m-2]$ and all $\gamma \in (0,\gamma_0]$, if $|P \cap A_{\geq j}| \geq \gamma_0 \lambda {/4}$ and $|P \cap A_{\geq j+1}| \geq \gamma \lambda$, then
				$$
				\operatorname{Pr}_{y \sim D(P)} \left[y \in A_{\geq j+1}\right] \geq (1+\delta)\gamma.
				$$ 
			\item[(0)] For each level $j \in [1..m-1]$, if $|P \cap A_{\geq j}| \geq \gamma_0 \lambda {/4}$, then
				$$
				\operatorname{Pr}_{y \sim D(P)} \left[y \in A_{\geq j+1}\right] \geq z_j.
				$$
			\item[(PS)] The population size $\lambda$ satisfies
				$$
				\lambda \geq {\frac{256}{\gamma_0 \delta} \ln \left(8 t_0 \right)}.
				$$
		\end{description}
		Then $T := \min \set{\lambda t}{P_t \cap A_m \neq \emptyset}$ satisfies
		\begin{align*}
		\Ew{T}  \leq 8\lambda t_0 = c_1 \frac{\lambda}{\delta} \left(m + \frac{1}{1-\gamma_0} \sum_{j=1}^{m-1} \log^0_2\left(\frac{2\gamma_0\lambda}{1+\frac{z_j \lambda}{D_0}}\right) + \frac 1 \lambda \sum_{j=1}^{m-1}\frac{1}{z_j} \right).
		\end{align*}
	
\end{theoremEN}

	Note that, with $z^* = \min_{j \in [1..m-1]} z_j$ and $\gamma_0$ a constant, \referDefined{(PS)} in the previous theorem is satisfied for some $\lambda$ with
	$$
	\lambda = O\left(\frac{1}{\delta}\log\left( \frac{m}{\delta z^*} \right)\right)
	$$
	as well as for all larger $\lambda$.

\subsection{Wormald's Method}
\label{subsec:advancedDriftTheorems:wormald}

	A very different approach to understanding random processes via their step-wise changes is given by Wormald \cite{Wormald:j:99}, tracking the processes via solutions of a system of differential equations. We briefly state a version of this theorem here.

	Consider a stochastic process $(Y^{(t)})_{t\geq 0}$, where each random variable $Y^{(t)}$ takes value in some set $S$. We use $H_t$ to denote a history of the process up to time $t$, i.e.~$H_t=(Y^{(0)},\ldots, Y^{(t)})$. And  $S^+$ denotes the set of all sequences $(Y^{(0)},\ldots, Y^{(t)})$ such that $Y^{(t)}\in S$. 

	We say that a function $f\colon \realnum^k \rightarrow \realnum$ satisfies a Lipschitz condition on $D \subseteq \realnum^k$ if there is $L \gt 0$ such that, for all $u=(u_1,\ldots,u_k), v=(v_1,\ldots,v_k)\in D$,  $$|f(u)-f(v)| \leq L \max_{1\leq i\leq k} |u_i-v_i|.$$

\begin{theoremEN}[label=thm:advancedDrift:Wormald]{Wormald's Method}{}
	
	For some $a \in \natnum$, let $(Y^{(t)}_i)_{1 \leq i \leq a, t \geq 0}$ be a stochastic process, such that there is $C \in \realnum_+$ so that for all $m \in \natnum_+$, $|Y_i^{(t)}| \lt m$ for all $H_t\in S^{+}$. Let $D$ be some bounded connected open set containing the closure of $$\Big\{(0,z_1,\ldots,z_a) \Big| \Pr{Y_i^{(0)}=z_im, 1\leq i \leq a} \neq 0 \mbox{ for some } m\Big\}.$$

	Assume the following three conditions hold, where for each $1 \leq i \leq a$ function $f_i\colon \realnum_{+} \times \realnum^{a} \rightarrow \realnum$ is continuous, and satisfies a Lipschitz condition on $D$ with the same Lipschitz constant $L$ for all $i$ (drift, boundedness).
	\begin{description}
		\item[(D)] $\Ew{Y_i^{(t+1)}-Y_i^{(t)} \mid H_t} = f_i(t/m,Y^{(t)}_1/m,\ldots,Y^{(t)}_a/m)$.
		\item[(B)] For all $t \geq 0$, $\max_{1\leq i\leq a} |Y_i^{(t+1)}-Y_i^{(t)}| \leq 1$.
	\end{description}
	Then the following are true.
	\begin{enumerate}
			\item For any $(0,\hat{z}_1,\ldots,\hat{z}_a)\in D$, the system of differential equations $$\frac{dz_i}{dx}=f_i(x,z_1,\ldots,z_a),\; i=1,\ldots,a$$ has a unique solution in $D$ for $z_i:\realnum \rightarrow \realnum$ passing through $z_i(0)=\hat{z}_i, 1\leq i\leq a$, and which extends to points arbitrary close to the boundary of $D$;
			\item Let $\lambda = \lambda(m) = o(1)$. For some constant $C > 0$, with probability $1-\bigO{\frac{1}{\lambda}\exp(-m\lambda^3)}$, 
		$$Y_i^{(t)}=mz_i(t/m)+\bigO{\lambda m}$$
		uniformly for $0\leq t\leq \sigma m$ and for each i, where $z_i(x)$ is the solution in given above with $\hat{z}_i=\frac{1}{m}Y_i^{(0)}$, and $\sigma=\sigma(m)$ is the supremum of those $x$ to which the solution can be extended before reaching within $l^{\infty}$-distance $C\lambda$ of the boundary of $D$.
	\end{enumerate}
	
\end{theoremEN}

	A few analyzes of randomized search heuristics with Wormald's theorem exist \cite{DBLP:conf/foga/LenglerS15,DBLP:conf/gecco/0001KM17,DBLP:journals/ec/Heredia18}.
	It has the advantage that it allows to track multiple interacting random variables (which, for other drift theorems, would have to be combined to a single potential). On the other hand, it requires solving a differential equation (well-known to be not an easy task) and the conclusion is typically deteriorating over time, since the variance is not averaged out but accumulates over time.

	Note that there are also theorems closer to the classic drift theorems for tracking multiple random variables in restricted settings \cite{DBLP:journals/tcs/Rowe18,DBLP:conf/ppsn/JanettL22}.

\clearpage

\section{No Going Back: The Fitness Level Method (FLM)}
\label{sec:fitnessLevelMethod}

	Some processes $(X_t)_{t \in \natnum}$ are \emph{monotone}, that is, we have $\forall t: X_t \leq X_{t+1}$. Monotone processes occur frequently in the analysis of heuristic optimization, since the best fitness found so far is a typical process considered. For some such processes, simpler (and sometimes stronger) analyses are possible than with drift theorems allowing for non-monotone processes (see \buildRef{thm:fitnessLevelMethod:pullProtocoll} for an example).

	Wegener \cite{DBLP:conf/icalp/Wegener01} proposed the following method, called the \emph{fitness level method} (FLM). We partition the search space into a number $m$ of sections (``levels'') in a linear fashion, so that all elements of later levels have better fitness than all elements of earlier levels. For the algorithm to be analyzed we regard the best-so-far individual and the level it is in. Since the best-so-far individual can never move to lower levels, it will visit each level at most once (possibly staying there for some time). Suppose we can show that, for any level $i \lt m$ which the algorithm is currently in, the probability to leave this level is at least $p_i$. Then, bounding the expected waiting for leaving a level~$i$ by $1/p_i$ (geometric distribution) and pessimistically assuming that we visit (and thus have to leave) each level $i \lt m$ before reaching the target level $m$, we can derive an upper bound for the optimization time of 
	$$
	\sum_{i=1}^{m-1} \frac{1}{p_i}.
	$$
	The fitness level method allows for simple and intuitive proofs and has therefore frequently been applied. Variations of it come with tail bounds~\cite{DBLP:journals/ipl/Witt14}, work for parallel EAs~\cite{DBLP:journals/ec/LassigS14} or regard populations \cite{DBLP:journals/ec/Witt06}. A similar analysis in levels can be made for non-elitist EAs, but here it is crucially possible (and sometimes not unlikely) to lose a level. See \buildRef{thm:advancedDrift:LevelBasedTheorem} for a corresponding theorem along a discussion.

	We state the fitness level method (FLM) formally as follows.

\begin{theoremEN}[label=thm:fitnessLevelMethod:upper]{Fitness Level Method (FLM)}{}
	
		Let $(X_t)_{t \in \natnum}$ be a monotone process on $[m]$. For all $i \in [m-1]$, let $p_i$ be a lower bound on the probability of a state change of $(X_t)_{t \in \natnum}$, conditional on being in state~$i$, formally: for all $t$ with $\Pr{X_t=i} \gt 0$, 
		$$
		\Pr{X_{t+1} \gt i \mid X_0,\ldots, X_t, X_t=i} \leq p_i.
		$$
		Let $T$ be the random variable describing the first time $t$ such that $X_t = m$. Then
		$$
		\Ew{T} \leq \sum_{i=1}^{m-1} \frac{1}{p_i}.
		$$

\end{theoremEN}
\begin{proof}
		For all $i \in [m-1]$, we let $S_i = \set{t \in \natnum}{X_t = i}$. Since $(X_t)_{t \in \natnum}$ is monotone, each $S_i$ is a discrete interval. Calling the leaving of $S_i$ a \emph{success event}, we can use \buildRef{thm:classicDrift:geometricDistribution} to see that the expected size of each $S_i$ is at most $1/p_i$. Since $T = \sum_{i=1}^{m-1}|S_i|$, the theorem follows.
	\end{proof}

\begin{noteEN}{}{}
	\textbf{Note.} The main strength of the fitness level method over drift theorems is that the chance to leave a level $i$, $p_i$, can depend arbitrarily on $i$. In contrast, in \buildRef{thm:advancedDrift:variableDrift}, one of the most general drift theorems, the drift has to depend monotonically on the state of the process.

	Conversely, the main strength of drift theorems over the fitness level method is that the \emph{process} is allowed to be non-monotone, so that we can work with other processes than those based on fitness.
\end{noteEN}

	The following example shows a toy application of the fitness level method where drift theorems are not easily applicable. The setting is borrowed from the area of \emph{rumor spreading}, see, for example, \cite{DBLP:conf/analco/DoerrK14} and also \cite{DBLP:journals/algorithmica/DoerrFFFKS19} for an application in the area of randomized search heuristics.

\begin{propositionEN}[label=thm:fitnessLevelMethod:pullProtocoll]{Application to Rumor Spreading}{}
	
		Let $n \in \natnum$. Suppose $n$ people each want to obtain a certain information, and suppose in iteration $0$ exactly one of them knows this information. In each iteration, one of the $n$ people is chosen uniformly at random and if this person does not know the information, the person will contact another person chosen uniformly at random. If this other person knows the information, then the calling person from now on also knows the information. Then it takes, in expectation, at most $2n(\ln(n-1)+1)$ iterations until all persons know the rumor.

\end{propositionEN}
\begin{proof}
		We let, for each $t \in \natnum$, $X_t$ be the number of persons who know the rumor after $t$ iterations. Then $(X_t)_{t \in \natnum}$ is a monotone process on $[n]$. Let some iteration $t$ be given. If, after $t$ iterations, exactly $i \in [n-1]$ persons know the rumor, then the probability that in the next iteration an uninformed person is chose to make a call is $(n-i)/n$. The probability that this person calls an informed person is independent of that probability and $i/(n-1)$. Thus, the chance $p_i$ to ``leave state $i$'' is 
		$$
		p_i = \frac{n-i}{n} \; \frac{i}{n-1}.
		$$
		By \buildRef{thm:fitnessLevelMethod:upper}, the total time until all people are informed is thus at most
		$$
		\sum_{i=1}^{n-1} \frac{1}{p_i} = \sum_{i=1}^{n-1} \frac{n(n-1)}{(n-i)i}.
		$$
		\inlineComment{For didactic reasons, we give two different ways of bounding this sum; the second uses calculus and gives the better bound.}
		A first way to bound the sum is by splitting it into two and using worst case estimates to simplify. Let $k = \lfloor n/2 \rfloor$; we have
		\begin{align*}
		\sum_{i=1}^{n-1} \frac{n(n-1)}{(n-i)i} & = \sum_{i=1}^{k} \frac{n(n-1)}{(n-i)i} + \sum_{i=k+1}^{n-1} \frac{n(n-1)}{(n-i)i}\\
		& \leq \sum_{i=1}^{k} \frac{n(n-1)}{(n-k)i} + \sum_{i=k+1}^{n-1} \frac{n(n-1)}{(n-i)k}\\
		& = \sum_{i=1}^{k} \frac{n(n-1)}{(n-k)i} + \sum_{i=1}^{n-k-1} \frac{n(n-1)}{i \; k}\\
		& = \frac{n(n-1)}{n-k} \sum_{i=1}^{k} \frac{1}{i} + \frac{n(n-1)}{k}\sum_{i=1}^{n-k-1} \frac{1}{i}\\
		& \leq \frac{n(n-1)}{n-k} \left(\ln(k)+1\right) + \frac{n(n-1)}{k}\left(\ln(n-k-1)+1\right).
		\end{align*}
		The last inequality uses \buildRef{lem:notation:boundHarmonicSum}.
		By bounding $1/(n-k) \leq 2/n$ and $1/k \leq 2/(n-1)$ we can further bound the term by
		$$
		2(n-1)\left(\ln(k)+1\right)+2n\left(\ln(n-k-1)+1\right)  \leq 4n\left(\ln(n)+1\right).
		$$
		We can get a tighter bound as follows by turning to calculus. The function 
		$$
		f\colon [1,n-1] \rightarrow \realnum, x \mapsto \frac{1}{(n-x)x}
		$$
		has a minimum at $n/2$, is before that monotone decreasing and afterwards monotone increasing. Thus, we can use \buildRef{lem:notation:sumByIntegral} twice, once on the interval $[1,\lfloor n/2 \rfloor]$ and once on $[\lceil n/2 \rceil,n-1]$ to bound
		$$
		\sum_{i=1}^{n-1} \frac{n(n-1)}{(n-i)i} \leq n + n + n(n-1) \int_1^{n-1} f(x) \mathrm{d}x.
		$$
		Note that the summands before the integral are the first and the last summand of the large sum (which are not covered by the cited lemma).
		The indefinite integral over $f$ is given by the function
		$$
		x \mapsto \frac{\ln(x)-\ln(n-x)}{n},
		$$
		which can be seen by taking the derivative of that function. Using the integral bounds, we arrive at
		$$
		\sum_{i=1}^{n-1} \frac{n(n-1)}{(n-i)i} \leq 2n + (n-1)\left(\ln(n-1)-\ln(1) - \ln(1)+\ln(n-1)\right) = 2(n-1)\ln(n-1) + 2n
		$$
		as desired.
	\end{proof}

	Note that due to the drift being strongest in the middle of the state space, no other drift theorem is directly available without losing asymptotically: one could apply the additive drift theorem to obtain a bound of $O(n^2)$ by using that the drift is $\Omega(1/n)$. Another choice is to us an argument in phases: since the process is monotone, one can analyze the time until reaching $n/2$ separately from the remainder. This would then potentially allow to use some version of multiplicative up-drift (for the first phase) and multiplicative drift (for the second phase). However, this would lead to an unnecessarily complicated analysis.

	While very effective for proving upper bounds, it seems much harder to use fitness level arguments to prove lower bounds. The first to devise a lower bound method based on fitness levels that gives competitive bounds was Sudholt \cite{DBLP:journals/tec/Sudholt13}. Next we see a lower bound from \cite{DBLP:conf/gecco/DoerrK21}.

\begin{theoremEN}[label=thm:fitnessLevelMethodVisit:lower]{Fitness Level Method with Visit Probabilities, Lower Bound}{}
	
		Let $(X_t)_t$ be a monotone process on $[m]$. For all $i \in [m-1]$, let $p_i$ be an upper bound on the probability of a state change of $(X_t)_t$, conditional on being in state~$i$. Furthermore, let $v_i$ be a lower bound on the probability of there being a $t$ such that $X_t = i$ (the \emph{visit probability} of level~$i$). Then the expected time for $(X_t)_t$ to reach the state $m$ is 
		$$
		\Ew{T} \geq \sum_{i=1}^{m-1} \frac{v_i}{p_i}.
		$$

\end{theoremEN}
\begin{proof}
		We proceed as in the proof for \buildRef{thm:fitnessLevelMethod:upper}. For all $i \in [m-1]$, we let $S_i = \set{t \in \natnum}{X_t = i}$. With probability at most $1-v_i$ we have that $S_i = \emptyset$. Again using \buildRef{thm:classicDrift:geometricDistribution}, we see that the expected size of each non-empty $S_i$ is at least $1/p_i$. Since $T = \sum_{i=1}^{m-1}|S_i|$, the theorem follows with linearity of expectation.
	\end{proof}

	A corresponding upper bound \cite{DBLP:conf/gecco/DoerrK21} follows with analogous arguments and shows the tightness of the approach. Note that the bounds required on $p_i$ and $v_i$ are naturally from the other side.

\begin{theoremEN}[label=thm:fitnessLevelMethodVisit:upper]{Fitness Level Method with Visit Probabilities, Upper Bound}{}
	
		Let $(X_t)_t$ be a monotone process on $[m]$. For all $i \in [m-1]$, let $p_i$ be an lower bound on the probability of a state change of $(X_t)_t$, conditional on being in state~$i$. Furthermore, let $v_i$ be an upper bound on the probability of there being a $t$ such that $X_t = i$. Then the expected time for $(X_t)_t$ to reach the state $m$ is 
		$$
		\Ew{T} \leq \sum_{i=1}^{m-1} \frac{v_i}{p_i}.
		$$

\end{theoremEN}
\begin{proof}
		Analogous to the proof of \buildRef{thm:fitnessLevelMethodVisit:lower}.
	\end{proof}

	In a typical application of the fitness level method, finding good estimates for the leaving probabilities is easy. It is more complicated to estimate the visit probabilities accurately, the following lemma from \cite{DBLP:conf/gecco/DoerrK21} offers an option.

\begin{lemmaEN}[label=thm:computeVisit]{Computing Visit Probabilities}{}
	
		Let $(X_t)_t$ be a monotone process on $[m]$. Further, suppose that $(X_t)_t$ reaches state $m$ after a finite time with probability $1$. Let $i \lt m$ be given. 		
		Suppose there is $v_i$ such that, for all $t \in \natnum$ with $\Pr{X_{t+1} \geq i \gt X_t} \gt 0$,
		$$
		\Pr{X_{t+1} = i \mid X_0,\ldots, X_t; X_{t+1} \geq i \gt X_t} \geq v_i,
		$$
		and 
		$$
		\Pr{X_0 = i \mid X_0 \geq i} \geq v_i.
		$$
		Then $v_i$ is a lower bound for visiting level~$i$ as required by \buildRef{thm:fitnessLevelMethodVisit:lower}.
	
\end{lemmaEN}

	An analogous bound for upper bounds on visit probabilities also holds.

\subsection{Applications}

	As a first application of these methods we now determine a lower bound for the coupon collector problem (see \buildRef{thm:classicDrift:couponCollector}).

\begin{theoremEN}[label=thm:fitnessLevelMethod:couponCollectorLowerBoundFLM]{Coupon Collector, Lower Bound via Fitness Levels}{}
	
		Suppose we want to collect at least one of each kind of $n \in \natnum_{\geq 1}$ coupons. Each round, we are given one coupon chosen uniformly at random from the~$n$ kinds. Then, in expectation, we have to collect for at least $n(1+ \ln n)$ iterations.

\end{theoremEN}
\begin{proof}
		Let $X_t$ be the number of coupons after $t$ iterations. Note that this process is monotone and, since it gains at most one in any iteration, visits all elements of $[n-1]$ before reaching the target of $n$. The probability of making progress (of $1$) with coupon $t+1$ is $p_i= (n-i)/n$, since $n-i$ coupons are missing and each has a probability of $1/n$, and all these events are disjoint. An application of both \buildRef{thm:fitnessLevelMethodVisit:lower} and \buildRef{thm:fitnessLevelMethodVisit:upper} gives an exact value of the expected time to find all $n$ coupons of
		$$
		\sum_{i=0}^{n-1} \frac{v_i}{p_i} = \sum_{i=0}^{n-1} \frac{n}{n-i} = nH_{n},
		$$
		Where, for each $m \in \natnum_+$, we use $H_m$ to denote the $m$th harmonic number.
	\end{proof}

	As a second  application  we now determine the precise run time of the \OneOneEA on \textsc{LeadingOnes} via the two fitness level theorems. This optimization time was first established in \cite{DBLP:conf/ppsn/BottcherDN10}.

\begin{theoremEN}[label=thm:fitnessLevelMethod:leadingOnes]{Run Time of \OneOneEA on \textsc{LeadingOnes}}{}

		Consider the \OneOneEA optimizing \textsc{LeadingOnes} with mutation rate $p$. Let $T$ be the (random) time for the \OneOneEA to find the optimum. Then
		$$
		\Ew{T} = \frac{1}{2} \sum_{i=0}^{n-1} \frac{1}{(1-p)^i p}.
		$$

\end{theoremEN}
\begin{proof}
		We want to apply \buildRef{thm:fitnessLevelMethodVisit:upper} and \buildRef{thm:fitnessLevelMethodVisit:lower} simultaneously.  For all $t \in \natnum$, we let $X_t$ be the \textsc{LeadingOnes}-value of the individual which the \OneOneEA has found after $t$ iterations. Now we need a precise result for the probability to leave a level and for the probability to visit a level.

		First, we consider the probability $p_i$ to leave a given level $i \lt n$. Suppose the algorithm has a current search point in level $i$, so it has $i$ leading $1$s and then a $0$. The algorithm leaves level $A_i$ now if and only if it flips the first $0$ of the bit string (probability of $p$) and no previous bits (probability $(1-p)^i$). Hence, $p_i = p(1-p)^i$.

		Next we consider the probability $v_i$ to visit a level~$i$. We claim that it is exactly $1/2$, following reasoning given in several places before~\cite{DBLP:journals/tcs/DrosteJW02,DBLP:journals/tec/Sudholt13}. We want to use \buildRef{thm:computeVisit} and its analogue for upper bounds. Let $i$ be given. For the initial search point, if it is at least on level~$i$ (the condition considered by the lemma), the individual is on level~$i$ if and only if the $i+1$st bit is a $0$, so exactly with probability $1/2$ as desired for both bounds. Before an individual with at least $i$ leading $1$s is created, the bit at position $i+1$ remains uniformly random (this can be seen by induction: it is uniform at the beginning and does not experience any bias in any iteration while no individual with at least~$i$ leading $1$s is created). Once such an individual is created, if the bit at position $i+1$ is $1$, the level~$i$ is skipped, otherwise it is visited. Thus, the algorithm skips level~$i$ with probability exactly $1/2$, giving $v_i = 1/2$. With these exact values for the $p_i$ and $v_i$, \buildRef{thm:fitnessLevelMethodVisit:upper} and~\buildRef{thm:fitnessLevelMethodVisit:lower} immediately yield the claim.	
	\end{proof}

	By computing the geometric series in \buildRef{thm:fitnessLevelMethod:leadingOnes}, we obtain as a (well-known) corollary that the \OneOneEA with the classic mutation rate $p = 1/n$  optimizes \textsc{LeadingOnes} in an expected run time of $n^2\frac{e-1}{2}(1\pm o(1))$.

\clearpage

\section{A Different Perspective: Fixed Budget Optimization}
\label{sec:fixedBudget}

	In the previous chapters we have seen many theorems regarding the \emph{first hitting time} of a process. This answers the question: ``How much time do I have to invest until a desired outcome?'' Sometimes we want to answer a different question: ``I have a \emph{fixed budget} $t_0$ of time available; what performance can I expect?'' 
	Furthermore, fixed-budget results that hold with high probability are crucial for the analysis of algorithm configurators~\cite{DBLP:conf/gecco/HallOS19}. These configurators test different algorithms for fixed budgets in order to make statements about their appropriateness in a given setting.

	In this chapter we want to discuss general tools for fixed-budget analyses. We still want to use knowledge about step-wise changes and translate them into the global view, just as for the drift theorems for first hitting times.

	We start by analyzing of the most basic setting in \buildRef{sec:fixedBudget:additive} and generalize it in \buildRef{sec:fixedBudget:variable}. We show sample results derived with these methods in \buildRef{subsec:fixedBudget:applications}.

\subsection{The Additive Case}
\label{sec:fixedBudget:additive}

	We start with the simple case of additive drift. If we expect to go down by $\delta$ in each iteration, then, after $t$ iterations, we expect to be down $t \delta$, as would be the case for a \emph{completely deterministic} process. Here the proof is simple and instructive.

\begin{theoremEN}[label=thm:fixedBudget:additive]{Additive Fixed-Budget Drift}{}
	
		Let $(X_t)_{t \in \natnum}$, be a stochastic process on $\realnum$. Suppose there is a $\delta \in \realnum_{+}$ so that we have the drift condition
		\begin{description}
			\item[(D)] $\Ew{X_t - X_{t+1} \mid X_0,\ldots, X_t} \geq \delta$.
		\end{description}
		Thus, the drift condition is equivalent to
		\begin{description}
			\item[(D')] $\Ew{X_{t+1} \mid X_0,\ldots, X_t} \leq X_t - \delta$.
		\end{description}
		Then, for all $t \geq 0$,
		$$
			\Ew{X_t \mid X_0} \leq X_0 - t \delta.
		$$

\end{theoremEN}
\begin{proof}
		We prove the theorem by induction on $t$, with a trivial induction basis. Suppose now the statement is true for some $t \geq 0$ (IH). Using the law of total expectation (LTE), we have
		\begin{align*}
			\Ew{X_{t+1} \mid X_0} 
				& \eqnComment{(LTE)}{=} \Ew{ \Ew{X_{t+1} \mid X_0,\ldots, X_t} \mid X_0}\\
				& \eqnComment{(D')}{\leq} \Ew{ X_t - \delta \mid X_0}\\
				& \eqnComment{(IH)}{\leq} X_0 - (t+1) \delta.
		\end{align*}
		This concludes the induction.
	\end{proof}

	Note that this version does not take into account that drift might only hold before a target has been reached. We refer to this setting as \emph{unlimited time}. The next theorem considers a potential end point. Note that the proof follows the proof of Theorem~1 in \cite{DBLP:series/ncs/Lengler20}.

\begin{theoremEN}[label=thm:fixedBudget:additiveTimeLimited]{Additive Fixed-Budget Drift, Limited Time}{}
	
		Let $X_t$, $t \geq 0$, be a stochastic process on $\realnum_{\geq 0}$ and let $T$ be any random variable on $\natnum$. Suppose that there is a $\delta \in \realnum_{+}$ so that we have the drift condition
		\begin{description}
			\item[(D)] $\Ew{X_t - X_{t+1} \mid t \lt T} \geq \delta$.
		\end{description}
		Then, for all $t \geq 0$,
		$$
			\Ew{X_t \mid X_0} \leq X_0 - t\delta \Pr{t \leq T}.
		$$

\end{theoremEN}
\begin{proof}
		We prove the theorem by induction on $t$, with a trivial induction basis for $t=0$. Suppose now the statement is true for some $t \geq 0$ (IH). Then, using the law of total expectation (LTE),
		\begin{align*}
			\Ew{X_{t+1} \mid X_0} 
				& = \Ew{ X_{t+1} - X_t + X_t \mid X_0}\\
				& = \Ew{ X_{t+1} - X_t \mid X_0} + \Ew{ X_t \mid X_0}\\
				& \eqnComment{(IH)}{\leq} \Ew{ X_{t+1} - X_t \mid X_0} + X_0 - t\delta \Pr{t \leq T}\\
				& \eqnComment{(LTE)}{=} \Ew{ X_{t+1} - X_t \mid X_0, t \lt T}\Pr{t \lt T} + \Ew{ X_{t+1} - X_t \mid X_0, t \geq T}\Pr{t \geq T}\\
				& \;\;\;\;\;\;+ X_0 - t\delta \Pr{t \leq T}\\
				& \leq -\delta\Pr{t \lt T} + X_0 - t\delta \Pr{t \leq T}\\
				& \leq -\delta\Pr{t+1 \leq T} + X_0 - t\delta \Pr{t+1 \leq T}\\
				& = X_0 - (t+1)\delta \Pr{t+1 \leq T}.
		\end{align*}
		This concludes the induction.
	\end{proof}

\subsection{Variable Fixed Budget Drift}
\label{sec:fixedBudget:variable}

	For the rest of this chapter, we want to generalize \buildRef{thm:fixedBudget:additive} to \emph{state}-varying drift. Suppose that, for some function $h$, in state $x$ we observe a drift of $h(x)$. In order to understand what kind of result to expect in this context, we consider a completely deterministic process starting in $x_0 \in \realnum$ and progressing down by $h(x)$ when in state $x$. Then, after one step, the process is in $x_0 - h(x_0)$, after two steps in $x_0 - h(x_0) - h(x_0-h(x_0))$ and so on. We write, for all $x$, $\tilde{h}(x) = x - h(x)$. Thus we can write the sequence of states of the process as $x_0, \tilde{h}(x_0), \tilde{h}(\tilde{h}(x_0))$  and so on. We write $\tilde{h}^t$ for the $t$-fold application of $\tilde{h}$, so after $t$ steps of the process the state is $\tilde{h}^t$. Thus, we want a theorem that shows that we get a similar expected value for a probabilistic process.

	The main question is now what we need to assume about $h$ to get a behavior similar to the deterministic process. Consider the following monotone process on $\{0,1,2\}$: $X_0$ is $2$ and the process moves to one of $\{0,1\}$ uniformly. State $0$ is the target state, from state $1$ there is only a very small probability to progress to $0$ (say $0.1$). Then it is better to stay in State $2$ instead of being trapped in State $1$. Here the drift is $1.5$ in State $2$ and only $0.1$ in State $1$. Thus, the \emph{expected next state} for State $2$ is $0.5$, which is less than the expected next state for State $1$, which is $0.9$! Intuitively, greedily going forward is a bad idea, if given the choice between States $1$ and $2$ one should choose (non-greedily) State $2$. It turns out that forbidding this kind of situation, formalized in the next definition, leads to a viable generalization of \buildRef{thm:fixedBudget:additive}.

\begin{definitionEN}{Greed-Admitting Functions}{}
	We say that a drift function $h\colon S \to \realnum_{\gt 0}$ is \emph{greed-admitting} if $\mathrm{id} - h$ (the function $x \mapsto x - h(x)$) is monotone non-decreasing.
\end{definitionEN}

	Intuitively, this formalizes the idea that being closer to the goal is always better (``greed is good''). The process described before the definition is, in a sense, badly designed: State $1$ is worse than State $2$, so it should not have a smaller value.

	We now give two different versions of fixed-budget drift theorems. The first considers \emph{unlimited time}, a very strong requirement, leading to a strong conclusion.

\begin{theoremEN}[label=thm:fixedBudget:variableUnlimitedTime]{Variable Fixed-Budget Drift, Unlimited Time}{}
	
		Let $X_t$, $t \geq 0$, be a stochastic process on $S \subseteq \realnum$, where $0 = \min S$. Let $h\colon S \to \realnum_{\geq 0}$ be a \emph{twice differentiable, convex and greed-admitting} function such that $\tilde{h}'(0) \in \; ]0,1]$ and we have the drift condition
		\begin{description}
			\item[(D-ut)] $\Ew{X_t - X_{t+1} \mid X_0,\ldots, X_t} \geq h(X_t)$.
		\end{description}
		Define $\tilde{h}(x) = x - h(x)$. Thus, the drift condition is equivalent to
		\begin{description}
			\item[(D-ut')] $\Ew{X_{t+1} \mid X_0,\ldots, X_t} \leq \tilde{h}(X_t)$.
		\end{description}
		Then, for all $t \geq 0$,
		$$
			\Ew{X_t \mid X_0} \leq \tilde{h}^t(X_0)
		$$
		and, in particular,
		$$
			\Ew{X_t} \leq \tilde{h}^t(\Ew{X_0}).
		$$
	
\end{theoremEN}

	Crucial for this theorem is that the drift condition is \emph{unlimited time}, by which we mean that the drift condition has to hold for all times $t$, not just (which is the typical case in the literature for drift theorems) those before the optimum is hit. This theorem is applicable if there is no optimum (and the optimization progresses indefinitely) of if the drift is $0$ in the optimum. In order to bypass these limitations we also give a variant which allows for \emph{limited time} drift, where the drift condition only needs to hold before the optimum is hit; however, in this case we pick up an additional error term in the result, derived from the possibility of hitting the optimum within the allowed time budget of $t$. Thus, in order to apply this theorem, one will typically need concentrations bounds for the time to hit the optimum.

	A special case of the previous theorem is given in~\cite{DBLP:conf/foga/LenglerS15}, where the drift is necessarily multiplicative. Note that in this case we can typically consider unlimited time, since after reaching the state $0$ the multiplicative drift holds vacuously.

	Now we give a version of the variable fixed-budget drift where the time is limited in the sense that the drift condition might no longer hold at some point in time.

\begin{theoremEN}[label=thm:fixedBudget:variableLimitedTime]{Variable Fixed-Budget Drift, Limited Time}{}
	
		Let $X_t$, $t \geq 0$, be a stochastic process on $S \subseteq \realnum$, where $0 = \min S$. Let $T = \min\set{ t \ge 0}{ X_t = 0 }$ and $h\colon S \to \realnum_{\geq 0}$ be a \emph{twice differentiable, convex and greed-admitting} function such that $\tilde{h}'(0) \in \; ]0,1]$ and we have, for all $t \lt T$, the drift condition
		\begin{description}
			\item[(D-lt)]  $\Ew{X_t - X_{t+1} \mid X_0,\ldots, X_t} \geq h(X_t)$.
		\end{description}
		Define $\tilde{h}(x) = x - h(x)$. Thus, the drift condition is equivalent to
		\begin{description}
			\item[(D-lt')] $\Ew{X_{t+1} \mid X_0,\ldots, X_t} \leq \tilde{h}(X_t)$.
		\end{description}
		Then, for all $t \geq 0$,
		$$
			\Ew{X_t \mid X_0} \leq \tilde{h}^t(X_0) + \frac{\tilde{h}(0)}{\tilde{h}'(0)}
		$$
		and, in particular,
		$$
			\Ew{X_t} \leq \tilde{h}^t(\Ew{X_0}) - \frac{\tilde{h}(0)}{\tilde{h}'(0)} \cdot \Pr{t \geq T \mid X_0}.
		$$
	
\end{theoremEN}

	For both these theorems, the drift function bounding the drift has to be convex and \emph{greed-admitting}, which intuitively says that being closer to the goal is always better in terms of the expected state after an additional iteration, while search points closer to the goal are required to have weaker drift. These conditions are fulfilled in many sample applications.

	In order to interpret the conclusions of the last two theorems properly, we need to estimate the term $\tilde{h}^t$. With the following theorem we give a general way of making this estimation.

\begin{theoremEN}[label=thm:fixedBudget:estimateTildeHwithIntegral]{Estimation of Iterated Functions}{}
	
		Let $h\colon \realnum \rightarrow \realnum_{+}$ be a monotone non-decreasing and integrable function. Let $\tilde{h} = \mathrm{id} - h$. Then, for all starting points $n$ and all target points $x \leq y$ and all time budgets $t$,
		$$
		\mbox{if } t \geq \int_{x}^{y} \frac{1}{h(z)}\mathrm{d}z \;\;\mbox{  then  }\;\; \tilde{h}^t(y) \leq x.
		$$
	
\end{theoremEN}

	We can specialize the previous theorem to the disc\emph{}rete case.

\begin{theoremEN}[label=thm:fixedBudget:estimateTildeHdiscrete]{Estimation of Iterated Functions, Sum Formula}{}
	
		Let $h\colon \natnum \rightarrow \realnum_{+}$ be a monotone non-decreasing function and let $\tilde{h} = \mathrm{id} - h$. Then, for all starting points $n \in \natnum$ and all target points $m \leq n$ and all time budgets $t$,
		$$
		\mbox{if } t \geq \sum_{i=m}^{n-1} \frac{1}{h(i)} \;\;\mbox{ then }\;\; \tilde{h}^t(n) \leq m.
		$$

\end{theoremEN}
\begin{proof}
		We apply \buildRef{thm:fixedBudget:estimateTildeHwithIntegral} to $h'\colon \realnum \rightarrow \realnum_{\gt 0}, x \mapsto h(\max(0,\lfloor x \rfloor))$ and use that, for all $i \in \natnum$,
		$$
			\frac{1}{h(i)} = \int_{i}^{i+1} \frac{1}{h'(z)}\mathrm{d}z.
		$$
	\end{proof}

\subsection{Applications to \textsc{OneMax} and \textsc{LeadingOnes}}
\label{subsec:fixedBudget:applications}

	In this section we show results from applications of \buildRef{thm:fixedBudget:variableUnlimitedTime} as given in \cite{DBLP:conf/ppsn/KotzingW20}. We consider the optimization of the \OneOneEA on \textsc{OneMax} and on \textsc{LeadingOnes} as examples. We start with \textsc{OneMax}, where we have multiplicative drift.

\begin{theoremEN}[label=thm:fixedBudget:oneMaxWithDirectDrift]{Fixed Budget for OneMax}{}
	
		Let $X_t$ be the number of $1$s which the \OneOneEA on $\OneMax$ has found after $t$ iterations of the algorithm. Then we have, for all $t$,
		$$
		\Ew{X_t} \geq 
		\begin{cases}
		\frac{n}{2} + \frac{t}{2\sqrt{e}}-O(1),				&\mbox{if }t = O(\sqrt{n});\\
		\frac{n}{2} + \frac{t}{2\sqrt{e}}(1-o(1)),			&\mbox{if }t = o(n).
		\end{cases}
		$$
		Furthermore, for all $t$, we have $\Ew{X_t} \geq n(1 - \exp(-t/(en))/2)$.
	
\end{theoremEN}

	For the \OneOneEA on \textsc{OneMax}, no concrete formula for a bound on the fitness value after $t$ iterations was known: The original work \cite{DBLP:conf/gecco/JansenZ12} could only handle RLS on \textsc{OneMax}, not the \OneOneEA. The multiplicative drift theorem of \cite{DBLP:conf/foga/LenglerS15} allows for deriving a lower bound of $n/2 + t/(2e)$ for $t = o(n)$, using a multiplicative drift constant of $(1-1/n)^n/n$. Since our drift theorem allows for variable drift, we can give the better bound of $n/2 + t/(2\sqrt{e}) - o(t)$ for the \OneOneEA on \textsc{OneMax} with $t = o(n)$. Note that \cite{DBLP:conf/foga/LenglerS15} also gives bounds for values of $t$ closer to the expected optimization time.

	Our second example shows the progress of the \OneOneEA on \textsc{LeadingOnes}, where we have additive drift. The result is summarized in the following theorem.

\begin{theoremEN}[label=thm:fixedBudget:leadingOnesWithDirectDrift]{Fixed Budget for LeadingOnes}{}
	
		Let $X_t$ be the number of leading $1$s which the \OneOneEA on \textsc{LeadingOnes} has found after $t$ iterations of the algorithm. We have, for all $t$,
		$$
		\Ew{X_t} \geq 
			\begin{cases}
			\frac{2t}{n} - O(1),									&\mbox{if }t=O(n^{3/2});\\
			\frac{2t}{n} \cdot (1 - o(1)),							&\mbox{if }t=o(n^2);\\
			n\ln(1+\frac{2t}{n^2}) - O(1),							&\mbox{if }t \leq \frac{e-1}{2} n^2 - n^{3/2}.
			\end{cases}
		$$
	
\end{theoremEN}

	For the \OneOneEA on \textsc{LeadingOnes} with a budget of $t = o(n^2)$ iterations, the paper \cite{DBLP:conf/gecco/JansenZ12} gives a lower bound of $2t/n - o(t/n)$ for the expected fitness after $t$ iterations, which are recovered with a simpler proof. The general theorems from this section also allow budgets closer to the expected optimization time, where we get a lower bound of $n\ln(1+2t/n^2) - O(1)$.

\subsection{Bibliographic Remarks}

	The setting of fixed-budget analysis was introduced to the analysis of randomized search heuristics by Jansen and Zarges~\cite{DBLP:conf/gecco/JansenZ12}, who derived fixed-budget results for the classical example functions \textsc{OneMax} and \textsc{LeadingOnes} by bounding the expected progress in each iteration. A different perspective was proposed by Doerr, Jansen, Witt and Zarges \cite{DBLP:conf/gecco/DoerrJWZ13}, who showed that fixed-budget statements can be derived from bounds on optimization times if these exhibit strong concentration. Lengler and Spooner~\cite{DBLP:conf/foga/LenglerS15} proposed a variant of multiplicative drift for fixed-budget results and the use of differential equations in the context of \textsc{OneMax} and general linear functions. Nallaperuma, Neumann and Sudholt~\cite{DBLP:journals/ec/NallaperumaNS17} applied fixed-budget theory to the analysis of evolutionary algorithms on the traveling salesman problem and Jansen and Zarges~\cite{DBLP:journals/tec/0001Z14} to artificial immune systems. The quality gains of optimal black-box algorithms on \textsc{OneMax} in a fixed-budget perspective were analyzed by Doerr, Doerr and Yang~\cite{DBLP:journals/tcs/DoerrDY20}. He, Jansen and Zarges \cite{DBLP:conf/gecco/00040Z19} consider the so-called unlimited budgets to estimate fitness values in particular for points of time larger than the expected optimization time. A survey by Jansen~\cite{DBLP:series/ncs/000120} summarizes the state of the art in the area of fixed-budget analysis.

	In contrast to the numerous drift theorems available for bounding the optimization time, there was no corresponding theorem for making a fixed-budget analysis apart from one for the multiplicative case given in~\cite{DBLP:conf/foga/LenglerS15}. This changed with~\cite{DBLP:conf/ppsn/KotzingW20}, introduced in this chapter, providing several such drift theorems.
	
	Note that a further fixed-budget drift theorem can be found in \cite{DBLP:conf/ppsn/KotzingW20}, where a detour of the computation of fixed budget results via first hitting times is made.

\clearpage

\section{Drift as an Average: A closer look on the conditioning of drift}
\label{sec:averagingDrift}

	Consider a deterministic process which starts at $100$ and goes down by $1$ in each iteration; we trivially see that the expected time until the process reaches $0$ is $100$. If the process goes down by $1$ or $1/2$, then a worst-case view would state an upper bound of $200$ until the process reaches $0$.

	Drift theory allows for an \emph{average case} view. In fact, the main strength of drift theory is that even the possibility of going \emph{away} from the target is incorporated, as long as \emph{in expectation} we have a bias towards the target. For example, if a process goes down by $10$ with probability $11/20$ and up by $10$ with probability $9/20$, then the progress is \emph{in expectation} also $1$, but not in the worst case. The additive drift theorem tells us that, also in this case, we expect to arrive at $0$ after $100$ steps. Thus, instead of a worst case bound, averaging different outcomes leads to a useful bound.

	In this section we investigate three questions.
	\begin{enumerate}
			\item How do drift theorems allow to exploit that drift is an average over a range of possibilities?
			\item Why do drift theorems in the literature condition on various different things?
			\item How do we account for insufficient drift after reaching the target?
	\end{enumerate}

\subsection{Drift as an Average}

	In \buildRef{thm:classicDrift:additiveDriftUpper} we have seen the standard (additive) drift condition to be
	\begin{description}
		\item[(D)] there is a $\delta > 0$ such that, for all $t \lt T$, it holds that $\Ew{X_t - X_{t+1} \mid X_0,\ldots,X_t} \geq \delta$.
	\end{description}
	In this section we want to take a closer look at the conditioning on $X_0,\ldots,X_t$. First we show that \emph{not} conditioning on anything leads to counterexamples.

\begin{beispielEN}[label=ex:averagingDrift:globalStop]{Global Averaging}{}
	Suppose we first flip a coin in secret. Then, if the coin shows tails, for all $i \in \natnum$ we let $X_i = 1$. If, on the other hand, the coin shows heads, we let $X_0 = 1$ and, for all $i \in \natnum$, we draw a uniformly random bit $B \in \{0,1\}$ and set $X_{i+1} = X_i - B$. For any $t$ with $t \lt T$ we now have $\Ew{X_t - X_{t+1}} \geq 1/4$ (we make a progress of exactly $1$ if the coin shows tails and the bit is $1$, and otherwise of $0$). Thus, the conclusion of the additive drift theorem states an upper bound of $4$ on the expected time to hit $0$, while, in fact, $\Ew{T}$ is infinite.
\end{beispielEN}

	The previous example shows that the conclusion of the additive drift theorem can be false while the drift condition holds, averaged over \emph{all} possible situations. What we \emph{can} do is average over \emph{all possible situations with the same history of $X_0, \ldots, X_t$}, as stated by the additive drift theorem. To illustrate this, we have the following example.

\begin{beispielEN}[label=ex:averagingDrift:localAveraging]{Local Averaging}{}
	Let us play a game where your goal is to draw a total number of $10$ red balls. In each iteration I randomly fill in secret a bag of balls of different colors, and draw a ball uniformly at random from that bag. Suppose I either fill the bag with $10$ balls, $9$ of which are blue and one is red, or with $100$ balls, where $99$ are blue and one is red. Suppose I choose either situation with equal probability of $1/2$. Then, on average, in any iteration your probability to pick a red ball is $11/200$. Thus, you arrive at a value of $10$ drawn red balls after an expected number of $2000/11$ iterations.
\end{beispielEN}

	Note that in this example there seem to be two different possible situations in each iteration with different drift. One way to address this is to bound drift by the smaller of the two drift values; but the drift theorem allows for averaging the drift of the different situations, since we are given the probabilities of the two values.

	We cannot average over global decisions that we can learn about from the history, see \buildRef{ex:averagingDrift:globalStop}. If we cannot learn about the global decisions from the history, we can use the \emph{principle of deferred decision} to model the random decision as a decision in that given iteration, as illustrated by the following example.

\begin{beispielEN}{Deferred Decisions}{}
	Let us again play a game where your goal is to draw a total number of $10$ red balls. This time, before the game starts, I fill an infinite sequence of bags, making for each the exact same decision as given in \buildRef{ex:averagingDrift:localAveraging}. In each iteration you get the next bag from this sequence. Since the outcome of one bag is independent of other bags, we cannot learn anything about future bags from the history. Thus, using the principle of deferred decision, we compute as if the bag was only packed in the current iteration, after all previous (random) decisions have been made. Thus, the analysis proceeds exactly as in \buildRef{ex:averagingDrift:localAveraging} and you arrive at a value of $10$ drawn red balls after an expected number of $2000/11$ iterations.
\end{beispielEN}

\subsection{The Conditioning of Drift}

	Let us take a closer look at the drift condition. For this section, we define four variants with different conditioning of the drift as follows.

\begin{definitionEN}{Variants on the Conditioning of Drift}{}
	Let $(X_t)_{t \in \natnum}$ be a sequence of random variables over $\realnum$ and let $(F_t)_{t \in \natnum}$ be a filtration such that $(X_t)_{t \in \natnum}$ is adapted to $(F_t)_{t \in \natnum}$, let $f$ be a measurable function and $t \in \natnum$.
	\begin{description}
		\item[(D-filtration)] $\Ew{X_t - X_{t+1} \mid F_t} \geq f(X_t)$ with probability $1$.
		\item[(D-history)] $\Ew{X_t - X_{t+1} \mid X_0,\ldots, X_t} \geq f(X_t)$ with probability $1$.
		\item[(D-events)] For all $s_0,\ldots, s_t$ with $\Pr{X_0=s_0, \ldots, X_t = s_t} \gt 0$,

 $\Ew{X_t - X_{t+1} \mid X_0=s_0,\ldots, X_t=s_t} \geq f(s_t)$.
		\item[(D-Markov)] For all $s_t$ with $\Pr{X_t = s_t} \gt 0$, $\Ew{X_t - X_{t+1} \mid X_t = s_t} \geq f(s_t)$.
	\end{description}
\end{definitionEN}

	Some researchers prefer to state all drift theorems in terms of filtrations, see, for example, \cite{DBLP:journals/tcs/Witt23}; some prefer conditioning on the history \cite{DBLP:journals/tcs/KotzingK19}.
	
	In Section~2.1.2 of \cite{DBLP:series/ncs/Lengler20}, Lengler discusses the differences between \referDefined{(D-filtration)} and \referDefined{(D-Markov)}, phrasing drift theorems in terms of \referDefined{(D-Markov)}. Many applications of drift theorems involve states of algorithms which typically behave as Markov chains. This ubiquity of Markov chains sometimes leads to drift theorems being stated for Markov chains only. However, the states of algorithms need to be mapped to real numbers in order to apply drift theorems (see \buildRef{sec:potentialFunctions} for a discussion on potential functions). If this mapping is not 1-to-1, then, in general, the resulting mapped process is not Markovian any more, so drift theorems applicable for Markov chains on $\realnum$ are no longer applicable. As we will see, \referDefined{(D-history)} is a necessary condition for \referDefined{(D-Markov)} (see \buildRef{thm:averagingDrift:bigImplication}), which is why in this work all drift theorems are stated analogously to \referDefined{(D-history)}.

	In the following we want to discuss how the four given conditions differ and in what sense they are equivalent. First we recall that conditioning on the history $X_0, \ldots, X_t$ is defined as conditioning on the $\sigma$-algebra $\sigma(X_0,\ldots, X_t)$, the canonical filtration. In this sense, the condition \referDefined{(D-history)} implies that there is a filtration such that \referDefined{(D-filtration)} holds.

\begin{propositionEN}[label=thm:averagingDrift:canonicalFiltration]{Canonical Filtration as Filtration}{}
	
		Let $(X_t)_{t \in \natnum}$ be a sequence of random variables over $\realnum$, $f$ a measurable function and $t \in \natnum$. Suppose \referDefined{(D-history)}. Then there is a filtration  $(F_t)_{t \in \natnum}$ such that $(X_t)_{t \in \natnum}$ is adapted to $(F_t)_{t \in \natnum}$ and such that \referDefined{(D-filtration)} holds.

\end{propositionEN}
\begin{proof}
		Using $F_t = \sigma(X_0,\ldots,X_t)$, the canonical filtration, \referDefined{(D-history)} and \referDefined{(D-filtration)} are identical.
	\end{proof}

	The question now arises whether anything can be gained from using other filtrations than the canonical filtration. Sometimes it can be easier to assume a different filtration, which gives more information for the analysis to work with; more outcomes of random variables can be fixed, allowing the analysis to proceed with these concrete outcomes (see, for example, the proof of \buildRef{thm:driftWithoutDrift:variance-two-barrier-hitting-time}). But how should the drift theorem be stated? Using the following theorem, we see that if the drift theorem is only stated conditional on the history, any other filtration (where the process is adapted to) can also be used, since \referDefined{(D-filtration)} implies \referDefined{(D-history)}.

\begin{theoremEN}[label=thm:averagingDrift:bigImplication]{Conditioning on Filtration vs.~History vs.~Events}{}
	
		Let $(X_t)_{t \in \natnum}$ be a sequence of random variables over $\realnum$ and $f$ a measurable function. Suppose further $(F_t)_{t \in \natnum}$ is a filtration such that $(X_t)_{t \in \natnum}$ is adapted to $(F_t)_{t \in \natnum}$. We then have, for all $t \in \natnum$, the following implications:

		\referDefined{(D-filtration)} $\Rightarrow$ \referDefined{(D-history)} $\Rightarrow$ \referDefined{(D-events)} $\Rightarrow$ \referDefined{(D-Markov)}

		where the last implication holds for discrete $(X_t)_{t \in \natnum}$.

\end{theoremEN}
\begin{proof}
		Suppose first, for ``\referDefined{(D-filtration)} $\Rightarrow$ \referDefined{(D-history)}'', $\Ew{X_{t+1} \mid F_t} \geq f(X_t)$ with probability~$1$.
	
		Since $(X_t)_{t \in \natnum}$ is adapted to $(F_t)_{t \in \natnum}$, we have that, for all $t$,
		$$
			\sigma(X_0,\ldots, X_t) \subseteq F_t.
		$$
		Using \buildRef{lem:notation:towerRule} in the second equality, we have, with probability $1$,
		\begin{align*}
			\Ew{X_t - X_{t+1} \mid X_0,\ldots, X_t} & = \Ew{X_t - X_{t+1} \mid \sigma(X_0,\ldots, X_t)}\\
			& = \Ew{\Ew{X_t - X_{t+1} \mid F_t} \mid \sigma(X_0,\ldots, X_t)}\\
			& \geq \Ew{f(X_t) \mid \sigma(X_0,\ldots, X_t)}\\
			& = f(X_t).
		\end{align*}
	
		Suppose now, for ``\referDefined{(D-history)} $\Rightarrow$ \referDefined{(D-events)}'', $\Ew{X_t - X_{t+1} \mid X_0,\ldots, X_t} \geq f(X_t)$ with probability $1$. Let $s_0,\ldots, s_t$ and let $A$ be the event such that $X_0=s_0, \ldots, X_t = s_t$. Suppose $\Pr{A} \gt 0$. Since $A \in \sigma(X_0,\ldots, X_t)$, we get (using \buildRef{lem:notation:conditionalAndIndicator} in the first step and \buildRef{def:notation:filtrationConditional} in the second),
		\begin{align*}
			\Ew{X_t - X_{t+1} \mid A}\Pr{A} &= \Ew{(X_t-X_{t+1}) \ind{A}}\\
			&= \Ew{\Ew{X_t-X_{t+1} \mid X_0,\ldots, X_t} \ind{A}}\\
			&\geq \Ew{f(X_t)\ind{A}}\\
			&= f(s_t)\Pr{A}.
		\end{align*}
		
		Regarding \referDefined{(D-events)} $\Rightarrow$ \referDefined{(D-Markov)}, we note that, for discrete $(X_t)_{t \in \natnum}$, and any $s_t$ with $\Pr{X_t = s_t} \gt 0$, we can find $s_0,\ldots,s_{t-1}$ such that $\Pr{X_0=s_0, \ldots, X_t = s_t} \gt 0$. For such $s_0,\ldots,s_t$ we then have
		$$\Ew{X_t - X_{t+1} \mid X_t=s_t} = \Ew{X_t - X_{t+1} \mid X_0=s_0,\ldots, X_t=s_t}.$$
		This gives the desired implication.
	\end{proof}

	Note that \referDefined{(D-events)} and \referDefined{(D-Markov)} implicitly consider the process to be discrete. As we see in the example given next, a drift theorem based on \referDefined{(D-events)} without the the requirement of a discrete process would in this generality be wrong.

\begin{beispielEN}[label=ex:averagingDrift:failedEventsDrift]{Conditioning on Events of Continuous Processes}{}
	Let $X_0$ be a uniformly random number from $[1,2]$ and let, for all $t \in \natnum$, $X_{t+1} = X_t$. Then we have
	\begin{description}
		\item[(D-events)] for all $s_0,\ldots, s_t$ with $\Pr{X_0=s_0, \ldots, X_t = s_t} \gt 0$, $E[X_t - X_{t+1} \mid X_0=s_0,\ldots, X_t=s_t] \geq 1$.
	\end{description}
	This follows since, for all $t \in \natnum$ and all $s_0,\ldots, s_t$, we have $\Pr{X_0=s_0, \ldots, X_t = s_t} = 0$ (we have a truly continuous random variable). Thus, \referDefined{(D-events)} is vacuously true. Furthermore, for all $t \in \natnum$, $X_t \geq 1$, so there is no $t$ such that $X_t \leq 0$.
\end{beispielEN}

	The example shows that the problem arises when considering continuous random variables. The next proposition shows that, for discrete processes, \referDefined{(D-events)} implies \referDefined{(D-history)}, making these two conditions equivalent. The proof is due to Marcus Pappik (private communication).

\begin{propositionEN}[label=thm:averagingDrift:discreteImpliesAll]{Equivalence of Conditionals for Discrete Search Spaces}{}
	
		Let $(X_t)_{t \in \natnum}$ be a \emph{discrete} process over $\realnum$, $f$ a measurable function and $t \in \natnum$. 
		Then \referDefined{(D-events)} implies \referDefined{(D-history)}.

\end{propositionEN}
\begin{proof}
		Fix $t \geq 0$, and let $R := \mathrm{range}(X_0) \times \ldots \times \mathrm{range}(X_t)$.
		For a tuple $s = (s_i)_{0 \le i \le t} \in R$, let $A(s) := \{\forall 0 \leq i \leq t: X_i = s_i\}$.
		Let $S := \set{s \in R}{\Pr{A(s)} \gt 0}$. 
		For all $s = (s_i)_{0 \le i \le t} \in S$ and $\omega \in A(s)$, \buildRef{lem:notation:conditioningHistoryEvents} and \referDefined{(D-event)} yield
		$$
			\Ew{X_t - X_{t+1} \mid X_0, \ldots, X_t}(\omega) = \Ew{X_t - X_{t+1} \mid A(s)} \geq f(s_t) = f(X_t(\omega)).
		$$
		Since further every $X_i$ has countable range, it holds in particular that $S$ and $R$ are countable. 
		Hence, we have that $\bigcup_{s \in S} A(s)$ is measurable and 
		$$
			\Pr{\bigcup_{s \in S} A(s)} \geq 1 - \sum_{s \in R \setminus S} \Pr{A(s)} = 1 .
		$$
		Therefore, we have that $\Ew{X_t - X_{t+1} \mid X_0, \ldots, X_t} \geq f(X_t)$ holds with probability $1$.
	\end{proof}

	Just as \referDefined{(D-events)} implicitly assumes a discrete space, \referDefined{(D-Markov)} assumes the process to be Markovian. The following theorem shows that, for discrete Markov chains, \referDefined{(D-Markov)} implies \referDefined{(D-events)}, making also these two conditions equivalent in this case.

\begin{propositionEN}[label=thm:averagingDrift:equivalenceMarkov]{Equivalence for Markov Chains}{}
	
		Let $(X_t)_{t \in \natnum}$ be a \emph{discrete} Markov chain over $\realnum$, $f$ a measurable function and $t \in \natnum$. 
		Then \referDefined{(D-Markov)} implies \referDefined{(D-events)}.

\end{propositionEN}
\begin{proof}
		This follows directly from the Markov property that, for all $s_0,\ldots,s_t$ with $\Pr{X_0=s_0,\ldots, X_t=s_t} \gt 0$, $\Ew{X_{t+1} \mid X_t=s_t} = \Ew{X_{t+1} \mid X_0=s_0,\ldots, X_t=s_t}$.
	\end{proof}

	We have seen that for discrete spaces and for Markov chains, one can give specialized formulations of drift theorems. However, drift theorems conditioning on the history are strictly more general. Furthermore, conditioning on the history is easy to state and understand for users of the theorem.

	Conditioning on a filtration results in drift theorems equally general as those conditioning on the history, since the history is one possible filtration, and in fact the least restrictive.

\subsection{Reaching the Target}

	Let us consider again the drift condition
	\begin{description}
		\item[(D)] there is a $\delta > 0$ such that, for all $t \lt T$, it holds that $\Ew{X_t - X_{t+1} \mid X_0,\ldots,X_t} \geq \delta$.
	\end{description}
	We want to take a closer look at the requirement ``for all $t \lt T$''. Since $T$ is a random variable, this does not properly define a range for $t$. This is desirable to allow for processes which naturally do not go down any more after reaching the target. One clean way to write it would be 
	\begin{description}
		\item[(D)] there is a $\delta > 0$ such that, for all $t \in \natnum$, $\Ew{X_t - X_{t+1} \mid X_0,\ldots,X_t} \ind{t \lt T} \geq \delta \ind{t \lt T}$.
	\end{description}
	This notation resorts to indicator random variables and, while quantifying $t$ over all of $\natnum$, effectively requires the inequality to hold only in case of $t \lt T$. This inspires the following convention.

\begin{conventionEN}{Drift While not at Target}{}
	We state inequalities that only need to hold for points in time when a random process did not reach its target yet.
	Formally, let~$T$ be a random variable over $\natnum \cup \{\infty\}$, let~$X$ and~$Y$ be random variables over~$\realnum$.
	
	Further, let~$\sim$ denote a relation symbol, such as $=$, $\leq$, or $\geq$.
	We define the phrase ``for all $t \lt T$, it holds that $X \sim Y$'' to be equivalent to ``for all $t \in \natnum$, it holds that $X \cdot \ind{t \lt T} \sim Y \cdot \ind{t \lt T}$''.
\end{conventionEN}

	Note that, alternatively, we can \emph{condition} on $t \leq T$, the way chosen in \buildRef{thm:advancedDriftTheorems:additiveTimeConditionedUpper} and \buildRef{thm:advancedDriftTheorems:additiveTimeConditionedLower}. This makes the dependence on $t \lt T$ very explicit (this was the reason for stating it in this way in the two named theorems). However, now one has to quantify over ``all $t \in \natnum$ such that $\Pr{t \lt T}\gt 0$'', which is again somewhat cumbersome.

\clearpage

\section{Notation}
\label{sec:notation}

	In this section we collect some algorithms, notation and lemmas used in this work.

\subsection{Algorithms}
\label{sec:algorithms}

	The most simple search heuristic is \emph{Random Local Search} (RLS). It starts with a random bit string and iteratively tries to improve it by changing the currently best search point in exactly one position (we use this as the ``flipOne'' function below). The pseudo code for RLS maximizing a given function $f\colon \BitStrings \rightarrow \realnum$ is given as follows.

\begin{center}
\includegraphics[width=80mm]{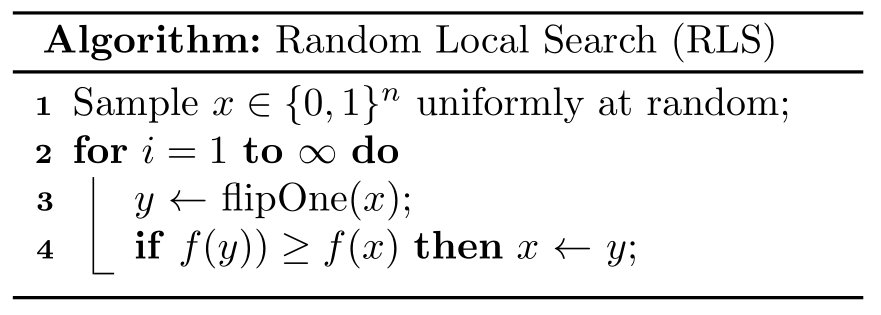}

\end{center}

	The algorithm is set up to maximize the given function $f$; by turning the inequality around, we get the analogous algorithm for minimization.

	RLS constitutes a simple and straightforward hill climber. A slightly more advanced algorithm allows for larger \emph{jumps}, that is, it also considers changing the currently best search point in more than one position. The most common way to achieve this is by flipping not exactly one bit, but each bit independently with some predefined probability $p$. This independently random flipping of bits is called \emph{mutation}, and the resulting algorithm is called the \OneOneEA; its pseudo code is given as follows.

\begin{center}
\includegraphics[width=80mm]{OneOneEA.png}

\end{center}

	Note that the standard bit flip probability is $p=1/n$, implying that, on average, exactly one bit changes.

	We consider in this document two concrete test functions and one function class as follows.

\begin{definitionEN}{Test Functions}{}
	\begin{itemize}
			\item $\OneMax$ is a function $\BitStrings \rightarrow \realnum$ mapping any bit string to the number of $1$s in the bit string.
			\item $\LeadingOnes$ is a function $\BitStrings \rightarrow \realnum$ mapping any bit string to the number of $1$s \emph{before the first $0$} (if any) in the bit string (the number of leading $1$s).
			\item A \emph{linear function} is any function $f\colon \BitStrings \rightarrow \realnum$ such that there exists $w_1,\ldots,w_n \in \realnum$ such that, for all $x \in \BitStrings$, $f(x) = \sum_{i=1}^n w_i \; x_i$.	
	\end{itemize}
\end{definitionEN}

\subsection{Notation}

	Next we give some non-standard notation.

\begin{definitionEN}{Discrete Intervals}{}
	For any $n,m \in \natnum$ with $n \leq m$, we use $[n..m]$ to denote the set $\set{i \in \natnum}{n \leq i \leq m}$. Furthermore, for any $n \in \natnum_+$, we will write $[n]$ for $[1..n]$.
\end{definitionEN}

\begin{definitionEN}{Function Self-Composition}{}
	For any function $f\colon X \rightarrow X$ and $i \geq 0$, we let $f^i$ denote the $i$-times self-composition of $f$ (with $f^0$ being the identity on $X$).
\end{definitionEN}

	We use the following notation regarding probabilities.

\begin{definitionEN}{Indicator Function}{}
	For any event~$A$, let $\ind{A}$ denote the indicator function for the event~$A$.
\end{definitionEN}

\begin{definitionEN}{Conditioning the Expectation}{}
	For any discrete random variable $X$ and any event $A$ such that $\Pr{A} \gt 0$, we have
	$$\Ew{X \mid A} = \sum_{x} x \frac{\Pr{\{X=x\} \cap A}}{\Pr{A}}.$$
\end{definitionEN}

	Note that we can only condition on the event $X=s$ if $\Pr{X = s} \gt 0$.

\begin{definitionEN}[label=def:notation:integrable]{Integrability}{}
	A random variable~$X$ is \emph{integrable} if and only if $\Ew{|X|} \lt \infty$. In general, a random process $(X_t)_{t \in \natnum}$ is \emph{integrable} if and only if, for all $t \in \natnum$, it holds that $X_t$ is inte\emph{}grable.
\end{definitionEN}

\begin{definitionEN}[label=def:notation:discrete]{Discrete Random Variable, Discrete Random Process}{}
	A random variable $X$ is \emph{discrete} if and only if it has a countable range. We call a random process $(X_t)_{t \in \natnum}$ disc\emph{}rete if each $X_t$ is disc\emph{}rete.
\end{definitionEN}

\begin{definitionEN}{Monotone Process}{}
	A random process~$(X_t)_{t \in \natnum}$ is \emph{monotone} (sometimes called \emph{monotone non-decreasing}) if and only if, for all $t$, $X_t \leq X_{t+1}$ (point-wise, i.e., for all atomic events $\omega \in \Omega$ we have $X_t(\omega) \leq X_{t+1}(\omega)$).
\end{definitionEN}

\begin{definitionEN}[label=def:notation:stochasticDominance]{Stochastic Dominance}{}
	Given two random variables $X,Y$ over $\realnum$, we say that \emph{$Y$ stochastically dominates $X$} if, for all $x \in \realnum$, $\Pr{X \leq x} \geq \Pr{Y \leq x}$; we write $X \preceq Y$.
\end{definitionEN}

\begin{definitionEN}[label=def:notation:markovChain]{Markov Chain}{}
	A discrete random process $(X_t)_{t \in \natnum}$ is a \emph{Markov chain} if and only if the outcome of each next step only depends on the current state and time point.
	Formally, for all $t \in \natnum$ as well as all $s \in \realnum$ and all $v \in \realnum^t$, it holds that $\Pr{X_{t+1} = s \mid X_t = v_t} = \Pr{X_{t+1} = s \mid \forall t' \in [0..t]\colon X_{t'} = v_{t'}}$.
	A Mar\emph{}kov chain $(X_t)_{t \in \natnum}$ is \emph{time-homogeneous} if and only if the outcome of each next step only depends on the current state but not the current time point.
	Formally, for all $t,k \in \natnum$ as well as all $s, u \in \realnum$, it holds that $\Pr{X_{t+1} = s \mid X_t = u} = \Pr{X_{t+k+1} = s \mid X_{t + k} = u}$.
\end{definitionEN}

\subsection{Lemmas}

	We will make use of the following lemmas. The first two pertain to bounding certain sums.

\begin{lemmaEN}[label=lem:notation:boundHarmonicSum]{Upper Bound on the Harmonic Sum}{}
	
		For all $n \in \natnum_{+}$ we have
		$$\sum_{i=1}^n \frac{1}{i} \leq \ln(n)+1.$$

\end{lemmaEN}
\begin{proof}
		See (1.4.12) in \cite{DBLP:series/ncs/Doerr20}.
	\end{proof}

\begin{lemmaEN}[label=lem:notation:sumByIntegral]{Upper Bounding Sum by Integral}{}
	
		Let $a,b \in \realnum$ with $a \leq b$ and let $f\colon [a,b] \rightarrow \realnum$ be integra\emph{}ble. If $f$ is monotone increasing, then
		$$\sum_{i=a}^{b-1} f(i) \leq \int_a^b f(x) \; \mathrm{d}x.$$
		If $f$ is monotone decreasing, then
		$$\sum_{i=a+1}^{b} f(i) \leq \int_a^b f(x) \; \mathrm{d}x.$$
	
\end{lemmaEN}

	The remaining lemmas pertain to random variables.

\begin{lemmaEN}[label=lem:notation:conditionalAndIndicator]{Conditional and Indicator}{}
	
		For all events $A$ we have
		$$\Ew{Y \mid A} \cdot \Pr{A} = \Ew{Y \ind{A}}.$$
	
\end{lemmaEN}

\begin{definitionEN}[label=def:notation:filtrationConditional]{Filtration Conditional}{}
	Let a sub-$\sigma$-algebra $F$ of the underlying probability space be given and let $X$ be a random variable. Then $\Ew{Y \mid F}$ refers to any random variable such that, for all $A \in F$ with probability $1$,
		$$\Ew{Y \cdot \ind{A}} = \Ew{\Ew{Y \mid F} \cdot \ind{A} }.$$
\end{definitionEN}

\begin{lemmaEN}[label=lem:notation:towerRule]{Tower Property for Sub-$\sigma$-Algebra}{}
	 Let two sub-$\sigma$-algebras $F_1 \subseteq F_2$ of the underlying probability space be given and let $X$ be a random variable. Then, with probability $1$,
		$$\Ew{\Ew{X \mid F_2} \mid F_1} = \Ew{X \mid F_1}.$$
	
\end{lemmaEN}

	Very much related to the preceding lemma is the law of total  expectation (also known under other names and with different formulations).

\begin{lemmaEN}[label=lem:notation:TotalExpectation]{Law of Total Expectation}{}
	 Let $X$ be a random variable and $A_i, \ldots, A_n$ disjoint measurable events with positive probability that partition the probability space. Then we have
		$$\Ew{X} = \sum_{i=1}^n \Ew{X \mid A_i}\Pr{A_i}.$$
	We will also have cause to use a different way to state this. Let $X,Y$ be two random variables. Then
		$$\Ew{X} = \Ew{\Ew{X \mid Y}}.$$
	
\end{lemmaEN}

	We need the following lemma to make a specific proof in this document rigorous. The proof is due to Marcus Pappik (private communication).

\begin{lemmaEN}[label=lem:notation:conditioningHistoryEvents]{Conditioning on History vs.~Events}{}

	    Let $Y, X_0, \ldots, X_t$ be random variables and let $Z$ be any version of the conditional expectation $\Ew{Y \mid X_0, \ldots, X_t}$. 
		For all $s_0, \ldots, s_t \in \realnum$ such that $\Pr{X_0 = s_0, \ldots, X_t = s_t} > 0$ and all $\omega \in \{X_0 = s_0, \ldots, X_t = s_t\}$ it holds that
		$$
			Z(\omega) = \Ew{Y \mid X_0 = s_0, \ldots, X_t = s_t}.
		$$

\end{lemmaEN}
\begin{proof}
	We start by proving the following claim.

	\textbf{Claim.} Consider two measurable spaces $(\Omega, \mathcal{F})$ and $(S, \mathcal{S})$ such that, for all $s \in S$, $\{s\} \in \mathcal{S}$.
    Let $X, Z\colon \Omega \to S$ be measurable. If $Z$ is $\sigma(X)$-measurable then $Z$ is constant on $\{X = s\}$ for every $s \in S$.

	\textbf{Proof of claim.} For the sake of contradiction, assume the statement is false.
    That is, suppose $Z$ is $\sigma(X)$-measurable and there is some $s_0 \in S$ such that $Z$ is not constant on $A_0 := \{X=s_0\}$.
    Let $s_1, s_2 \in S$ be distinct values of $Z$ on $A_0$.
    Since $Z$ is $\sigma(X)$-measurable, it holds that $A_1 := \{Z = s_1\} \in \sigma(X)$.
    Consequently, we have that $A_0 \cap A_1 \in \sigma(X)$ and, by construction, $\emptyset \subset A_0 \cap A_1 \subset A_0$.
    We now show that these two properties of $A_0 \cap A_1$ pose a contradiction.
    To this end, suppose there exists any set $A \subseteq \Omega$ with $\emptyset \subset A \subset A_0$ and $A \in \sigma(X)$.
    Note that
    $$
        \sigma(X) = \{X^{-1}(B) \mid B \in \mathcal{S}\} .
    $$
    Hence, if $A \in \sigma(X)$ then $B := \{X(\omega) \mid \omega \in A\}$ must be in $\mathcal{S}$.
    But since $\emptyset \subset A \subset A_0$ it holds that
    $$
        \emptyset \subset B \subset \{X(a) \mid a \in A_0\} = \{s_0\} . 
    $$
    However, since both inclusions are strict, such a set $B$ cannot exist.
	\textbf{End of proof of claim.}

	To prove the lemma, suppose now that $Z$ is any version of $\Ew{Y \mid X_0, \ldots, X_t}$ and let $A = \{X_0 = s_0, \ldots, X_t = s_t\}$ for any $s_0, \ldots, s_t \in \realnum$.
    Since $Z$ is $\sigma(X_0, \ldots, X_t)$-measurable, the claim above yields that $Z$ is constant on $A$.
    Let $z$ be the value of $Z$ on $A$ and note that, for all $\omega \in A$, we have
    $$
        Z(\omega) \cdot \Pr{A} = \Ew{z \cdot \ind{A}} = \Ew{Z \cdot \ind{A}} = \Ew{Y \cdot \ind{A}} ,
    $$
    where the first equality is due to $Z(\omega) = z$ and $\Pr{A} = \Ew{\ind{A}}$, the second equality follows from the fact that $z \cdot \ind{A} = Z \cdot \ind{A}$ point-wise, and the last equality is due to the fact that $Z$ is a version of $\Ew{Y \mid X_0, \ldots, X_t}$.
    Further, by \buildRef{lem:notation:conditionalAndIndicator}, we have 
    $$
        \Ew{Y \cdot \ind{A}} = \Ew{Y \mid A} \cdot \Pr{A} .
    $$
    Combining both steps yields
    $$
        Z(\omega) \cdot \Pr{A} = \Ew{Y \mid A} \cdot \Pr{A}.
    $$
    Provided $\Pr{A} > 0$, this implies $Z(\omega) = \Ew{Y \mid A}$ for all $\omega \in A$, which proves the lemma.
	\end{proof}

\clearpage

\section{Acknowledgments}
\label{sec:acknowledgments}

	While many people were important for this habilitation document, three people come to mind in particular. First of all Tobias Friedrich, who had my back for many years, providing an outstanding work environment for me to thrive. Most inspirational for my development as a researcher was Benjamin Doerr, who acted as a mentor while pretending to be merely a collaborator; I learned a lot from Benjamin and enjoyed our many projects together. Furthermore, Martin Krejca was my PhD student and we had many fruitful collaborations; now he is a good friend, an impeccable discussion partner and generally a person I can rely on in many situations.

	The community of researchers working on the Theory of Randomized Search Heuristics gave me a home as a scientist. The reviewing is thorough, the discussions interesting and the atmosphere very friendly. I want to point out particularly Carsten Witt, Carola Doerr, Johannes Lengler and Dirk Sudholt (all are also collaborators), but many other people are also significant for the community.

	Finally, many researchers from the research group I am in contributed to this thesis either with discussions or by reading and commenting on drafts. In no particular order, ``thank you'' to Anton Krohmer, Marcus Pappik, Nicolas Klodt, Aishwarya Radhakrishnan, Xiaoyue Sherry Li, Nadym Mallek, Stefan Neubert, Janosch Ruff, Ziena Zeif and Sam Baguley.

\bibliographystyle{alpha}
\bibliography{knownKeys.bib,knownKeysStatic.bib}

\newcommand{\etalchar}[1]{$^{#1}$}
\begin{thebibliography}{FKKS15b}

\bibitem[BDN10]{DBLP:conf/ppsn/BottcherDN10}
S{\"{u}}ntje B{\"{o}}ttcher, Benjamin Doerr, and Frank Neumann.
\newblock Optimal fixed and adaptive mutation rates for the leadingones
  problem.
\newblock In Robert Schaefer, Carlos Cotta, Joanna Kolodziej, and G{\"{u}}nter
  Rudolph, editors, {\em Parallel Problem Solving from Nature - {PPSN} XI, 11th
  International Conference, Krak{\'{o}}w, Poland, September 11-15, 2010,
  Proceedings, Part {I}}, volume 6238 of {\em Lecture Notes in Computer
  Science}, pages 1--10. Springer, 2010.

\bibitem[BLM{\etalchar{+}}20]{DBLP:conf/esa/BertschingerLMM20}
Daniel Bertschinger, Johannes Lengler, Anders Martinsson, Robert Meier,
  Angelika Steger, Milos Trujic, and Emo Welzl.
\newblock An optimal decentralized ({\(\Delta\)} + 1)-coloring algorithm.
\newblock In Fabrizio Grandoni, Grzegorz Herman, and Peter Sanders, editors,
  {\em 28th Annual European Symposium on Algorithms, {ESA} 2020, September 7-9,
  2020, Pisa, Italy (Virtual Conference)}, volume 173 of {\em LIPIcs}, pages
  17:1--17:12. Schloss Dagstuhl - Leibniz-Zentrum f{\"{u}}r Informatik, 2020.

\bibitem[CDEL18]{DBLP:journals/tec/CorusDEL18}
Dogan Corus, Duc{-}Cuong Dang, Anton~V. Eremeev, and Per~Kristian Lehre.
\newblock Level-based analysis of genetic algorithms and other search
  processes.
\newblock {\em {IEEE} Trans. Evol. Comput.}, 22(5):707--719, 2018.

\bibitem[DDK15]{DBLP:conf/gecco/DoerrDK15}
Benjamin Doerr, Carola Doerr, and Timo K{\"{o}}tzing.
\newblock Solving problems with unknown solution length at (almost) no extra
  cost.
\newblock In Sara Silva and Anna~Isabel Esparcia{-}Alc{\'{a}}zar, editors, {\em
  Proceedings of the Genetic and Evolutionary Computation Conference, {GECCO}
  2015, Madrid, Spain, July 11-15, 2015}, pages 831--838. {ACM}, 2015.

\bibitem[DDK16]{DBLP:conf/gecco/DoerrDK16}
Benjamin Doerr, Carola Doerr, and Timo K{\"{o}}tzing.
\newblock The right mutation strength for multi-valued decision variables.
\newblock In Tobias Friedrich, Frank Neumann, and Andrew~M. Sutton, editors,
  {\em Proceedings of the 2016 on Genetic and Evolutionary Computation
  Conference, Denver, CO, USA, July 20 - 24, 2016}, pages 1115--1122. {ACM},
  2016.

\bibitem[DDK18]{DBLP:journals/algorithmica/DoerrDK18}
Benjamin Doerr, Carola Doerr, and Timo K{\"{o}}tzing.
\newblock Static and self-adjusting mutation strengths for multi-valued
  decision variables.
\newblock {\em Algorithmica}, 80(5):1732--1768, 2018.

\bibitem[DDY20]{DBLP:journals/tcs/DoerrDY20}
Benjamin Doerr, Carola Doerr, and Jing Yang.
\newblock Optimal parameter choices via precise black-box analysis.
\newblock {\em Theor. Comput. Sci.}, 801:1--34, 2020.

\bibitem[DFF{\etalchar{+}}19]{DBLP:journals/algorithmica/DoerrFFFKS19}
Benjamin Doerr, Philipp Fischbeck, Clemens Frahnow, Tobias Friedrich, Timo
  K{\"{o}}tzing, and Martin Schirneck.
\newblock Island models meet rumor spreading.
\newblock {\em Algorithmica}, 81(2):886--915, 2019.

\bibitem[DG10]{DBLP:conf/ppsn/DoerrG10a}
Benjamin Doerr and Leslie~Ann Goldberg.
\newblock Drift analysis with tail bounds.
\newblock In Robert Schaefer, Carlos Cotta, Joanna Kolodziej, and G{\"{u}}nter
  Rudolph, editors, {\em Parallel Problem Solving from Nature - {PPSN} XI, 11th
  International Conference, Krak{\'{o}}w, Poland, September 11-15, 2010,
  Proceedings, Part {I}}, volume 6238 of {\em Lecture Notes in Computer
  Science}, pages 174--183. Springer, 2010.

\bibitem[DJW02]{DBLP:journals/tcs/DrosteJW02}
Stefan Droste, Thomas Jansen, and Ingo Wegener.
\newblock On the analysis of the {(1+1)} evolutionary algorithm.
\newblock {\em Theor. Comput. Sci.}, 276(1-2):51--81, 2002.

\bibitem[DJW10]{DBLP:conf/gecco/DoerrJW10}
Benjamin Doerr, Daniel Johannsen, and Carola Winzen.
\newblock Multiplicative drift analysis.
\newblock In Martin Pelikan and J{\"{u}}rgen Branke, editors, {\em Genetic and
  Evolutionary Computation Conference, {GECCO} 2010, Proceedings, Portland,
  Oregon, USA, July 7-11, 2010}, pages 1449--1456. {ACM}, 2010.

\bibitem[DJW12]{DBLP:journals/algorithmica/DoerrJW12}
Benjamin Doerr, Daniel Johannsen, and Carola Winzen.
\newblock Multiplicative drift analysis.
\newblock {\em Algorithmica}, 64(4):673--697, 2012.

\bibitem[DJWZ13]{DBLP:conf/gecco/DoerrJWZ13}
Benjamin Doerr, Thomas Jansen, Carsten Witt, and Christine Zarges.
\newblock A method to derive fixed budget results from expected optimisation
  times.
\newblock In Christian Blum and Enrique Alba, editors, {\em Genetic and
  Evolutionary Computation Conference, {GECCO} '13, Amsterdam, The Netherlands,
  July 6-10, 2013}, pages 1581--1588. {ACM}, 2013.

\bibitem[DK14]{DBLP:conf/analco/DoerrK14}
Benjamin Doerr and Marvin K{\"{u}}nnemann.
\newblock Tight analysis of randomized rumor spreading in complete graphs.
\newblock In Michael Drmota and Mark~Daniel Ward, editors, {\em 2014
  Proceedings of the Eleventh Workshop on Analytic Algorithmics and
  Combinatorics, {ANALCO} 2014, Portland, Oregon, USA, January 6, 2014}, pages
  82--91. {SIAM}, 2014.

\bibitem[DK21a]{DBLP:conf/gecco/DoerrK21}
Benjamin Doerr and Timo K{\"{o}}tzing.
\newblock Lower bounds from fitness levels made easy.
\newblock In Francisco Chicano and Krzysztof Krawiec, editors, {\em {GECCO}
  '21: Genetic and Evolutionary Computation Conference, Lille, France, July
  10-14, 2021}, pages 1142--1150. {ACM}, 2021.

\bibitem[DK21b]{DBLP:journals/algorithmica/DoerrK21}
Benjamin Doerr and Timo K{\"{o}}tzing.
\newblock Multiplicative up-drift.
\newblock {\em Algorithmica}, 83(10):3017--3058, 2021.

\bibitem[DKLL17]{DBLP:conf/gecco/DoerrKLL17}
Benjamin Doerr, Timo K{\"{o}}tzing, J.~A.~Gregor Lagodzinski, and Johannes
  Lengler.
\newblock Bounding bloat in genetic programming.
\newblock In Peter A.~N. Bosman, editor, {\em Proceedings of the Genetic and
  Evolutionary Computation Conference, {GECCO} 2017, Berlin, Germany, July
  15-19, 2017}, pages 921--928. {ACM}, 2017.

\bibitem[DKLL20]{DBLP:journals/tcs/DoerrKLL20}
Benjamin Doerr, Timo K{\"{o}}tzing, J.~A.~Gregor Lagodzinski, and Johannes
  Lengler.
\newblock The impact of lexicographic parsimony pressure for {ORDER/MAJORITY}
  on the run time.
\newblock {\em Theor. Comput. Sci.}, 816:144--168, 2020.

\bibitem[DL16]{DBLP:journals/algorithmica/DangL16}
Duc{-}Cuong Dang and Per~Kristian Lehre.
\newblock Runtime analysis of non-elitist populations: From classical
  optimisation to partial information.
\newblock {\em Algorithmica}, 75(3):428--461, 2016.

\bibitem[Doe20]{DBLP:series/ncs/Doerr20}
Benjamin Doerr.
\newblock Probabilistic tools for the analysis of randomized optimization
  heuristics.
\newblock In Benjamin Doerr and Frank Neumann, editors, {\em Theory of
  Evolutionary Computation - Recent Developments in Discrete Optimization},
  Natural Computing Series, pages 1--87. Springer, 2020.

\bibitem[Dro04]{DBLP:conf/gecco/Droste04}
Stefan Droste.
\newblock Analysis of the {(1+1)} {EA} for a noisy onemax.
\newblock In Kalyanmoy Deb, Riccardo Poli, Wolfgang Banzhaf, Hans{-}Georg
  Beyer, Edmund~K. Burke, Paul~J. Darwen, Dipankar Dasgupta, Dario Floreano,
  James~A. Foster, Mark Harman, Owen Holland, Pier~Luca Lanzi, Lee Spector,
  Andrea Tettamanzi, Dirk Thierens, and Andrew~M. Tyrrell, editors, {\em
  Genetic and Evolutionary Computation - {GECCO} 2004, Genetic and Evolutionary
  Computation Conference, Seattle, WA, USA, June 26-30, 2004, Proceedings, Part
  {I}}, volume 3102 of {\em Lecture Notes in Computer Science}, pages
  1088--1099. Springer, 2004.

\bibitem[FK13]{DBLP:conf/foga/FeldmannK13}
Matthias Feldmann and Timo K{\"{o}}tzing.
\newblock Optimizing expected path lengths with ant colony optimization using
  fitness proportional update.
\newblock In Frank Neumann and Kenneth A.~De Jong, editors, {\em Foundations of
  Genetic Algorithms XII, {FOGA} '13, Adelaide, SA, Australia, January 16-20,
  2013}, pages 65--74. {ACM}, 2013.

\bibitem[FKK16]{DBLP:conf/gecco/FriedrichKK16}
Tobias Friedrich, Timo K{\"{o}}tzing, and Martin~S. Krejca.
\newblock Edas cannot be balanced and stable.
\newblock In Tobias Friedrich, Frank Neumann, and Andrew~M. Sutton, editors,
  {\em Proceedings of the 2016 on Genetic and Evolutionary Computation
  Conference, Denver, CO, USA, July 20 - 24, 2016}, pages 1139--1146. {ACM},
  2016.

\bibitem[FKKS15a]{DBLP:conf/isaac/FriedrichKKS15}
Tobias Friedrich, Timo K{\"{o}}tzing, Martin~S. Krejca, and Andrew~M. Sutton.
\newblock The benefit of recombination in noisy evolutionary search.
\newblock In Khaled~M. Elbassioni and Kazuhisa Makino, editors, {\em Algorithms
  and Computation - 26th International Symposium, {ISAAC} 2015, Nagoya, Japan,
  December 9-11, 2015, Proceedings}, volume 9472 of {\em Lecture Notes in
  Computer Science}, pages 140--150. Springer, 2015.

\bibitem[FKKS15b]{DBLP:conf/gecco/FriedrichKKS15}
Tobias Friedrich, Timo K{\"{o}}tzing, Martin~S. Krejca, and Andrew~M. Sutton.
\newblock Robustness of ant colony optimization to noise.
\newblock In Sara Silva and Anna~Isabel Esparcia{-}Alc{\'{a}}zar, editors, {\em
  Proceedings of the Genetic and Evolutionary Computation Conference, {GECCO}
  2015, Madrid, Spain, July 11-15, 2015}, pages 17--24. {ACM}, 2015.

\bibitem[FKKS17]{DBLP:journals/tec/0001KKS17}
Tobias Friedrich, Timo K{\"{o}}tzing, Martin~S. Krejca, and Andrew~M. Sutton.
\newblock The compact genetic algorithm is efficient under extreme gaussian
  noise.
\newblock {\em {IEEE} Trans. Evol. Comput.}, 21(3):477--490, 2017.

\bibitem[FKL{\etalchar{+}}17]{DBLP:conf/foga/0001KLNS17}
Tobias Friedrich, Timo K{\"{o}}tzing, Gregor Lagodzinski, Frank Neumann, and
  Martin Schirneck.
\newblock Analysis of the {(1+1)} {EA} on subclasses of linear functions under
  uniform and linear constraints.
\newblock In Christian Igel, Dirk Sudholt, and Carsten Witt, editors, {\em
  Proceedings of the 14th {ACM/SIGEVO} Conference on Foundations of Genetic
  Algorithms, {FOGA} 2017, Copenhagen, Denmark, January 12-15, 2017}, pages
  45--54. {ACM}, 2017.

\bibitem[FKL{\etalchar{+}}20]{DBLP:journals/tcs/FriedrichKLNS20}
Tobias Friedrich, Timo K{\"{o}}tzing, J.~A.~Gregor Lagodzinski, Frank Neumann,
  and Martin Schirneck.
\newblock Analysis of the (1+1) {EA} on subclasses of linear functions under
  uniform and linear constraints.
\newblock {\em Theor. Comput. Sci.}, 832:3--19, 2020.

\bibitem[FKM17]{DBLP:conf/gecco/0001KM17}
Tobias Friedrich, Timo K{\"{o}}tzing, and Anna Melnichenko.
\newblock Analyzing search heuristics with differential equations.
\newblock In Peter A.~N. Bosman, editor, {\em Genetic and Evolutionary
  Computation Conference, Berlin, Germany, July 15-19, 2017, Companion Material
  Proceedings}, pages 313--314. {ACM}, 2017.

\bibitem[FKN{\etalchar{+}}23]{DBLP:conf/gecco/0001KN0R23}
Tobias Friedrich, Timo K{\"{o}}tzing, Aneta Neumann, Frank Neumann, and
  Aishwarya Radhakrishnan.
\newblock Analysis of {(1+1)} {EA} on leadingones with constraints.
\newblock In Sara Silva and Lu{\'{\i}}s Paquete, editors, {\em Proceedings of
  the Genetic and Evolutionary Computation Conference, {GECCO} 2023, Lisbon,
  Portugal, July 15-19, 2023}, pages 1584--1592. {ACM}, 2023.

\bibitem[FKS16]{DBLP:conf/ppsn/FriedrichKS16}
Tobias Friedrich, Timo K{\"{o}}tzing, and Andrew~M. Sutton.
\newblock On the robustness of evolving populations.
\newblock In Julia Handl, Emma Hart, Peter~R. Lewis, Manuel
  L{\'{o}}pez{-}Ib{\'{a}}{\~{n}}ez, Gabriela Ochoa, and Ben Paechter, editors,
  {\em Parallel Problem Solving from Nature - {PPSN} {XIV} - 14th International
  Conference, Edinburgh, UK, September 17-21, 2016, Proceedings}, volume 9921
  of {\em Lecture Notes in Computer Science}, pages 771--781. Springer, 2016.

\bibitem[GK14]{DBLP:conf/gecco/GiessenK14}
Christian Gie{\ss}en and Timo K{\"{o}}tzing.
\newblock Robustness of populations in stochastic environments.
\newblock In Dirk~V. Arnold, editor, {\em Genetic and Evolutionary Computation
  Conference, {GECCO} '14, Vancouver, BC, Canada, July 12-16, 2014}, pages
  1383--1390. {ACM}, 2014.

\bibitem[GK16]{DBLP:journals/algorithmica/GiessenK16}
Christian Gie{\ss}en and Timo K{\"{o}}tzing.
\newblock Robustness of populations in stochastic environments.
\newblock {\em Algorithmica}, 75(3):462--489, 2016.

\bibitem[GKS99]{DBLP:journals/ec/GarnierKS99}
Josselin Garnier, Leila Kallel, and Marc Schoenauer.
\newblock Rigorous hitting times for binary mutations.
\newblock {\em Evol. Comput.}, 7(2):173--203, 1999.

\bibitem[GLR20]{DBLP:journals/rsa/GoldbergLR20}
Leslie~Ann Goldberg, John Lapinskas, and David Richerby.
\newblock Phase transitions of the moran process and algorithmic consequences.
\newblock {\em Random Struct. Algorithms}, 56(3):597--647, 2020.

\bibitem[Haj82]{hajek_1982}
Bruce Hajek.
\newblock Hitting-time and occupation-time bounds implied by drift analysis
  with applications.
\newblock {\em Advances in Applied Probability}, 14(3):502--525, 1982.

\bibitem[Her18]{DBLP:journals/ec/Heredia18}
Jorge~P{\'{e}}rez Heredia.
\newblock Modelling evolutionary algorithms with stochastic differential
  equations.
\newblock {\em Evol. Comput.}, 26(4), 2018.

\bibitem[HJZ19]{DBLP:conf/gecco/00040Z19}
Jun He, Thomas Jansen, and Christine Zarges.
\newblock Unlimited budget analysis.
\newblock In Manuel L{\'{o}}pez{-}Ib{\'{a}}{\~{n}}ez, Anne Auger, and Thomas
  St{\"{u}}tzle, editors, {\em Proceedings of the Genetic and Evolutionary
  Computation Conference Companion, {GECCO} 2019, Prague, Czech Republic, July
  13-17, 2019}, pages 427--428. {ACM}, 2019.

\bibitem[HOS19]{DBLP:conf/gecco/HallOS19}
George~T. Hall, Pietro~S. Oliveto, and Dirk Sudholt.
\newblock On the impact of the cutoff time on the performance of algorithm
  configurators.
\newblock In Anne Auger and Thomas St{\"{u}}tzle, editors, {\em Proceedings of
  the Genetic and Evolutionary Computation Conference, {GECCO} 2019, Prague,
  Czech Republic, July 13-17, 2019}, pages 907--915. {ACM}, 2019.

\bibitem[HY01]{DBLP:journals/ai/HeY01}
Jun He and Xin Yao.
\newblock Drift analysis and average time complexity of evolutionary
  algorithms.
\newblock {\em Artif. Intell.}, 127(1):57--85, 2001.

\bibitem[HY04]{DBLP:journals/nc/HeY04}
Jun He and Xin Yao.
\newblock A study of drift analysis for estimating computation time of
  evolutionary algorithms.
\newblock {\em Nat. Comput.}, 3(1):21--35, 2004.

\bibitem[Jan20]{DBLP:series/ncs/000120}
Thomas Jansen.
\newblock Analysing stochastic search heuristics operating on a fixed budget.
\newblock In Benjamin Doerr and Frank Neumann, editors, {\em Theory of
  Evolutionary Computation - Recent Developments in Discrete Optimization},
  Natural Computing Series, pages 249--270. Springer, 2020.

\bibitem[JL22]{DBLP:conf/ppsn/JanettL22}
Duri Janett and Johannes Lengler.
\newblock Two-dimensional drift analysis: - optimizing two functions
  simultaneously can be hard.
\newblock In G{\"{u}}nter Rudolph, Anna~V. Kononova, Hern{\'{a}}n~E. Aguirre,
  Pascal Kerschke, Gabriela Ochoa, and Tea Tusar, editors, {\em Parallel
  Problem Solving from Nature - {PPSN} {XVII} - 17th International Conference,
  {PPSN} 2022, Dortmund, Germany, September 10-14, 2022, Proceedings, Part
  {II}}, volume 13399 of {\em Lecture Notes in Computer Science}, pages
  612--625. Springer, 2022.

\bibitem[Joh10]{Johannsen:thesis:10}
Daniel Johannsen.
\newblock {\em Random Combinatorial Structures and Randomized Search
  Heuristics}.
\newblock PhD thesis, Universität des Saarlandes, 2010.

\bibitem[JZ12]{DBLP:conf/gecco/JansenZ12}
Thomas Jansen and Christine Zarges.
\newblock Fixed budget computations: a different perspective on run time
  analysis.
\newblock In Terence Soule and Jason~H. Moore, editors, {\em Genetic and
  Evolutionary Computation Conference, {GECCO} '12, Philadelphia, PA, USA, July
  7-11, 2012}, pages 1325--1332. {ACM}, 2012.

\bibitem[JZ14]{DBLP:journals/tec/0001Z14}
Thomas Jansen and Christine Zarges.
\newblock Reevaluating immune-inspired hypermutations using the fixed budget
  perspective.
\newblock {\em {IEEE} Trans. Evol. Comput.}, 18(5):674--688, 2014.

\bibitem[KK18]{DBLP:conf/ppsn/KotzingK18}
Timo K{\"{o}}tzing and Martin~S. Krejca.
\newblock First-hitting times for finite state spaces.
\newblock In Anne Auger, Carlos~M. Fonseca, Nuno Louren{\c{c}}o, Penousal
  Machado, Lu{\'{\i}}s Paquete, and L.~Darrell Whitley, editors, {\em Parallel
  Problem Solving from Nature - {PPSN} {XV} - 15th International Conference,
  Coimbra, Portugal, September 8-12, 2018, Proceedings, Part {II}}, volume
  11102 of {\em Lecture Notes in Computer Science}, pages 79--91. Springer,
  2018.

\bibitem[KK19]{DBLP:journals/tcs/KotzingK19}
Timo K{\"{o}}tzing and Martin~S. Krejca.
\newblock First-hitting times under drift.
\newblock {\em Theor. Comput. Sci.}, 796:51--69, 2019.

\bibitem[KLW15]{DBLP:conf/foga/KotzingLW15}
Timo K{\"{o}}tzing, Andrei Lissovoi, and Carsten Witt.
\newblock {(1+1)} {EA} on generalized dynamic onemax.
\newblock In Jun He, Thomas Jansen, Gabriela Ochoa, and Christine Zarges,
  editors, {\em Proceedings of the 2015 {ACM} Conference on Foundations of
  Genetic Algorithms XIII, Aberystwyth, United Kingdom, January 17 - 20, 2015},
  pages 40--51. {ACM}, 2015.

\bibitem[KM12]{DBLP:conf/ppsn/KotzingM12}
Timo K{\"{o}}tzing and Hendrik Molter.
\newblock {ACO} beats {EA} on a dynamic pseudo-boolean function.
\newblock In Carlos A.~Coello Coello, Vincenzo Cutello, Kalyanmoy Deb,
  Stephanie Forrest, Giuseppe Nicosia, and Mario Pavone, editors, {\em Parallel
  Problem Solving from Nature - {PPSN} {XII} - 12th International Conference,
  Taormina, Italy, September 1-5, 2012, Proceedings, Part {I}}, volume 7491 of
  {\em Lecture Notes in Computer Science}, pages 113--122. Springer, 2012.

\bibitem[K{\"{o}}t16]{DBLP:journals/algorithmica/Kotzing16}
Timo K{\"{o}}tzing.
\newblock Concentration of first hitting times under additive drift.
\newblock {\em Algorithmica}, 75(3):490--506, 2016.

\bibitem[Kre19]{Krejca:thesis:19}
Martin~S. Krejca.
\newblock {\em Theoretical analyses of univariate estimation-of-distribution
  algorithms}.
\newblock PhD thesis, Universität Potsdam, 2019.

\bibitem[KSNO12]{DBLP:conf/gecco/KotzingSNO12}
Timo K{\"{o}}tzing, Andrew~M. Sutton, Frank Neumann, and Una{-}May O'Reilly.
\newblock The max problem revisited: the importance of mutation in genetic
  programming.
\newblock In Terence Soule and Jason~H. Moore, editors, {\em Genetic and
  Evolutionary Computation Conference, {GECCO} '12, Philadelphia, PA, USA, July
  7-11, 2012}, pages 1333--1340. {ACM}, 2012.

\bibitem[KST11]{DBLP:conf/gecco/KotzingST11}
Timo K{\"{o}}tzing, Dirk Sudholt, and Madeleine Theile.
\newblock How crossover helps in pseudo-boolean optimization.
\newblock In Natalio Krasnogor and Pier~Luca Lanzi, editors, {\em 13th Annual
  Genetic and Evolutionary Computation Conference, {GECCO} 2011, Proceedings,
  Dublin, Ireland, July 12-16, 2011}, pages 989--996. {ACM}, 2011.

\bibitem[KU18]{DBLP:conf/podc/KosowskiU18}
Adrian Kosowski and Przemyslaw Uznanski.
\newblock Brief announcement: Population protocols are fast.
\newblock In Calvin Newport and Idit Keidar, editors, {\em Proceedings of the
  2018 {ACM} Symposium on Principles of Distributed Computing, {PODC} 2018,
  Egham, United Kingdom, July 23-27, 2018}, pages 475--477. {ACM}, 2018.

\bibitem[KW20]{DBLP:conf/ppsn/KotzingW20}
Timo K{\"{o}}tzing and Carsten Witt.
\newblock Improved fixed-budget results via drift analysis.
\newblock In Thomas B{\"{a}}ck, Mike Preuss, Andr{\'{e}}~H. Deutz, Hao Wang,
  Carola Doerr, Michael T.~M. Emmerich, and Heike Trautmann, editors, {\em
  Parallel Problem Solving from Nature - {PPSN} {XVI} - 16th International
  Conference, {PPSN} 2020, Leiden, The Netherlands, September 5-9, 2020,
  Proceedings, Part {II}}, volume 12270 of {\em Lecture Notes in Computer
  Science}, pages 648--660. Springer, 2020.

\bibitem[Leh11]{DBLP:conf/gecco/Lehre11}
Per~Kristian Lehre.
\newblock Fitness-levels for non-elitist populations.
\newblock In Natalio Krasnogor and Pier~Luca Lanzi, editors, {\em 13th Annual
  Genetic and Evolutionary Computation Conference, {GECCO} 2011, Proceedings,
  Dublin, Ireland, July 12-16, 2011}, pages 2075--2082. {ACM}, 2011.

\bibitem[Len20]{DBLP:series/ncs/Lengler20}
Johannes Lengler.
\newblock Drift analysis.
\newblock In Benjamin Doerr and Frank Neumann, editors, {\em Theory of
  Evolutionary Computation - Recent Developments in Discrete Optimization},
  Natural Computing Series, pages 89--131. Springer, 2020.

\bibitem[LS14]{DBLP:journals/ec/LassigS14}
J{\"{o}}rg L{\"{a}}ssig and Dirk Sudholt.
\newblock General upper bounds on the runtime of parallel evolutionary
  algorithms.
\newblock {\em Evol. Comput.}, 22(3):405--437, 2014.

\bibitem[LS15]{DBLP:conf/foga/LenglerS15}
Johannes Lengler and Nicholas Spooner.
\newblock Fixed budget performance of the {(1+1)} {EA} on linear functions.
\newblock In Jun He, Thomas Jansen, Gabriela Ochoa, and Christine Zarges,
  editors, {\em Proceedings of the 2015 {ACM} Conference on Foundations of
  Genetic Algorithms XIII, Aberystwyth, United Kingdom, January 17 - 20, 2015},
  pages 52--61. {ACM}, 2015.

\bibitem[LW14]{DBLP:conf/isaac/LehreW14}
Per~Kristian Lehre and Carsten Witt.
\newblock Concentrated hitting times of randomized search heuristics with
  variable drift.
\newblock In Hee{-}Kap Ahn and Chan{-}Su Shin, editors, {\em Algorithms and
  Computation - 25th International Symposium, {ISAAC} 2014, Jeonju, Korea,
  December 15-17, 2014, Proceedings}, volume 8889 of {\em Lecture Notes in
  Computer Science}, pages 686--697. Springer, 2014.

\bibitem[McD93]{mcdiarmid_1993}
Colin McDiarmid.
\newblock A random recolouring method for graphs and hypergraphs.
\newblock {\em Combinatorics, Probability and Computing}, 2(3):363--365, 1993.

\bibitem[MRC09]{DBLP:journals/ijicc/MitavskiyRC09}
Boris Mitavskiy, Jonathan~E. Rowe, and Chris Cannings.
\newblock Theoretical analysis of local search strategies to optimize network
  communication subject to preserving the total number of links.
\newblock {\em Int. J. Intell. Comput. Cybern.}, 2(2):243--284, 2009.

\bibitem[MU05]{MitzenmacherUpfal:2005:ProbComp}
Michael Mitzenmacher and Eli Upfal.
\newblock {\em Probability and computing: randomized algorithms and
  probabilistic analysis}.
\newblock Cambridge University Press, 2005.

\bibitem[NNS17]{DBLP:journals/ec/NallaperumaNS17}
Samadhi Nallaperuma, Frank Neumann, and Dirk Sudholt.
\newblock Expected fitness gains of randomized search heuristics for the
  traveling salesperson problem.
\newblock {\em Evol. Comput.}, 25(4), 2017.

\bibitem[OW11]{DBLP:journals/algorithmica/OlivetoW11}
Pietro~S. Oliveto and Carsten Witt.
\newblock Simplified drift analysis for proving lower bounds in~evolutionary
  computation.
\newblock {\em Algorithmica}, 59(3):369--386, 2011.

\bibitem[OW12]{DBLP:journals/corr/OlivetoW12}
Pietro~S. Oliveto and Carsten Witt.
\newblock Erratum: Simplified drift analysis for proving lower bounds in
  evolutionary computation.
\newblock {\em CoRR}, abs/1211.7184, 2012.

\bibitem[OW14]{DBLP:journals/tcs/OlivetoW14}
Pietro~S. Oliveto and Carsten Witt.
\newblock On the runtime analysis of the simple genetic algorithm.
\newblock {\em Theor. Comput. Sci.}, 545:2--19, 2014.

\bibitem[Pap91]{DBLP:conf/focs/Papadimitriou91}
Christos~H. Papadimitriou.
\newblock On selecting a satisfying truth assignment (extended abstract).
\newblock In {\em 32nd Annual Symposium on Foundations of Computer Science, San
  Juan, Puerto Rico, 1-4 October 1991}, pages 163--169. {IEEE} Computer
  Society, 1991.

\bibitem[Row18]{DBLP:journals/tcs/Rowe18}
Jonathan~E. Rowe.
\newblock Linear multi-objective drift analysis.
\newblock {\em Theor. Comput. Sci.}, 736:25--40, 2018.

\bibitem[RS14]{DBLP:journals/tcs/RoweS14}
Jonathan~E. Rowe and Dirk Sudholt.
\newblock The choice of the offspring population size in the (1, {\(\lambda\)})
  evolutionary algorithm.
\newblock {\em Theor. Comput. Sci.}, 545:20--38, 2014.

\bibitem[Sem03]{DBLP:journals/ec/Semenov03}
Mikhail~A. Semenov.
\newblock Analysis of convergence of an evolutionary algorithm with
  self-adaptation using a stochastic lyapunov function.
\newblock {\em Evol. Comput.}, 11(4):363--379, 2003.

\bibitem[STW04]{DBLP:journals/jmma/ScharnowTW04}
Jens Scharnow, Karsten Tinnefeld, and Ingo Wegener.
\newblock The analysis of evolutionary algorithms on sorting and shortest paths
  problems.
\newblock {\em J. Math. Model. Algorithms}, 3(4):349--366, 2004.

\bibitem[Sud13]{DBLP:journals/tec/Sudholt13}
Dirk Sudholt.
\newblock A new method for lower bounds on the running time of evolutionary
  algorithms.
\newblock {\em {IEEE} Trans. Evol. Comput.}, 17(3):418--435, 2013.

\bibitem[Weg01]{DBLP:conf/icalp/Wegener01}
Ingo Wegener.
\newblock Theoretical aspects of evolutionary algorithms.
\newblock In Fernando Orejas, Paul~G. Spirakis, and Jan van Leeuwen, editors,
  {\em Automata, Languages and Programming, 28th International Colloquium,
  {ICALP} 2001, Crete, Greece, July 8-12, 2001, Proceedings}, volume 2076 of
  {\em Lecture Notes in Computer Science}, pages 64--78. Springer, 2001.

\bibitem[Wit06]{DBLP:journals/ec/Witt06}
Carsten Witt.
\newblock Runtime analysis of the (mu + 1) {EA} on simple pseudo-boolean
  functions.
\newblock {\em Evol. Comput.}, 14(1):65--86, 2006.

\bibitem[Wit13]{DBLP:journals/cpc/Witt13}
Carsten Witt.
\newblock Tight bounds on the optimization time of a randomized search
  heuristic on linear functions.
\newblock {\em Comb. Probab. Comput.}, 22(2):294--318, 2013.

\bibitem[Wit14]{DBLP:journals/ipl/Witt14}
Carsten Witt.
\newblock Fitness levels with tail bounds for the analysis of randomized search
  heuristics.
\newblock {\em Inf. Process. Lett.}, 114(1-2):38--41, 2014.

\bibitem[Wit23]{DBLP:journals/tcs/Witt23}
Carsten Witt.
\newblock How majority-vote crossover and estimation-of-distribution algorithms
  cope with fitness valleys.
\newblock {\em Theor. Comput. Sci.}, 940(Part):18--42, 2023.

\bibitem[Wor99]{Wormald:j:99}
Nicholas~C. Wormald.
\newblock The differential equation method for random graph processes and
  greedy algorithms.
\newblock {\em Lectures on Approximation and Randomized Algorithms}, pages
  73--155, 1999.

\end{thebibliography}

\end{document}